\documentclass[reqno,11pt]{amsart}

\usepackage{amssymb}
\usepackage{amsmath}
\usepackage[all]{xy}
\usepackage{caption}
\usepackage[pdftex]{hyperref}
\usepackage{graphicx}
\usepackage{array}
\usepackage{physics}
\usepackage{color}
\usepackage{xy}
\usepackage{bm}
\usepackage{comment}
\usepackage{mathtools}
\usepackage{tikz-cd}
\usetikzlibrary{arrows.meta,decorations.markings}
\usepackage{float}
\usetikzlibrary{knots, hobby}

\allowdisplaybreaks

\newtheorem{thm}{Theorem}[section]
\newtheorem{lem}[thm]{Lemma}
\newtheorem{prop}[thm]{Proposition}
\newtheorem{cor}[thm]{Corollary}
\theoremstyle{definition}
\newtheorem{defn}[thm]{Definition}
\newtheorem{ex}[thm]{Example}
\newtheorem{rem}[thm]{Remark}

\newtheorem{introthm}{Theorem}

\providecommand\sslash{\mathbin{/\mkern-5.5mu/}}

\newcommand{\sldos}{\SL_2(\mathbb{C})}

\newcommand{\CC}{{\mathbb C}}
\newcommand{\ZZ}{{\mathbb Z}}
\newcommand{\RR}{{\mathbb R}}

\newcommand{\R}{{\mathcal{R}}}
\newcommand{\X}{{\mathcal{X}}}

\newcommand{\cA}{{\mathcal{A}}}
\newcommand{\cR}{{\mathcal{R}}}

\newcommand{\Free}[1]{\mathbf{F}_{#1}}
\newcommand{\SVB}{\textbf{\textup{SVectB}}}

\newcommand{\Varc}{\mathbf{Var}_{\mathbb{C}}}

\newcommand{\Arr}[1]{#1^{\to}}

\newcommand{\Qdl}{\textbf{\textup{Qdl}}}
\newcommand{\Grp}{\textbf{\textup{Grp}}}

\newcommand{\GL}[1]{{\textup{GL}_{#1}}}
\newcommand{\AGL}[1]{{\textup{AGL}_{#1}}}

\newcommand{\aglC}{\AGL{1}(\mathbb{C})}

\DeclareMathOperator{\Spec}{Spec}
\DeclareMathOperator{\id}{Id}
\DeclareMathOperator{\Id}{Id}
\DeclareMathOperator{\Hom}{Hom\,}           %
\DeclareMathOperator{\Stab}{Stab}

\DeclareMathOperator{\SL}{SL}

\DeclareMathOperator{\diag}{diag}
\DeclareMathOperator{\Conj}{Conj}
\DeclareMathOperator{\Span}{Span}
\DeclareMathOperator{\Vect}{\mathbf{Vect}}
\DeclareMathOperator{\As}{As}
\newcommand{\Tang}{\mathsf{Tang}}

\newcommand{\upcurvearrowright}{\rotatebox[origin=c]{180}{$\curvearrowleft$}}
\newcommand{\upcurvearrowleft}{\rotatebox[origin=c]{180}{$\curvearrowright$}}

\title[Parabolic character varieties of knots via quandle representations]{Parabolic character varieties of knots\\ via quandle representations}

\author{\'Angel Gonz\'alez-Prieto}
\address{Departamento de \'Algebra, Geometr\'ia y Topolog\'ia, Facultad de Ciencias Matem\'aticas, Universidad Complutense de Madrid, Plaza Ciencias 3, 28040 Madrid, Spain; and Instituto de Ciencias Matem\'aticas (CSIC-UAM-UCM-UC3M), C.\ Nicolás Cabrera 13-15, 28049 Madrid, Spain.}\email{angelgonzalezprieto@ucm.es}

\author{Javier Mart\'inez}
\address{Departamento de Matem\'atica Aplicada a la Ingeniería Industrial, Escuela Técnica Superior de Ingeniería y Diseño Industrial, Universidad Politécnica de Madrid, Ronda de Valencia 3, 28012 Madrid, Spain.}
\email{javier.martinezmartinez@upm.es}

\begin{document}

\begin{abstract}
    In this paper, we study parabolic character varieties of links by means of representation varieties of link quandles. For this purpose, we construct a functor quantizing quandle representation varieties for any complex algebraic quandle, as well as its relative version. Applying this construction to the morphism arising from the Chevalley map of a complex reductive group, we are able to recover parabolic representation varieties as slice varieties for $G = \textup{SL}_2$. Furthermore, building on the description introduced by Jo and Kim of the quandle of parabolic matrices of $\textup{SL}_2$ as the parabolic symplectic quandle, we reformulate the Riley polynomial of a $2$-bridge knot in terms of the generator of the ideal of the coordinate ring of the parabolic character variety. As application, we provide explicit descriptions of these parabolic representation and character varieties for arbitrary pretzel knots.
\end{abstract}

\maketitle
\section{Introduction}\label{sec:introduction}

Let $K \subseteq S^3$ be a knot and let $G$ be a complex reductive group. The $G$-representation variety of $K$ is the space
 $$
    \R_G(K) = \Hom(\Gamma_K,G)
 $$
 of representations of the knot group $\Gamma_K=\pi_1(S^3-K)$, the fundamental group of the complement of $K$, into $G$, which is naturally endowed with the structure of a complex algebraic variety. Furthermore, the representation variety $\R_G(K)$ inherits an action of $G$ by conjugation, and the corresponding Geometric Invariant Theory (GIT) quotient is the so-called character variety
 $$
    \X_G(K) = \R_G(K) \sslash G.
 $$

 Since the foundational works in the 1980s and 1990s \cite{CullerShalen1983, MorganShalen1984, Cooper1994Planecurves}, character varieties were established as a powerful tool for exploring the geometry and the topology of $3$-manifolds. They have played a central role in knot theory, due to their ubiquity and through their connections with essential surfaces and Dehn surgery, skein theory and quantization \cite{Bullock1997,PrzSikora2000} or hyperbolic geometry.
 
 An interesting case arises for $G=\sldos$, where we restrict to the so-called parabolic representations, i.e.\ representations for which each meridian is sent to a parabolic element. The pioneering work of Riley \cite{RileyparabolicI,RileyparabolicII,RileyparabolicIII} provided an effective algebraic description of such representations for $2$-bridge knots. In this setting, Riley showed that the representation equations could be reduced to a single polynomial equation, nowadays known as the {Riley polynomial} $R_K(u) \in \ZZ[u]$, whose roots encode the non-abelian parabolic representations up to conjugacy. Riley's approach was also fundamental in some of the first explicit constructions of hyperbolic structures on knot complements. Riley's work also influences more recent developments, such as the study of twisted Alexander polynomials associated with parabolic representations of two-bridge knots \cite{MorifujiTran}.

The starting point of this manuscript is the work of Jo and Kim \cite{symplecticquandle1}, in which parabolic representations were related to the symplectic quandle framework, building on the parabolic quandle approach developed in \cite{InoueKabaya2014,Cho2018}. Their method led to explicit computations for $2$-bridge knots and links, as well as further arithmetic and algebraic properties of the Riley polynomial, including the factorization
$$
R_K(-u^2)=\pm g(u)g(-u),
\qquad \textrm{for some } g(u)\in\mathbb Z[u],
$$
which was further investigated in \cite{JoKimContinuantChebyshevRiley}. This approach was later extended to general knots with zero-dimensional parabolic character variety by Cho, Kim, Kim and Yoon \cite{symplecticquandle3}. For a knot diagram with a specified crossing $c$, they introduced a generalized Riley polynomial $R_c(u)\in\mathbb Q[u]$, whose roots correspond to conjugacy classes of parabolic representations, in the spirit of Riley's original construction. Finally, Jo and Kim have recently extended the symplectic quandle method beyond the parabolic setting to general $\sldos$-representations of $2$-bridge ``kmot'' groups, a class of two-generator one-relator groups modelled on the standard presentations of $2$-bridge knot groups \cite{symplecticquandle2}.

Regarding representations of knot groups, the present authors, together with Mu\~{n}oz, developed in \cite{AGL1Alexanderinvariants} the functorial constructions
$$
\cA:\Tang \longrightarrow \Span(\SVB),
\qquad
\cA_t:\Tang \longrightarrow \Span(\Vect_{\CC}),
$$
which encode the $\aglC$-representation varieties of tangles and knots and their associated Alexander invariants from an algebro-geometric perspective using categories of spans. Here, $\Tang$ is the category of oriented tangles and $\SVB$ is the category of singular vector bundles. This provides, in particular, a geometric extension of the classical Burde-de Rham theorem, which characterizes the roots of the Alexander polynomial in terms of non-abelian affine representations of the knot group. More precisely, the Alexander module of $K$ is recovered directly from the coordinate ring of the representation variety
so that the Alexander polynomial is encoded in the algebraic geometry of the $\mathrm{AGL}_1(\mathbb C)$-representation space. This functorial approach also allows explicit computations and has subsequently led to closed formulas for the Alexander polynomials of rational and pretzel knots \cite{Martinezaffine}.

An important ingredient in \cite{AGL1Alexanderinvariants} and \cite{Martinezaffine} is the adjunction between groups and quandles \cite{racks}, which allows to interpret the classical representation varieties as quandle representation varieties into the conjugacy quandle $\Conj \aglC$ and the Alexander quandle. Concretely, if we orient $K$, then the natural projection $\aglC \to \CC^*$ induces a natural fibration
\begin{equation}\label{eq:Alexander-morphism}
    \R_{\aglC}(K) \to \R_{\CC^*}(K) = \CC^*.
\end{equation}
In \cite{AGL1Alexanderinvariants}, it is proven that the fiber of this map at $t \in \CC^*$ is precisely the quandle representation variety $\R_{Q_t}(K) = \Hom_{\Qdl}(Q_K, Q_t)$, where $Q_K$ is the so-called knot quandle, as introduced in \cite{Joycequandle}, and $Q_t$ is the $t$-Alexander quandle. 

In this manuscript, quandle representation varieties are treated as first-class citizens, generalizing the previous framework to more general quandles. To be precise, given a complex algebraic quandle $X$, the $X$-quandle representation variety is the complex algebraic variety
$$
    \R_X(K) = \Hom_{\Qdl}(Q_K, X).
$$
These quandle representation varieties generalize standard representation varieties of groups in the case $X = \Conj G$ for a group $G$, since in this case we have a natural identification $\R_{\Conj G}(K) \cong \R_{G}(K)$.

Going a step further, we will be interested in keeping track of the regular morphism $\R_X(K) \to \R_{X'}(K)$ induced by a morphism $\alpha: X \to X'$. For instance, for $X = \Conj \aglC$ and $X' = \CC^*$ with the trivial quandle structure, the natural quandle homomorphism $\Conj \aglC \to \CC^*$ induces a regular morphism $\R_{\Conj \aglC}(K) = \R_{\aglC}(K) \to \R_{\CC^*}(K)$, recovering (\ref{eq:Alexander-morphism}). 

In this direction, generalizing the previous construction, we prove the following result.

\begin{introthm}[{Theorem \ref{thm:Zalphafunctor}}]
Let $\alpha\colon X\rightarrow X'$ be a morphism of complex algebraic quandles. There exists a monoidal functor
$$
Z_\alpha \colon \Tang \longrightarrow \Span(\Arr{\Varc})
$$
which computes the quandle representation variety morphism $\R_X(T) \to \R_{X'}(T)$, for any tangle $T\in \Tang$.
\end{introthm}

To be precise, the functor above assigns to any object $\epsilon \in \Tang$, which is a finite set of oriented points, the morphism of representation varieties $\alpha_{\#}: \R_X(\epsilon) \to \R_{X'}(\epsilon)$, and to any tangle $T\colon\epsilon_1\rightarrow \epsilon_2$ between $\epsilon_1$ and $\epsilon_2$, the double span 
\begin{equation} \label{eqn:Zalphaarrowspan}
\begin{split}
\xymatrix{
& \R_{X}(T) \ar[d] \ar[dl]_{i_1} \ar[dr]^{i_2} & \\
\R_X(\epsilon_1) \ar[d] & \R_{X'}(T) \ar[dl]_{i_1'} \ar[dr]^{i_2'} & \R_X(\epsilon_2) \ar[d] \\
\R_{X'}(\epsilon_1)  & & \R_{X'}(\epsilon_2) 
}
\end{split}
\end{equation}
which is a span in the arrow category $\Arr{\Varc}$ of $\Varc$. It is worth noticing that this construction is a representation of the bordered quandle category in the sense of \cite{cattabriga}.

A particular case of interest arises when the morphism of quandles is the Chevalley map $\chi:G\rightarrow C_G:=G\sslash G$, interpreted as a quandle morphism from the conjugation quandle $\Conj G$ to $G\sslash G$, endowed with the trivial quandle structure. The Chevalley base $\cR_{C_G}(T)$ keeps track of the spectral meridional local data of each representation in \eqref{eqn:Zalphaarrowspan}. In particular, the fibers $G_t = \chi^{-1}(t)$ of the Chevalley map, where $t\in C_G$, inherit a quandle structure, so that $Z_{G_t}$ can be seen as the fibers of the functor $Z_{\chi}$ at $t\in C_G$, as it occurred in the Alexander case where $C_{\aglC} = \CC^*$.

When $G=\SL_2(\CC)$, we have that the Chevalley map $\chi: \SL_2(\CC) \to C_{\SL_2(\CC)}=\CC$ is the trace map, and the quandle structure of $(\SL_2(\CC))_t = \chi^{-1}(t) \subseteq \Conj \SL_2(\CC)$, for $t \in \CC$, can be expressed in terms of generalized symplectic quandles, as shown in \cite{symplecticquandle1,symplecticquandle2}. Of particular interest is the case $t=2$, since $(\SL_2(\CC))_2$ is precisely the set of parabolic matrices of $\SL_2(\CC)$, which turn out to be isomorphic as quandle to the symplectic quandle $Q_2:=\CC^2/\ZZ_2$, where $x\mapsto -x$, with quandle operations for $x,y \in \CC^2$
\begin{equation*}
[x] \triangleright [y] = \left[y+\langle x,y \rangle x \right], \qquad [y] \triangleleft [x] = \left[y - \langle x,y \rangle x\right].
\end{equation*}
Here $\langle x,y\rangle = \det(x,y)$ denotes the standard complex symplectic form in $\CC^2$.
In this way, the quandle representation variety $\R_{Q_2}(K) = \R_{(\SL_2(\CC))_2}(K)$ is precisely the above-mentioned parabolic representation variety of $K$.

Following these ideas, after choosing generators of $\Tang$, we provide a complete description of the functor $Z_{\chi}$ for $Q_2$ and its double cover $\widetilde{Q}_2:=\CC^2$. Thanks to this, for any knot $K$, we explicitly describe the rings of regular functions of the representation varieties $\R_{Q_2}(K)$ and $\R_{\widetilde{Q}_2}(K)$ and their associated character varieties $ \X_{Q_2}(K)$ and $\X_{\widetilde{Q}_2}(K)$. For this purpose, we shall describe them as closed subschemes of $\R_{Q_2}(\Free{n})$ (resp.\ $\R_{\widetilde{Q}_2}(\Free{n})$) for the representation variety, and $\X_{Q_2}(\Free{n})$ (resp.\ $\X_{\widetilde{Q}_2}(\Free{n})$) for the character variety, where $\Free{n}$ is the free quandle in $n$ generators. 
It is worth mentioning that, during preparation of this manuscript, the scheme structure of the parabolic character variety of $\X_{Q_2}(K)$ was independently described in \cite{ChoKimYoonschemes} using sign-refined arc-colorings.

In this context, the algebra of the parabolic character variety $\CC[\X_{Q_2}(K)]$ arises as a knot invariant, which we refer to as the {Riley algebra}. We also show that it provides the right context to understand the Riley polynomial: For rational knots, $\CC[\X_{Q_2}(K)]$ is a $0$-dimensional algebra of the form $\CC[\X_{Q_2}(K)] \cong \CC[u]/(uR_K(u))$, where $R_K(u)$ is the classical Riley polynomial, obtained after a choice of presentation of $Q_K$. Notice that, in order to recover the Riley polynomial, we must remove the point $u = 0$, which corresponds to the trivial character not considered in the work of Riley. However, although the Riley algebra $\CC[\X_{Q_2}(K)]$ is a knot invariant, the Riley polynomial $R_K(u)$ is not, as it corresponds to the well-known fact that different presentations of the same rational knot may lead to different Riley polynomials (see \cite{KitanoMorifuji}). 

\begin{figure}[h]
  \centering
  \includegraphics[width=0.28\textwidth]{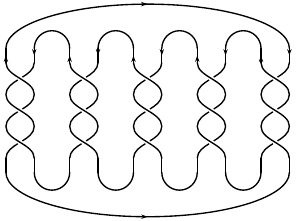}\hfill
  \caption{Pretzel knot.}
  \label{fig:pretzelknot-ex}
\end{figure}

Finally, we apply all these techniques to a particular family of knots, the so-called pretzel knots. Recall that, given integers $q_1, \ldots, q_n \in \ZZ$, the pretzel knot $P(q_1,\ldots,q_n) \subseteq S^3$ is a knot constructed by considering $n$ braidings of two strands, where the $i$-th pair is twisted $q_i$ times, and then closed sequentially as in Figure \ref{fig:pretzelknot-ex}.

For this class of knots, we prove the following theorem.

\begin{introthm}[{Theorem \ref{thm:pullback-rep-var}}] \label{introthm2}
Let $P(q_1,\ldots,q_n) \subseteq S^3$ be an oriented pretzel knot. Then, there exist morphisms $f_{1}, \ldots, f_{n}: \widetilde{Q}_2 \times \widetilde{Q}_2 \to \widetilde{Q}_2 \times \widetilde{Q}_2$ and $\iota\colon \widetilde{Q}_2^n\longrightarrow \left(\widetilde{Q}_2 \times \widetilde{Q}_2\right)^n$ such that the $\widetilde{Q}_2$-representation variety of $P(q_1,\ldots,q_n)$ is given as the pullback
$$
\xymatrix{
\R_{\widetilde{Q}_2}(P(q_1,\ldots,q_n)) \ar[r]\ar[d] & \widetilde{Q}_2^n\ar[d]^{(f_{1} \times \ldots \times f_{n}) \circ \iota} \\
\widetilde{Q}_2^n \ar[r]_{\iota} & \left(\widetilde{Q}_2 \times \widetilde{Q}_2\right)^n
}
$$
\end{introthm}

It is important to point out that the morphisms $f_1, \ldots, f_n$ and $\iota$ appearing in the previous result are completely explicit, and can be described in terms of powers of matrices with coefficients on a polynomial ring, see Sections \ref{sec:codirected-braid} and \ref{sec:alternating-braid}, and particularly Table \ref{tab:summary-tangles}.

Furthermore, using the quandle covering $\tilde{Q}_2 = \CC^2 \rightarrow Q_2 = \CC^2/\ZZ_2$, Theorem \ref{introthm2} allows us to describe the classical parabolic representation and character varieties of any pretzel knot. For this purpose, using the explicit matrix maps of Theorem \ref{introthm2}, it is straightforward to construct a vector $(r_1, \ldots, r_{2n})$ whose entries are elements $r_k \in \CC^2$ that can be written as can be
written as vector-valued polynomials in the variables $a_1 = (a_1^1, a_1^2), \ldots, a_n = (a_n^1, a_n^2) \in \CC^2$.

\begin{introthm}[{Proposition \ref{prop:Q2pretzelrepvariety}}]
For any pretzel knot $P(q_1, \ldots, q_n) \subseteq S^3$, the parabolic representation variety is given by
\begin{align*}
    \R_{Q_2}(P(q_1, \ldots, q_n)) & = \left\{([a_1], \ldots, [a_n]) \in (\CC^2/\ZZ_2)^n\,\mid\, r_pr_p^t = r_qr_q^t \textrm{ for } (p,q) \in \Delta\right\} \\
    & = \Spec \frac{\CC[x_k, y_k, z_k]}{(r_pr_p^t-r_qr_q^t, x_kz_k - y_k^2)},
\end{align*}
for $1 \leq k \leq n$ and $(p,q) \in \Delta = \{(2,3),(4,5), \ldots, (2n-2, 2n-1), (1, 2n)\}$.
\end{introthm}

In the last description of $\R_{Q_2}(P(q_1, \ldots, q_n))$, the new variables are given by $x_k = (a_k^1)^2, y_k = a_k^1a_k^2$ and $z_k = (a_k^2)^2$ for $1 \leq k \leq n$, for which we prove that the entries of the products $r_pr_p^t-r_qr_q^t$ can be easily written as a polynomial in these variables. 

Furthermore, this description is so explicit that it also allows us to describe the corresponding parabolic character variety. For this purpose, we introduce the new variables
$$
    f_{pk} = \langle r_p, a_k \rangle, \qquad s_{k\ell} = \langle a_k, a_\ell\rangle,
$$
for $1 \leq p \leq 2n$ and $1 \leq k, \ell \leq n$. Now, consider the new variables $u_{ij} = s_{ij}^2$ and $t_{ijk} = s_{ij}s_{jk}s_{ik}$ for $1 \leq i < j < k \leq n$. Moreover, for $\gamma = (i_1, i_2, \ldots, i_s)$, we set $s_\gamma = s_{i_1i_2} s_{i_2i_3} \cdots s_{i_{s-1}i_s}$.

\begin{introthm}[{Corollary \ref{cor:char-var-Q2-pretzel}}]
 For any pretzel knot $P(q_1, \ldots, q_n) \subseteq S^3$,
    $$
        \mathcal{X}_{{Q}_2}(P(q_1, \ldots, q_n)) = \Spec \frac{\CC[u_{ij}, t_{ijk}]}{H_n + \left(s_{\gamma}(f_{pk}f_{p\ell} - f_{qk}f_{q\ell}), f_{pk}^2 - f_{qk}^2\right)},
    $$
    for $1 \leq i < j < k < \ell \leq n, (p,q) \in \Delta$ and $\gamma$ a simple path between $k$ and $\ell$. Here $H_n \subseteq \CC[u_{ij}, t_{ijk}]$ is an ideal that only depends on $n$.
\end{introthm}

The exact expressions for $f_{pk}, r_i$ and $s_\gamma$ are directly derived from Theorem \ref{introthm2} and can be checked in Sections \ref{sec:Q2tilde-pretzel} and \ref{sec:Q2-pretzel}. Additionally, a SageMath code is provided to automatically conduct these calculations for any pretzel knot, as publicly available in \cite{code_repository}.

As a final application of these techniques, we compute the representation and character variety of a non-rational example, $P(-2,-3,1,-3)$, and provide a pretzel family whose character variety has arbitrarily large dimension, matching \cite{PortiPaoluzzi2013}.

\subsection*{Structure of the manuscript}
After introducing in Section \ref{sec:quandlerepvarieties} quandle representation varieties and the conventions used throughout the paper, the construction of the functor $Z_\alpha$ is developed in Section \ref{sec:quandleTQFT}. Sections \ref{sec:SL2repvarietyknots} and \ref{sec:Rileyalgebra} particularize the previous construction to $\SL_2$ and its symplectic quandle description. The Riley algebra is introduced in Section \ref{sec:Rileyalgebra}. Finally, Section \ref{sec:parabolicpretzelknots} describes the representation and character varieties of the pretzel knot family.
\subsection*{Acknowledgements}

The first-named and the second-named authors have been partially supported by Ministerio de Ciencia e Innovaci\'on Project PID2024-156578NB-I00 funded by MICIU /AEI /10.13039/501100011033 / FEDER, EU. The second-named author has also been partially supported by the bilateral AEI-DFG project Celestial Mechanics, Hydrodynamics, and Turing Machines / Himmelsmechanik, Hydrodynamik und Turing-Maschinen (AQUACELL), with reference codes PCI2024-
155042-2 and PCI2024-155062-2.

\subsection*{AI declaration} Generative Artificial Intelligence, particularly OpenAI's GPT-5.6 and GPT-6 in their versions Sol and Astra, has been used throughout the development of this work for literature review, brainstorming, including preliminary calculations for the description of Section \ref{sec:closing-pretzel}, and generation and documentation of an initial draft of the accompanying SageMath script. Additionally, GPT-6 Astra was used in the final stage for mathematical assessment and proofreading of
the manuscript.

The research design, formulation and proofs of the results, and writing of the manuscript were entirely carried out by the authors, who take full responsibility for the content of the work. No funding from AI companies was received during the development of this work.

\section{Quandle representation varieties} \label{sec:quandlerepvarieties}

Let $\Gamma$ be a finitely generated group, and consider a complex reductive group $G$. The \emph{representation variety} is the set
$$
    \R_G(\Gamma) = \Hom_{\Grp}(\Gamma, G)
$$
of group homomorphisms $\rho: \Gamma \to G$. This set has a natural complex algebraic structure given as follows. If $\Gamma = \langle \gamma_1, \ldots, \gamma_n \,\mid\, r_\alpha(\gamma_1, \ldots, \gamma_n)\rangle$, for some relations $r_\alpha$, then $\R_G(\Gamma)$ can be identified with the algebraic set of $G^n$ given by
$$
\R_G(\Gamma) = \{(g_1, \ldots, g_n) \in G^n \,\mid\, r_\alpha(g_1, \ldots, g_n) = 1\} \subseteq G^n.
$$

Furthermore, notice that the group $G$ acts on $\R_G(\Gamma)$ by conjugation, i.e.\ $(g \cdot \rho)(\gamma) = g\rho(\gamma)g^{-1}$ for $g \in G$, $\rho: \Gamma \to G$ and $\gamma \in \Gamma$. In this sense, the Geometric Invariant Theory (GIT) quotient
$$
    \X_G(\Gamma) = \R_G(\Gamma) \sslash G, 
$$
is the so-called \emph{character variety}, and parametrizes representations of $\Gamma$ to $G$ up to isomorphism.

A particularly interesting case arises when we consider a link $L \subseteq S^3$. The fundamental group of its complement
$$
    \Gamma_L := \pi_1(S^3-L)
$$
is known as the \emph{link group}. The associated representation and character varieties
$$
    \R_G(L) := \R_G(\Gamma_L) = \Hom_{\Grp}(\Gamma_L, G), \qquad \X_G(L) := \R_G(L) \sslash G,
$$
are important spaces that can be exploited to understand the geometry of the link $L$, see for instance \cite{CullerShalen1983,Cooper1994Planecurves}. 

In this section, we will explore a generalization of these representation and character varieties by using quandles instead of groups. These are objects of intense research due to their intricate geometry and their relation with classical knot invariants, such as the Alexander polynomial \cite{AGL1Alexanderinvariants}.

\subsection{Review of quandles} In this section, we shall briefly review the theory of quandles, with a view towards its application to knot theory. For a more detailed account, see for instance \cite{Joycequandle}.

\begin{defn}
A \emph{quandle} is a set $Q$ equipped with a binary operation $\triangleright: Q\times Q \longrightarrow Q$ satisfying:
\begin{itemize}
\item[(P1)] $x\triangleright x= x$, for all $x\in Q$.
\item[(P2)] For all $x,z\in Q$, there exists a unique $y\in Q$ such that $x\triangleright y=z$.
\item[(P3)] $ x\triangleright (y\triangleright z) = (x\triangleright y)\triangleright (x\triangleright z)$, for all $x,y,z\in Q$.
\end{itemize}
\end{defn} 
We note that uniqueness in (P2) leads to a second operation $\triangleleft$, defined as $z \triangleleft x = y$ if and only if $x \triangleright y =z$.
We will also say that a \emph{quandle homomorphism} between two quandles $Q$ and $Q'$ is a map $f: Q \to Q'$ such that $f(x \triangleright y) = f(x) \triangleright f(y)$. With these morphisms, quandles define a category, denoted by $\Qdl$. Given a set $S$, we can always endow $S$ with a quandle operation by setting $x \triangleright y = y$ (so also $y \triangleleft x = y$), called the \emph{trivial quandle} on $S$ and also denoted by $S$.

The representation theory of groups and quandles is intimately related by an adjunction between the category of groups $\mathbf{Grp}$ and the category of quandles $\Qdl$. To be precise, a group $G$ gives rise to a quandle, the so-called \emph{conjugation quandle} $\Conj G$, with underlying set $G$ and quandle operation given by conjugation, that is,
$$
x \triangleright y = xyx^{-1}, \quad \text{for } x,y \in G.
$$
Notice that, in this case, $y \triangleleft x = x^{-1}yx$.
Conversely, any quandle $Q$ defines an \emph{associated group} as the group generated by the elements of $Q$, subject to the relation that $x \triangleright y = xyx^{-1}$ for all $x, y \in Q$. Explicitly,
$$
\As Q = \left\langle x \in Q \,\middle|\, x \triangleright y = xyx^{-1},\ \text{for all } x,y \in Q \right\rangle.
$$
These two functors form an adjunction (we refer to \cite{racks} for details)
\begin{equation} \label{eqn:adjunctiongrqdl}
\Hom_{\mathbf{Grp}}(\As Q, G) \cong \Hom_{\Qdl}(Q, \Conj G).
\end{equation}

\begin{ex} \label{ex:quandlefiber}
Let $G$ be a group, and consider a map $f: G \to B$ (not necessarily a group homorphism) that is invariant under conjugation, i.e.\ $f(xyx^{-1}) = f(y)$ for all $x,y \in G$. Then, for any $b \in B$, the fiber $G_b = f^{-1}(b)$ inherits a subquandle structure from the one on $\Conj G$ since, for $x, y \in G_b$, both $x \triangleright y = xyx^{-1} \in G_b$ and $y \triangleleft x = x^{-1}yx \in G_b$. We will call these quandles the \emph{slice conjugation quandles}.
\end{ex}

\subsection{The link quandle}

\begin{figure}[h]
    \centering
    \includegraphics[width=0.4\linewidth]{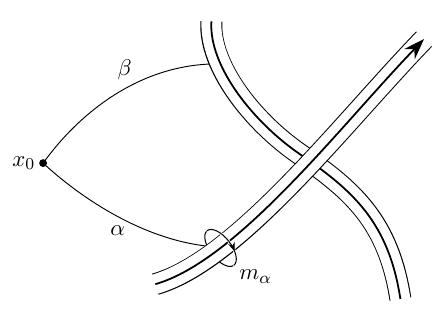}
    \caption{The quandle operation $[\alpha]\triangleright[\beta]$.}
    \label{fig:quandle_crossing}
\end{figure}

A particularly interesting quandle in knot theory is the so-called \emph{link quandle} $Q_L$ of a link $L \subseteq S^3$, as introduced by \cite{Matveev,Joycequandle}, or in general of a \emph{generalized link}, that is, a cooriented submanifold $L \subseteq M$ of codimension $2$. The coorientation means that the normal bundle of $L$ is oriented and thus it induces an orientation on the meridians of a tubular neighbourhood $N_L \subseteq M$ of $L$. Denote by $E_L=M-N_L$ the exterior of this neighbourhood. Fixing a base point $x_0\in E_L$, we define $\mathcal{P}_L$ as the set of homotopy classes of paths in $E_L$ from $x_0$ to $\partial N_L$, where homotopies are required to keep the initial point $x_0$ fixed while allowing the endpoint to move freely in $\partial N_L$. Furthermore, note that any point of $\partial N_L$ belongs to a unique meridian circle of the normal circle bundle, so any $\alpha \in \mathcal{P}_L$ defines a unique loop $m_\alpha$ in $\partial N_L$ based at the endpoint of $\alpha$ and positively oriented with respect to the orientation of $N_L$. This allows us to define an operation on $\mathcal{P}_L$ that, given paths $\alpha, \beta \in \mathcal{P}_L$, assigns
$$
\alpha \triangleright \beta = \alpha \star m_{\alpha} \star \bar{ \alpha}\star \beta,
$$
where $\star$ denotes the standard concatenation of paths, $\bar{\alpha}$ denotes the reverse path to $\alpha$, and we choose the same path representative
for $\alpha$ and $\bar{\alpha}$, and the path $m_{\alpha}$ at the endpoint of $\alpha$ (see Figure \ref{fig:quandle_crossing}). 
This operation endows $Q_L = (\mathcal{P}_L,\triangleright)$ with a quandle structure, called the \emph{quandle of $L$}. Notice that $Q_L$ only depends on the chosen basepoint up to isomorphism.

The construction of the quandle $Q_L$ is functorial in the following sense. Let $L \subseteq M$ as above, and suppose that $N \subseteq M$ is an oriented submanifold that intersects $L$ transversely. Then $L \cap N$ is a submanifold of $N$ of codimension two with a naturally induced orientation. Then the quandle $Q_{L \cap N}$ is well defined and the inclusion induces a quandle morphism $Q_{L \cap N} \to Q_L$.

\begin{figure}[h]
    \centering
    \includegraphics[width=0.5\linewidth]{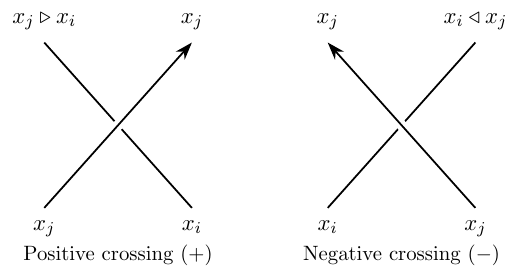}
    \caption{The quandle operation for simple crossings}
    \label{fig:quandle_operations}
\end{figure}

In the case of links $L \subseteq S^3$, it is possible to give a Wirtinger-type presentation of the quandle of $L$, as proven in \cite{fenn1992racks}. Consider a regular link diagram of $L$ with arcs $x_1, \ldots, x_n$. Then, $Q_L$ is the quandle generated by $x_1, \ldots, x_n$ and the relations that $x_j\triangleright x_i=x_{i+1}$, for each positive crossing $L_+$ as in Figure \ref{fig:quandle_operations}, and $x_i\triangleleft x_j=x_{i+1}$ for each negative crossing $L_-$. In particular, this proves that $Q_L$ is finitely generated. In fact, for knots $K \subseteq S^3$, $Q_K$ is known to be a complete knot invariant up to orientation \cite{Joycequandle,Matveev}. 

This link quandle is strongly related to the link group, i.e.\ the fundamental group of the complement $\Gamma_L := \pi_1(S^3-L)$. Indeed, from the fact that $\Gamma_L$ can be presented as the group generated by the arcs $x_1, \ldots, x_n$ with relations $x_jx_ix_j^{-1}=x_{i+1}$ for each positive crossing $L_+$, and $x_j^{-1}x_i x_j=x_{i+1}$ for each negative crossing $L_-$, we directly get an isomorphism of groups
$$
    \Gamma_L \cong \As Q_L.
$$

\subsection{Quandle representation varieties}

A quandle $X$ is said to be \emph{complex linear algebraic} if it is a quandle object in the category of complex affine varieties, i.e.\ if $X$ is a complex affine variety and the operations $\triangleright$ and $\triangleleft$ are both regular morphisms.

\begin{defn}
    Let $X$ be a complex linear algebraic quandle and let $L \subseteq M$ be a generalized link whose link quandle $Q_L$ is finitely generated. The \emph{$X$-quandle representation variety} is the space
 $$
\R_{X}(L) := \Hom_{\Qdl}(Q_L, X).
$$ 
\end{defn}

As in the case of representation varieties of groups, by considering a presentation of $Q_L$, it is easy to check that $\R_X(L)$ is equipped with a natural structure of a complex affine algebraic variety, and that this structure is independent of the chosen presentation up to isomorphism. Notice in particular that if $f: X \to X'$ is a homomorphism of complex linear algebraic quandles, then we get a regular morphism
$$
    f_{\#}: \R_X(L) = \Hom_{\Qdl}(Q_L, X) \to \R_{X'}(L) = \Hom_{\Qdl}(Q_L, X'), \qquad \rho \mapsto f \circ \rho.
$$

The quandle representation variety is particularly interesting in the case that $X = \Conj G$ is the conjugation quandle of a complex linear group $G$. In that case, for any generalized link $L$ such that $\Gamma_L = \As Q_L$ (for instance, for links $L \subseteq S^3$), by adjunction we get an isomorphism
\begin{equation}\label{eq:iso-qdl-grp}
    \R_{\Conj G}(L) = \Hom_{\Qdl}(Q_L, \Conj G) = \Hom_{\Grp}(\As Q_L, G) = \Hom_{\Grp}(\Gamma_L, G) = \R_G(L),
\end{equation}
between the $\Conj G$-quandle representation variety $\R_{\Conj G}(L)$ and the standard $G$-representation variety $\R_G(L)$.

The isomorphism (\ref{eq:iso-qdl-grp}) can be exploited in the following setting. 
Let $G$ be a connected reductive group and choose a Borel subgroup $B \subseteq G$, which determines a maximal torus $T \subseteq G$ and its Weyl group $W$ acting on $T$. Let $C_G:=G\sslash G= \Spec \mathbb{C}[G]^G \cong T/W$, with respect to the action of $G$ on itself by conjugation. The \emph{Chevalley map} is the regular morphism induced by the quotient
$$
    \chi: G \to C_G.
$$
Notice that, in general, $\chi$ is not a group homomorphism. For instance, for $G = \SL_{r}$, we have $C_{\SL_r} = T/W = \mathbb{A}^{r-1}$ and the Chevalley map is given by the coefficients of the characteristic polynomial of any element of $G$. Since $C_G$ parametrizes closed semisimple conjugacy classes and not an actual conjugacy class, we will refer to it throughout this paper as the \emph{spectral base}. 

Now, let us endow $C_G$ with the trivial quandle structure. Since the Chevalley map $\chi: G \to C_G = G \sslash G$ is invariant for the conjugacy action, the Chevalley map becomes a quandle homomorphism
\begin{equation}\label{eq:chevalley-quandles}
    \chi: \Conj G \to C_G.
\end{equation}
In particular, this homomorphism induces a regular morphism of representation varieties
\begin{equation}\label{eq:chevalley}
    \chi_{\#}: \R_{G}(L) \cong \R_{\Conj G}(L) \longrightarrow \R_{C_G}(L).
\end{equation}

\begin{prop}
    Let $K \subseteq S^3$ be a knot with knot quandle $Q_K$. For any trivial quandle $S$, we have an isomorphism
    $$
        \R_S(K) \cong S.
    $$
    If we orient $K$, this isomorphism is canonical and is given by the image of any oriented meridian.
\end{prop}

\begin{proof}
Consider the associated knot group $\Gamma_K = \pi_1(S^3-K)$. Then, $\Gamma_K$ is generated by the meridians $x_1, \ldots, x_n,$ and all of them are conjugate. This implies that, in $Q_K$, which is also generated by $x_1, \ldots, x_n$, for all $i = 1, \ldots n$, we have a relation of the form $x_{j} \triangleright x_{i} = x_{i+1}$ or $x_{i} \triangleleft x_{j} = x_{i+1}$, where we denote $x_{n+1} = x_1$. This implies that a quandle homomorphism $\rho: Q_K \to S$ must satisfy, for each relation of the form $x_{j} \triangleright x_{i} = x_{i+1}$, the equality
$$
    \rho(x_{i+1}) = \rho(x_{j} \triangleright x_{i}) = \rho(x_{j}) \triangleright \rho(x_{i}) = \rho(x_i),
$$
and similarly $\rho(x_{i+1}) = \rho(x_i)$ for a relation $x_{i} \triangleleft x_{j} = x_{i+1}$. Therefore, $\rho(x_1) = \rho(x_2) = \ldots = \rho(x_n)$ and thus we have an isomorphism
$$
    \R_{S}(K) \cong S.
$$
By construction, this isomorphism is given by the image of a meridian and, thus it is canonical provided that we choose an orientation of $K$, and thus of an oriented meridian.
\end{proof}

In particular, for $S = C_G$, we have $\R_{C_G}(K) \cong C_G$ for a knot $K$. More generally, if $L$ is a link with $s$ connected components and link quandle $Q_L$, then all the generators of $Q_L$ corresponding to the same connected component are conjugate and thus we get an isomorphism 
    $$
        \R_{C_G}(L) \cong C_G^s.
    $$
This isomorphism is canonical once we choose an orientation for $L$.

Now, let us fix an oriented knot $K$, so the Chevalley map (\ref{eq:chevalley}) becomes a morphism
$$
    \chi_{\#}: \R_{G}(K) \to C_G.
$$
Given $t \in C_G$, let us denote the fiber of $\chi_{\#}$ at $t$ by $\R_{G}(K)_t := \chi_{\#}^{-1}(t)$. Recall also that, by Example \ref{ex:quandlefiber}, for $t \in C_G$, the slice $G_t = \chi^{-1}(t) \subseteq G$ acquires a quandle structure as subquandle of $\Conj G$. It turns out that this fiber can be identified itself as a quandle representation variety.

\begin{prop}\label{prop:quandle-slice}
Let $K \subseteq S^3$ be a knot. For any $t \in C_G$, we have an isomorphism
$$
    \R_{G}(K)_t \cong \R_{G_t}(K).
$$
\end{prop}

\begin{proof}
Observe that we have a natural inclusion given by
\begin{equation}\label{eq:inclusion-quandle}
    \R_{G_t}(K) = \Hom_{\Qdl}(Q_K, G_t) \hookrightarrow \Hom_{\Qdl}(Q_K, \Conj G) = \R_G(K).
\end{equation}

Furthermore, if $\rho \in \R_{G_t}(K)$, then for any meridian $m$ we have that $\rho(m) \in G_t$, so $\chi(\rho(m)) = t$. This implies that $\R_{G_t}(K) \subseteq \R_G(K)_t$ under this inclusion. For the opposite inclusion, notice that in a knot, the knot quandle $Q_K$ is generated by meridians, so if a group representation $\rho: \Gamma_K \to G$ satisfies that $\rho(m) \in G_t$ for every meridian $m$ then, seen as a quandle representation $\rho: Q_K \to \Conj G$, it also satisfies $\textup{Im}(\rho) \subseteq G_t$.
\end{proof}

\begin{rem}
    The proof of Proposition \ref{prop:quandle-slice} can be readily generalized to any link $L \subseteq S^3$ with $s$ components. Indeed, there is still an inclusion
    $$
        \R_{G_t}(L) \hookrightarrow \R_G(L).
    $$
    The image of this inclusion coincides with the fiber of $\chi_{\#}: \R_G(L) \to \R_{C_G}(L) = C_G^s$ over the point $(t, \ldots, t) \in C_G^s$. However, in this setting, not every fiber of $\chi_{\#}$ admits such a quandle interpretation, since elements outside the diagonal of $C_G^s$ correspond to representations sending meridians in different components of $L$ to different slices $G_t \subseteq G$.
\end{rem}

\section{The quandle Topological Quantum Field Theory} \label{sec:quandleTQFT}

In this section, we shall construct a Topological Quantum Field Theory that computes the quandle representation variety of a link, as a functor with values in the category of spans of algebraic varieties. As we will see, the construction is explicit and allows a precise description of the image of the generators in terms of the quandle structure.

\subsection{Tangles and the bordered quandle category}
We briefly recall the definition of the tangle category $\Tang$; for a detailed account see \cite{Kasselquantumgroups}. 

For $n\geq 1$, let $[n]=\{1,\dots,n\}$ and $[0]=\emptyset$. A \emph{tangle} $T$ is a finite union of disjoint oriented polygonal arcs (possibly including closed loops) embedded in $\mathbb{R}^2\times[0,1]$, whose boundary satisfies
\[
\partial T = T \cap (\mathbb{R}^2\times\{0,1\}) = ([k]\times\{0\}\times\{0\}) \cup ([\ell]\times\{0\}\times\{1\})
\]
for some $k,\ell\geq 0$.

Tangles form a category $\Tang$ as follows. Objects are finite (possibly empty) sequences of signs $\epsilon=(\epsilon_1,\dots,\epsilon_n)$ with $\epsilon_i=\pm$, and morphisms are isotopy classes of oriented tangles. Each tangle $T$ determines a source $s(T)$ and target $b(T)$ given by the signed boundary points at the bottom and top, respectively, where the sign records whether the strand is incoming or outgoing. Composition is defined by vertical concatenation of tangles, defined whenever $b(T)=s(T')$. The identity morphism $\mathrm{id}_\epsilon$ consists of $|\epsilon|$ vertical oriented strands.

From a structural point of view, the category $\Tang$ is monoidal: the tensor product on objects is given by concatenation of sequences, with unit the empty sequence, and on morphisms by horizontal juxtaposition. This monoidal structure is not symmetric but braided, with braiding given by the usual crossings of strands.

\begin{figure}[h]
    \centering
    \includegraphics[width=0.4\linewidth]{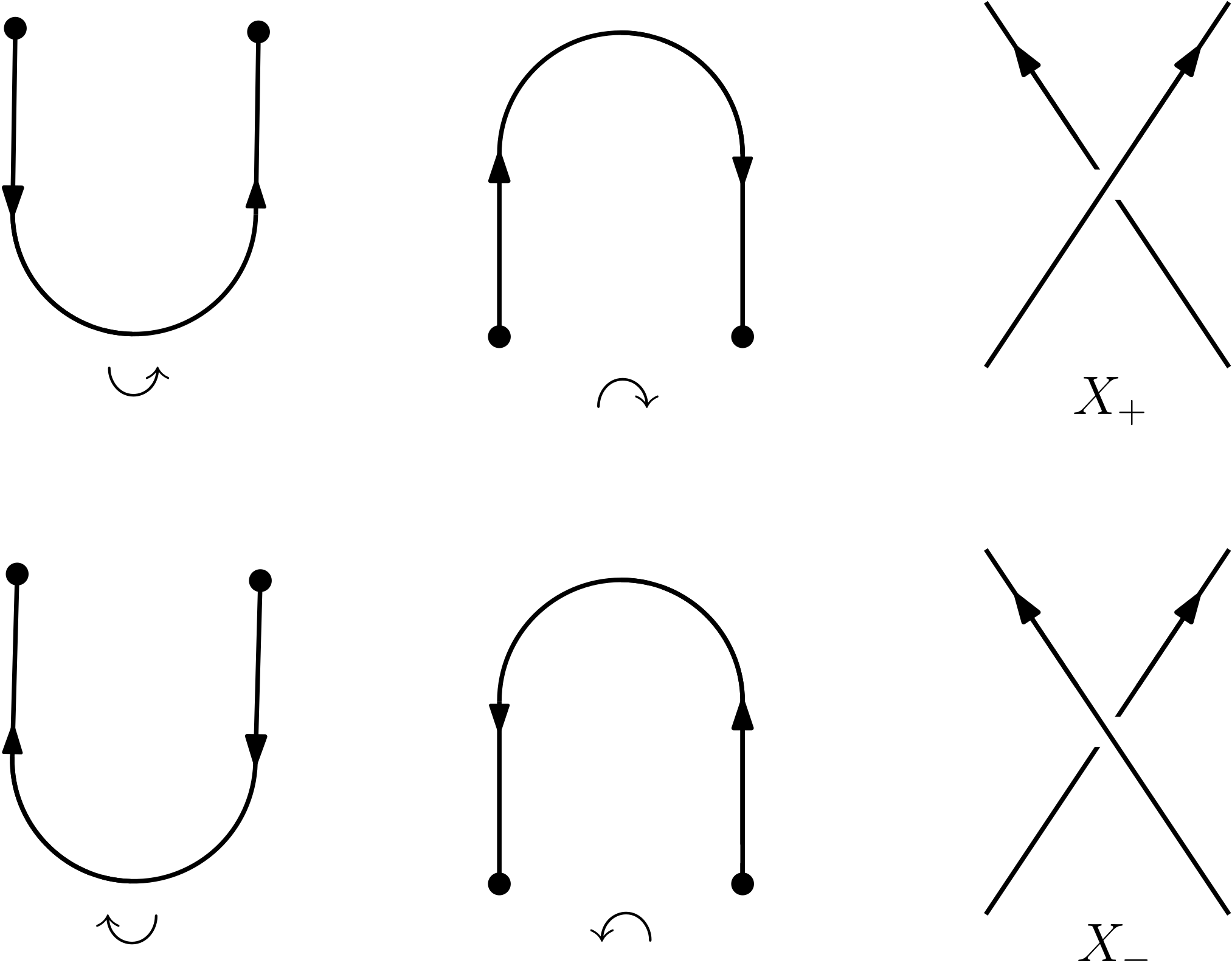}
    \caption{Generating tangles in $\Tang$.}
    \label{fig:tangles}
\end{figure}

Tangles admit planar diagrammatic representations in $\mathbb{R}\times[0,1]$, where crossings are recorded, and isotopy classes correspond to equivalence under Reidemeister moves. Finally, $\Tang$ admits a presentation as a monoidal category by generators and relations. It is generated by six elementary morphisms
\[
\upcurvearrowright,\ \upcurvearrowleft,\ \curvearrowright,\ \curvearrowleft,\ X_{+},\ X_{-},
\]
corresponding to cups, caps, and positive and negative crossings, as depicted in Figure \ref{fig:tangles}.

These generators satisfy relations arising from planar isotopies and the Reidemeister moves (see \cite[Theorem XII.2.2]{Kasselquantumgroups}), and every tangle can be obtained by composing and tensoring these elementary pieces.

\subsection{Representations of the bordered tangle category} \label{subsec:representationsborderedtanglecat} 
As shown in \cite{cattabriga}, any tangle $T: \epsilon_1 \to \epsilon_2$ gives rise to two quandle morphisms
$$
    i_1: Q_{\epsilon_1} \to Q_T, \qquad i_2: Q_{\epsilon_2} \to Q_T,
$$
where $Q_{\epsilon_1}, Q_{\epsilon_2}$ are the link quandles of the boundaries of $T$ inside $\RR^2$, and $Q_T$ is the link quandle of $T$ inside $\RR^2 \times [0,1]$. Observe that both $Q_{\epsilon_1}$ and $Q_{\epsilon_2}$ are free quandles with as many generators as the number of points $|\epsilon_1|$ and $|\epsilon_2|$.

\begin{rem}
This construction is categorized in \cite{cattabriga} as the \emph{bordered quandle category} $\mathcal{BQ}$, where free quandles (with finitely many generators) are objects, and morphisms between these objects $F_1,F_2$ are given by triples $F_1 \stackrel{i_1}{\rightarrow} Q \stackrel{i_2}{\leftarrow} F_2$, where $Q$ is an additional quandle.
\end{rem}

Now, if we fix a complex linear algebraic quandle $X$, and represent any of these previous triples onto $X$, we obtain that we can associate to any tangle $T$ the triple
\begin{equation}\label{eq:span}
\xymatrix{
& \R_{X}(T) = \Hom_{\Qdl}(Q_T,X) \ar[dl]_{i_1} \ar[dr]^{i_2} & \\
\R_X(\epsilon_1) = \Hom_{\Qdl}(Q_{\epsilon_1},X) & & \R_X(\epsilon_2) = \Hom_{\Qdl}(Q_{\epsilon_2},X)
}
\end{equation}

Notice that these spaces and morphisms are objects in the category $\Varc$ of complex algebraic varieties. In this manner, the aforementioned information can be seen as a morphism in the category $\Span(\Varc)$ of spans of complex varieties. Recall that given a category $\mathcal{C}$ with pullbacks, the category of spans, $\Span(\mathcal{C})$ is given by:
\begin{itemize}
    \item Objects: The objects are objects of $\mathcal{C}$.
    \item Morphisms: Given objects $c_1$ and $c_2$ of $\mathcal{C}$, a morphism between them is an equivalence class of pairs of morphisms in $\mathcal{C}$, known as a span
\begin{equation*}
\xymatrix{
& d \ar[dl]_{f} \ar[dr]^{g} & \\
c_1 & & c_2
}
\end{equation*}
    Two spans $c_1 \stackrel{f}{\leftarrow} d \stackrel{g}{\rightarrow}c_2$ and $c_1 \stackrel{f'}{\leftarrow} d' \stackrel{g'}{\rightarrow}c_2$ are declared to be equivalent if there exists an isomorphism $h: d \to d'$ such that $f = f'\circ h$ and $g = g' \circ h$. Composition of spans is given by pullback.
\end{itemize}
Notice that if $\mathcal{C}$ is monoidal, then $\Span(\mathcal{C})$ also inherits a monoidal structure. In particular, $\Span(\Varc)$ is a monoidal category with the cartesian products of both varieties and spans.

With this language, fixed a complex linear algebraic quandle $X$, we can pack the previous construction as a monoidal functor
$$
    Z_X: \Tang \longrightarrow \Span(\Varc),
$$
which assigns, to every object $\epsilon$ of $\Tang$, the variety $Z_X(\epsilon)=\R_X(\epsilon)$ and, to a tangle $T: \epsilon_1 \to \epsilon_2$, the span (\ref{eq:span}). The same proof of \cite{AGL1Alexanderinvariants} can be adapted here verbatim to prove that $Z_X$ is actually a functor, which is also monoidal.

\subsection{Relative construction and spectral data}\label{sec:TQFT-relative}

In this section, we are going to construct a relative version of the functor built in Section \ref{subsec:representationsborderedtanglecat} that will allow us to fix the spectral data in quandle representation varieties.

In order to do so, we are going to work on the arrow category $\Arr{\Varc}$, whose objects are regular morphisms $Z \to Z'$ and a morphism between $f_1: Z_1 \to Z_1'$ and $f_2: Z_2 \to Z_2'$ is a commutative square
$$
\xymatrix{
Z_1 \ar[d]^{f_1} \ar[r] & Z_2 \ar[d]^{f_2} \\
Z_1' \ar[r] & Z_2'
}
$$
In other words, the category $\Arr{\Varc}$ is the comma category $\left(\textup{Id}_{\Varc} \downarrow \textup{Id}_{\Varc}\right)$.

Now, let us consider complex algebraic quandles $X$ and $X'$ and fix a quandle morphism
$$
    \alpha: X \to X'.
$$
We construct a monoidal functor
$$
    Z_\alpha: \Tang \to \Span(\Arr{\Varc}).
$$
This functor is given as follows:
\begin{itemize}
    \item Given an object $\epsilon$ of $\Tang$, we assign the associated morphism at the level of representation varieties
    $$
        \alpha_{\#}: \R_X(\epsilon) \to \R_{X'}(\epsilon).
    $$
    \item For a tangle $T: \epsilon_1 \to \epsilon_2$, we assign the double span
    $$
\xymatrix{
& \R_{X}(T) \ar[d] \ar[dl]_{i_1} \ar[dr]^{i_2} & \\
\R_X(\epsilon_1) \ar[d] & \R_{X'}(T) \ar[dl]_{i_1} \ar[dr]^{i_2'} & \R_X(\epsilon_2) \ar[d] \\
\R_{X'}(\epsilon_1)  & & \R_{X'}(\epsilon_2) 
}
    $$
    Here, the vertical maps are the morphisms of representation varieties induced by $\alpha$ and the diagonal maps are the natural restriction maps to the boundaries. Notice that this double span can be seen as span of morphisms of $\Arr{\Varc}$.
\end{itemize}

\begin{thm} \label{thm:Zalphafunctor}
    The assignment $Z_\alpha: \Tang \to \Span(\Arr{\Varc})$ is a monoidal functor.
\end{thm}
\begin{proof}
    Let us first prove that $Z_\alpha$ preserves composition. Let $T: \epsilon_1 \to \epsilon_2$ and $T': \epsilon_2 \to \epsilon_3$ be tangles. The respective complements $\left(\RR^2 \times [0,1]\right) - T$ and $\left(\RR^2 \times [1,2]\right) - T'$ are connected, and their intersection is $\left(\RR^2 \times \{1\} \right) - \epsilon_2$, which is also connected. Then, by the Seifert-van Kampen theorem for quandles \cite[Theorem 13.1]{Joycequandle}, we have that the quandle of the composition $T' \circ T: \epsilon_1 \to \epsilon_3$ is the pushout in $\Qdl$
    $$
        \xymatrix{
        Q_{T' \circ T}   & Q_{T'} \ar[l] \\
        Q_{T} \ar[u] & Q_{\epsilon_2} \ar[l]\ar[u]
    }
    $$
    Therefore, since $\Hom_{\Qdl}(-, X)$ and $\Hom_{\Qdl}(-, X')$ are continuous functors, we get a commutative diagram
    $$
        \xymatrix{
        \R_X(T' \circ T) \ar[rr]\ar[d]\ar@{--{>}}[ddr]  & & \R_X(T') \ar[d] \ar@{--{>}}[ddr] & \\
        \R_X(T) \ar[rr]\ar@{--{>}}[ddr]  & & \R_X(\epsilon_2) \ar@{--{>}}[ddr]& \\
        & \R_{X'}(T' \circ T) \ar[rr]\ar[d]  & & \R_{X'}(T') \ar[d]  \\
        & \R_{X'}(T) \ar[rr] & &  \R_{X'}(\epsilon_2) 
        }
    $$
    where the two solid squares are cartesian. We therefore get that the composition $Z_\alpha(T') \circ Z_\alpha(T)$ corresponds to the solid lines in the following commutative diagram
    $$
        \xymatrix{
        \R_X(T' \circ T) \ar[rr]\ar[d]\ar[ddr]  & & \R_X(T') \ar@{--{>}}[d] \ar@{--{>}}[ddr] \ar[r] & \R_{X}(\epsilon_3)  \ar[ddr] & \\
        \R_X(T) \ar@{--{>}}[rr]\ar@{--{>}}[ddr] \ar[d]  & & \R_X(\epsilon_2) \ar@{--{>}}[ddr]& & \\
        \R_{X}(\epsilon_1) \ar[ddr] & \R_{X'}(T' \circ T) \ar[rr]\ar[d]  & & \R_{X'}(T') \ar@{--{>}}[d]\ar[r] & \R_{X'}(\epsilon_3)  \\
        & \R_{X'}(T) \ar@{--{>}}[rr] \ar[d] & &  \R_{X'}(\epsilon_2) &  \\
        & \R_{X'}(\epsilon_1)  & & & 
        }
    $$
    which agrees with $Z_\alpha(T' \circ T)$, proving the functoriality. Now, monoidality follows similarly from the fact that, for any objects $\epsilon$ and $\epsilon'$ of $\Tang$, we have
    $$
        Z_\alpha(\epsilon \otimes \epsilon') = Z_\alpha(\epsilon \cup \epsilon') = \left[\R_X(\epsilon \cup \epsilon') \to \R_{X'}(\epsilon \cup \epsilon')\right] = \left[\R_X(\epsilon) \to \R_{X'}(\epsilon)\right] \times \left[\R_X(\epsilon') \to \R_{X'}(\epsilon')\right],
    $$
    and similarly for tangles.
\end{proof}

In particular, let us now consider the Chevalley map (\ref{eq:chevalley-quandles}),
$$
    \chi: \Conj G \to C_G = G\sslash G.
$$
The previous construction produces a monoidal functor
$$
    Z_{\chi}: \Tang \to \Span(\Arr{\Varc})
$$
which assigns, for an object $\epsilon$ of $\Tang$, the morphism
$$
    Z_{\chi}(\epsilon): \R_{G}(\epsilon) = G^{|\epsilon|} \longrightarrow \R_{C_G}(\epsilon) = C_G^{|\epsilon|},
$$
where $|\epsilon|$ is the number of points of $\epsilon$. Furthermore, for a tangle $T: \epsilon_1 \to \epsilon_2$ with $s$ connected components, we get a double span
\begin{equation}\label{eq:double-span}
\begin{split}
\xymatrix{
& \R_{G}(T) \ar[d] \ar[dl]_{i_1} \ar[dr]^{i_2} & \\
\R_G(\epsilon_1) = G^{|\epsilon_1|} \ar[d] & \R_{C_G}(T) = C_G^s \ar[dl]_{i_1'} \ar[dr]^{i_2'} & \R_G(\epsilon_2) = G^{|\epsilon_2|} \ar[d] \\
\R_{C_G}(\epsilon_1) = C_G^{|\epsilon_1|} & & \R_{C_G}(\epsilon_2) = C_G^{|\epsilon_2|}
}
\end{split}
\end{equation}
Observe that all the identifications $\R_G(\epsilon_1) = G^{|\epsilon_1|}, \R_G(\epsilon_2) = G^{|\epsilon_2|}, \R_{C_G}(\epsilon_1) = C_G^{|\epsilon_1|}, \R_{C_G}(\epsilon_2) = C_G^{|\epsilon_2|}$ and $\R_{C_G}(T) = C_G^s$ are canonical since $T$ and its boundaries are oriented.

Now, given $t \in C_G$, let us also denote by $t$ the diagonal element $(t, \ldots, t) \in C_G^s$ for any $s \geq 1$ and for a morphism $\chi_{\#}:\R_G(L) \to \R_{C_G}(L) = C_G^s$, let $\R_G(L)_t = \chi_{\#}^{-1}(t)$ be its fiber at $t = (t, \ldots, t) \in C_G^s$. Furthermore, recall that this fiber can be identified with the quandle representation variety $\R_G(L)_t = \R_{G_t}(L)$ corresponding to the slice quandle $G_t = \chi^{-1}(t) \subseteq \Conj G$ by Proposition \ref{prop:quandle-slice}. In this manner, taking the fiber at $t$ in the double span (\ref{eq:double-span}), we get a span
$$
\xymatrix{
& \R_{G}(T)_t = \R_{G_t}(T) \ar[dl]_{i_1} \ar[dr]^{i_2} & \\
\R_G(\epsilon_1)_t = \R_{G_t}(\epsilon_1) = G_t^{|\epsilon_1|} & & \R_G(\epsilon_2)_t = \R_{G_t}(\epsilon_2) = G_t^{|\epsilon_2|}
}
$$
Thanks to this quandle interpretation as representation varieties of $G_t$, the previous span is actually the span $Z_{G_t}(T)$ corresponding to the functor constructed in Section \ref{subsec:representationsborderedtanglecat} for the quandle $X = G_t$. In particular, this allows us to see the functors $Z_{G_t}: \Tang \to \Span(\Varc)$ as the ``fibers'' of the functor $Z_{\chi}: \Tang \to \Span(\Arr{\Varc})$ at $t \in C_G$.

\subsection{Explicit operation on the generators of $\Tang$.}\label{sec:explicit-functor}
Let $X$ be an complex linear algebraic quandle. In this section, we will explicitly describe the image of the functor $Z_{X}: \Tang \to \Span(\Varc)$ on some tangles by fixing natural representatives for the generators of the fundamental groups associated with objects and morphisms in $\Tang$. 

For an object $\epsilon$ with $n=|\epsilon|$, the complement $\epsilon^c = \RR^2 - \epsilon \subseteq \mathbb{R}^2$ has as fundamental group $\Gamma_\epsilon=\pi_1(\epsilon^c)\cong \Free{n}$, the free group in $n$ generators. Each generator is a simple loop around each point of $\epsilon$, oriented according to the sign of each point, so that from a top view the loops wind clockwise or counterclockwise. This convention also determines compatible representatives in the quandle $Q_\epsilon$, giving a natural identification
$$
\R_{X}(\epsilon) = \Hom_{\Qdl}(Q_{\epsilon}, X) \cong X^n.
$$

\begin{figure}[h]
    \centering
    \includegraphics[width=0.6\linewidth]{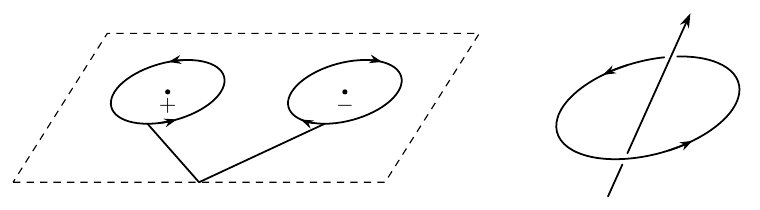}
    \caption{Generators of $\Gamma_\epsilon$ and $\Gamma_T$.}
    \label{fig:pi1generators}
\end{figure}

Similarly, for a tangle $T:\epsilon_1 \to \epsilon_2$, the complement $T^c = \RR^2 \times [0,1]-T$ has fundamental group $\Gamma_T=\pi_1(T^c)$. We take as generators the meridians around the strands of $T$, oriented by the right-hand rule (see Figure \ref{fig:pi1generators}).

These choices again determine representatives in the quandle $Q_T$, and hence a concrete description of $Z_{X}$. In particular, we will consider the following tangles, as depicted in Figure \ref{fig:tangles},
$$
\begin{array}{lll}
    \upcurvearrowright: \emptyset \to (-,+), & \curvearrowright: (-,+) \to \emptyset, & X_+: (+,+) \to (+,+), \\
    \upcurvearrowleft: \emptyset \to (+,-), & \curvearrowleft: (+,-) \to \emptyset, & X_-: (+,+) \to (+,+).
\end{array}
$$

In this way, applying the functor $Z_X$, we directly obtain:
\begin{itemize}
    \item Evaluations and coevaluations: It is shown that $Z_{X}(\curvearrowleft) = Z_{X}(\curvearrowright)$ and also $Z_{X}(\upcurvearrowleft)=Z_{X}(\upcurvearrowright)$, which are respectively given by the spans
$$
\xymatrix{
& X \ar[ld]_{\Delta} \ar[rd]& \\
X \times X & & \ast
}\quad
\xymatrix{
& X \ar[rd]^{\Delta} \ar[ld]& \\
\ast & & X \times X
}
$$
 where $\Delta: X \longrightarrow X \times X$ is the diagonal map $\Delta(x)=(x,x)$ for $x \in X$.
\item Simple crossings $X_+$ and $X_-$: Since the inclusion of $\epsilon^c$ into $T^c$ is a homotopy equivalence, we may choose generators for $Q_T$ so that $Z_{X}(X_+)$ and $Z_{X}(X_-)$ are respectively given by the spans
    $$
 \xymatrix{
& X \times X \ar[ld]_{\id} \ar[rd]^{f_+}& \\
X \times X  & & X \times X 
} \quad  \xymatrix{
& X \times X \ar[ld]_{\id} \ar[rd]^{f_-}& \\
X \times X & & X \times X 
}
$$
where $f_+,f_-: X \times X \longrightarrow X \times X$ are defined as
\begin{equation}
f_+(x,y) =(x \triangleright y, x),  \qquad f_-(x,y) =(y, x \triangleleft y). \label{eqn:Xlinearmap}
\end{equation}
\end{itemize}

Finally, let $B: \epsilon \to \epsilon$ be a braid, i.e.\ a tangle whose strands are never horizontal. In that case, the inclusion $\epsilon^c \hookrightarrow B^c$ is a homotopy equivalence. This implies that $\R_{X}(B) = \R_{X}(\epsilon) = X^{|\epsilon|}$ and, furthermore, the span $Z_X(B)$ can be taken to have the form
$$
\xymatrix{
    & X^{|\epsilon|} \ar[rd]^{f_B} \ar[ld]_{\id_{X^{|\epsilon|}}} & \\
    X^{|\epsilon|} & & X^{|\epsilon|}
}
$$
We will refer to the spans of this form as \emph{basic spans}, and denote them by $S_{f_B}$. Notice that basic spans are compatible with composition, since $S_{f_{B}} \circ S_{f_{B'}} = S_{f_{B} \circ f_{B'}}$. 

\section{$\SL_2$-representation varieties of knots} \label{sec:SL2repvarietyknots}

For the rest of this work, we will focus on the case of the group $G = \SL_2(\CC)$, simply denoted by $\SL_2$. With the natural choice of Borel subgroup given by upper triangular matrices, we have that the maximal torus is $T = \CC^*$ and the Weyl group is $W = \ZZ_2$ with action $\lambda \mapsto \lambda^{-1}$. Hence, $C_{\SL_2} = T/W = \CC$ and the quotient morphism $T = \CC^* \to T/W = \CC$ is $\lambda \mapsto \lambda + \lambda^{-1}$. This implies that the Chevalley map is given by the trace
$$
    \chi: \SL_2 \to C_{\SL_2} = \CC, \quad \chi(A) = \tr(A).
$$
For $t \in \CC-\{\pm 2\}$, the slice quandle $(\SL_2)_t \subseteq \Conj \SL_2 $ is
$$
    (\SL_2)_t = \{A \in \SL_2\,\mid\, \tr(A) = t\}.
$$
Writing $t = \lambda + \lambda^{-1}$, with $\lambda \neq \lambda^{-1}$, we have that $(\SL_2)_t$ is exactly equal to the conjugacy class of diagonalizable matrices $A \in \SL_2$ with eigenvalues $\lambda$ and $\lambda^{-1}$. The explicit identification of the quandle structure on $(\SL_2)_t$ for $t \in \CC -\{\pm 2\}$ is due to \cite{symplecticquandle2} and we describe it in Subsection \ref{subsec:generalizedsymplectic}.

Furthermore, if $t = 2$, the quandle $(\SL_2)_2$ corresponds to matrices $A \in \SL_2$ with a double eigenvalue $\lambda = 1$. Hence, $(\SL_2)_2$ is the union of two conjugacy classes, the identity $\{I\}$ and the matrices conjugate to the Jordan form
$$
    J_+ = \begin{pmatrix}
        1 & 1 \\ 0 & 1
    \end{pmatrix}.
$$
Analogously, $t = -2$ corresponds to matrices with double eigenvalue $\lambda = -1$ and $(\SL_2)_{-2}$ is the union of $\{-I\}$ with the conjugacy class of the corresponding Jordan form $J_-$. Notice that there is a quandle isomorphism $(\SL_2)_{2} \to (\SL_2)_{-2}$ given by $A \mapsto -A$. The description of the quandle structure on $(\SL_2)_2$ (and thus on $(\SL_2)_{-2})$ is due to \cite{symplecticquandle1} and we review it in Subsection \ref{subsec:parabolicquandle}.

\subsection{The parabolic symplectic quandle} \label{subsec:parabolicquandle}

In this section, we will describe explicitly the quandle structure of $(\SL_2)_2 \subseteq \Conj \SL_2$ by means of the so-called parabolic symplectic quandle.

\begin{defn}\label{def:parabolic-quandle}
    Let $Q_2 = \CC^2/\ZZ_2$, with action $x \mapsto -x$. The \emph{parabolic symplectic quandle}, or just \emph{parabolic quandle}, is the quandle structure on $Q_2$ given by
\begin{equation} \label{eqn:symplectictangleC2}
[x] \triangleright [y] = \left[y+\langle x,y \rangle x \right], \qquad [y] \triangleleft [x] = \left[y - \langle x,y \rangle x\right]
\end{equation}
for $[x], [y] \in \CC^2/\ZZ_2$, where $\langle x,y\rangle=\det (x,y)$ is the standard symplectic form on $\CC^2$.
\end{defn}

We rephrase and include for completeness the relation between the parabolic symplectic quandle and the set of parabolic elements in $\sldos$, as developed in \cite{symplecticquandle1}.

\begin{prop}\label{prop:iso-Q2}
    The symplectic quandle $Q_2$ is isomorphic to the slice quandle $(\SL_2)_2 \subseteq \Conj \SL_2$.
\end{prop}
\begin{proof}
Consider the map
\begin{equation} \label{eqn:C2paraboliciso}
\phi: Q_2 \longrightarrow \SL_2, \qquad [x] \mapsto P_x:=I+x\hat{x}, \qquad \textrm{where } x=\begin{pmatrix} x_1 \\ x_2 \end{pmatrix} \in \CC^2 \textrm{ and } \hat{x}=(-x_2,x_1).
\end{equation}

It is straightforward to check that this map is injective. Furthermore, its image is $(\SL_2)_2$. Indeed, $\tr(I+x\hat{x}) = 2$ for any $x \in \CC^2$. Reciprocally, any non-trivial parabolic matrix \(A\in (\SL_2)_2\) can be conjugated to the Jordan form $J_+$ and thus,
$$
A=\left(\begin{array}{ll}
a_{11} & a_{12} \\
a_{21} & a_{22}
\end{array}\right)\left(\begin{array}{ll}
1 & 1 \\
0 & 1
\end{array}\right)\left(\begin{array}{ll}
a_{11} & a_{12} \\
a_{21} & a_{22}
\end{array}\right)^{-1}=\left(\begin{array}{cc}
1-a_{11} a_{21} & a_{11}^2 \\
-a_{21}^2 & 1+a_{11} a_{21} 
\end{array}\right)= I+a\hat{a},
$$
where $a = (a_{11}, a_{21})^t \in \CC^2-\{(0,0)\}$. But, with the same identification, we have that the identity matrix corresponds exactly to $(0,0) \in \CC^2$. Finally, a direct calculation shows that if $A=\phi([a])$ and $B=\phi([b])$, for $[a],[b] \in Q_2$, then \(\phi([b] \triangleright [a])=BAB^{-1}\), showing that $\phi$ is a quandle isomorphism.
\end{proof}

Given a finitely presented quandle $Q$, the \emph{parabolic representation variety} is
$$
    \R_{Q_2}(Q) = \Hom_{\Qdl}(Q, Q_2) = \Hom_{\Qdl}(Q, (\SL_2)_2),
$$
where the last equality follows from the identification of Proposition \ref{prop:iso-Q2}. Furthermore, since $Q_2 = (\SL_2)_2$, the parabolic representation variety also inherits a conjugation action of $\SL_2$, and thus we can also consider the \emph{parabolic character variety}
$$
    \X_{Q_2}(Q) = \R_{Q_2}(Q) \sslash \SL_2.
$$
Actually, unraveling the isomorphism of Proposition \ref{prop:iso-Q2}, we observe that the action of $\SL_2$ on $Q_2 = \CC^2/\ZZ_2$ is $P \cdot [x] = [Px]$ for $P \in \SL_2$ and $x \in \CC^2$.  

In the particular case where $Q = Q_L$ is the link quandle of a generalized link $L$, we shall denote the corresponding varieties by $\R_{Q_2}(L)$ and $\X_{Q_2}(L)$.

\subsection{The generalized symplectic quandle} \label{subsec:generalizedsymplectic}

Although it will not be analyzed in this paper, the description of Section \ref{subsec:parabolicquandle} can be extended for the quandles $(\SL_2)_t$ for $t \neq \pm 2$, as done in \cite{symplecticquandle2}. 

Fix $M\in \CC-\{0,\pm 1\}$, and let us define
$$
\widetilde{\mathcal{D}}_M:=\left\{A \in \GL{2}(\mathbb{C}) \mid \operatorname{det}(A)=-M+M^{-1}\right\}
$$
Consider the equivalence relation $\sim$ on $\widetilde{\mathcal{D}}_M$ that declares that $A \sim B$ if and only if  $A^{-1}B\in \mathcal{D}$, where $\mathcal{D}$ is the set of diagonal elements in $\SL_2$.
\begin{defn}
The \emph{generalized symplectic quandle} is the quandle structure on \(\widetilde{\mathcal{D}}_M/\sim\) given by
\begin{equation} \label{eqn:symplectictangleC3}
[B] \triangleright [A] = \left[Ae_{M^{-1}}
- B\mathcal{D}(A,B)\right],
\end{equation}
where $e_M=\diag(M,M^{-1})$ and $\mathcal{D}(A,B)=\diag (\langle a_1,b_2 \rangle, \langle a_2,b_1 \rangle)$, where we write $A$ and $B$ as $A=(a_1\mid a_2)$ and $B=(b_1
 \mid b_2)$.
\end{defn}

Setting now $t=M+M^{-1} \neq \pm 2$, we have
$$
(\SL_2)_t= \lbrace A \in \SL_2 \mid \tr A = M+M^{-1} \rbrace,
$$
Any matrix $A\in (\SL_2 )_t$ can be conjugated to $J_M:=\left( \begin{smallmatrix}
    M & 1 \\ 0 & M^{-1}
\end{smallmatrix} \right)$, so that as sets, $(\SL_2 )_t\cong \SL_2/\Stab(J_M)$. Via this isomorphism, the following result is established in \cite{symplecticquandle2}.

\begin{prop}
The generalized symplectic quandle \(\widetilde{\mathcal{D}}_M/\sim\) is isomorphic to the slice quandle \( (\SL_2)_{t} \subseteq \Conj \SL_2\) for $t = M + M^{-1}$.
\end{prop}

\section{The Riley algebra and the determinant map} \label{sec:Rileyalgebra}

In \cite{RileyparabolicI}, Riley associated a polynomial $R(u) \in \ZZ[u]$ to study parabolic representations of a rational knot $K \subseteq S^3$. The construction works as follows. Let $K \subseteq S^3$ be a rational knot. Then Schubert showed that there is a presentation of the knot group $\Gamma_K = \pi_1(S^3-K)$ of $K$ with two meridian generators $x_1$ and $x_2$, and a single relation $r(x_1, x_2)$ that can be obtained from a suitable planar diagram of the knot. Therefore
$$
    \Gamma_K =\langle x_1, x_2 \,\mid\, r(x_1, x_2) = 1\rangle.
$$
Notice that, in particular, the knot quandle $Q_K$ admits a similar presentation. Hence, if one considers an irreducible parabolic representation $\rho: Q_K \to (\SL_2)_2$, using the $\SL_2$-action by conjugation, we can conjugate it to one of the form
\begin{equation} \label{eqn:Rileyconjugation}
\rho(x_1)=\begin{pmatrix} 1 & 1 \\ 0 & 1 \end{pmatrix}, \qquad \rho(x_2)=\begin{pmatrix} 1 & 0 \\ -u & 1 \end{pmatrix}.
\end{equation}
Riley showed that these two matrices satisfy $r(\rho(x_1), \rho(x_2)) = 1$, and thus define a parabolic representation of $Q_K$, if and only if $R(u) = 0$ for a certain polynomial $R(u) \in \CC[u]$, known as the \emph{Riley polynomial}. As mentioned in the introduction, this idea was further explored in \cite{symplecticquandle1} for kmot groups using quandles, where several additional properties of $R(u)$ were found. However, as noted in \cite{KitanoMorifuji}, the Riley polynomial is not a knot invariant: it depends on the chosen presentation of $\Gamma_K$.

\subsection{Parabolic character varieties}\label{sec:parabolic-character}

Let us provide a more intrinsic interpretation of these Riley invariants by means of parabolic character varieties, in the spirit of the classical result by Fricke-Vogt-Klein that states that the $\sldos$-character variety of the free group in two generators is isomorphic to $\CC^3$. Recall from Proposition \ref{prop:iso-Q2} that the slice quandle $(\SL_2)_2$ is isomorphic to $Q_2 = \CC^2/\ZZ_2$ so, we can consider the parabolic representation and character varieties
$$
    \R_{Q_2}(K) = \Hom_{\Qdl}(Q_K, Q_2), \qquad \X_{Q_2}(K) = \R_{Q_2}(K) \sslash \SL_2.
$$

Observe that $\R_{Q_2}(K)$ and $\X_{Q_2}(K)$ are knot invariants, since they only depend on $K$ through its knot quandle. Furthermore, since these varieties are affine, this GIT quotient can be explicitly described. To be precise, if $\R_{Q_2}(K) = \Spec A_K$, where $A_K$ is a finitely generated $\CC$-algebra (the algebra of regular functions on $\R_{Q_2}(K)$), then
$$
    \X_{Q_2}(K) = \Spec A_K^{\SL_2},
$$
where $A_K^{\SL_2} \subseteq A_K$ is the subalgebra of $\SL_2$-invariant functions, which can be interpreted as the algebra of regular functions on $\X_{Q_2}(K)$. Since $\SL_2$ is reductive, this algebra $A_K^{\SL_2}$ is a finitely generated $\CC$-algebra, which we will refer to as the \emph{Riley algebra} of $K$. Notice that $A_K^{\SL_2}$ is also a knot invariant.

Now, suppose that we provide a presentation of $\Gamma_K$, and thus of $Q_K$, with two generators and one relation, as in (\ref{eqn:Rileyconjugation}). This provides an inclusion $\R_{Q_2}(K) \subseteq \R_{Q_2}(\Free{2}) := \Hom_{\Qdl}(\Free{2}, Q_2)$, where $\Free{2}$ is the free quandle in two generators. Since $\R_{Q_2}(\Free{2}) = (\CC^2/\ZZ_2)^2 = \CC^4/(\ZZ_2 \times \ZZ_2)$, we get that
$$
    \R_{Q_2}(\Free{2}) = \Spec \CC[a_{11}, a_{12}, a_{21}, a_{22}]^{\ZZ_2 \times \ZZ_2},
$$
where the first $\ZZ_2$-factor acts by changing the sign in $a_{11}$ and $a_{12}$, whereas the second one acts on $a_{21}$ and $a_{22}$. Now, for the action of $\SL_2$ on $\R_{Q_2}(\Free{2})$, we have
$$
    \X_{Q_2}(\Free{2}) = \R_{Q_2}(\Free{2}) \sslash \SL_2 = \Spec \left(\CC[a_{ij}]^{\ZZ_2 \times \ZZ_2}\right)^{\SL_2}, \qquad 1 \leq i,j \leq 2.
$$
It can be proven (see Section \ref{sec:determinant-free} below) that $\left(\CC[a_{ij}]^{\ZZ_2 \times \ZZ_2}\right)^{\SL_2} \cong \CC[u]$ under the isomorphism $u \mapsto \left(a_{11}a_{22}-a_{12}a_{21}\right)^2$. In other words, we have an isomorphism
$$
    \X_{Q_2}(\Free{2}) \cong \mathbb{A}^1, \qquad \left[[a_1]_{\ZZ_2}, [a_2]_{\ZZ_2}\right]_{\SL_2} \mapsto u :=\det(a_1, a_2)^2,
$$
for $a_1, a_2 \in \CC^2$ with classes $[a_1]_{\ZZ_2}, [a_2]_{\ZZ_2} \in \CC^2/\ZZ_2$. In this manner, the inclusion $\X_{Q_2}(K) \subseteq \X_{Q_2}(\Free{2}) = \mathbb{A}^1 = \Spec \CC[u]$ as a closed subscheme implies that
$$
    \X_{Q_2}(K) \cong \Spec \CC[u]/(P(u))
$$
for some polynomial $P(u) \in \CC[u]$. In other words, $A_K^{\SL_2} \cong \CC[u]/(P(u))$. Notice that this polynomial satisfies $P(0) = 0$, since $u = 0$ corresponds to $a_1$ and $a_2$ being linearly dependent and, thus, to the class of the trivial representation. Following the definition by Riley in \cite{RileyparabolicI}, we see that (up to re-scaling) the original Riley polynomial of $K$ associated to the given presentation of the knot group of $K$ is $R(u) = P(u)/u$, after removing the trivial character corresponding to $u = 0$. 

\subsection{The determinant map of a free quandle}\label{sec:determinant-free}

The construction of Section \ref{sec:parabolic-character} can be much better understood in terms of a covering. To be precise, let $\widetilde{Q}_2 = \CC^2$ with quandle structure
\begin{equation*}
x \triangleright y = y + \langle x,y \rangle x, \qquad y \triangleleft x = y - \langle x,y \rangle x,
\end{equation*}
for $x, y \in \CC^2$, where $\langle x,y\rangle=\det (x,y)$ is the standard symplectic form on $\CC^2$. The quandle $(\widetilde{Q}_2, \triangleright)$ will be referred to as the \emph{semiparabolic quandle}. 

By its very description, we have a quandle homomorphism of index two onto the parabolic quandle of Definition \ref{def:parabolic-quandle},
$$
    \widetilde{Q}_2 \longrightarrow Q_2, \qquad x \mapsto [x].
$$
In particular, for a quandle $Q$, we get a morphism
$$
    \R_{\widetilde{Q}_2}(Q) = \Hom_{\Qdl}(Q, \widetilde{Q}_2) \longrightarrow \R_{Q_2}(Q) = \Hom_{\Qdl}(Q, Q_2),
$$
which we will refer to as the \emph{determinant map}, whereas the variety $\R_{\widetilde{Q}_2}(Q)$ will be the \emph{semiparabolic representation variety}. Notice that the action of $\SL_2$ on $Q_2$ lifts to an action on $\widetilde{Q}_2$ by $P \cdot x = Px$, for $P \in \SL_2$ and $x \in \widetilde{Q}_2 = \CC^2$, so we can also consider the \emph{semiparabolic character variety}
$$
    \X_{\widetilde{Q}_2}(Q) = \R_{\widetilde{Q}_2}(Q) \sslash \SL_2.
$$

The reason for the name ``determinant map'' comes from the following observation. Let $\Free{n}$ be the free quandle in $n$ generators, and consider the induced morphism $\R_{\widetilde{Q}_2}(\Free{n}) \to \R_{Q_2}(\Free{n})$, which gives rise to a morphism
$$
    \X_{\widetilde{Q}_2}(\Free{n}) = \R_{\widetilde{Q}_2}(\Free{n}) \sslash \SL_2 \longrightarrow \X_{Q_2}(\Free{n}) = \R_{Q_2}(\Free{n}) \sslash \SL_2.
$$

In order to understand this morphism, consider the complex polynomial algebra
$$
    \widetilde{A}_n = \CC[a_{1}^1, a_1^2, \ldots, a_n^2].
$$
It is the algebra of regular functions on $\R_{\widetilde{Q}_2}(\Free{n}) = (\CC^2)^n$, where a point is $(a_1, \ldots, a_n) \in \R_{\widetilde{Q}_2}(\Free{n})$ with $a_k = (a_k^1, a_k^2) \in \CC^2$.

Furthermore, the action of $\SL_2$ on $\R_{\widetilde{Q}_2}(\Free{n})$ induces an action on $\widetilde{A}_n$ by $P \cdot a_k = P^{-1}a_k$ for $P \in \SL_2$ and $a_k \in \CC^2$, and the ring of regular functions on $\X_{\widetilde{Q}_2}(\Free{n})$ is the subring of $\SL_2$-invariant functions $\widetilde{A}_n^{\SL_2}$. This ring is well known in the literature through the so-called first fundamental theorem of $\SL_2$.

\begin{prop}[{\cite[Chapter 2]{dolgachev2003lectures}}]\label{prop:tildeQ2-SL2}
We have an identification
$$
\widetilde{A}_n^{\SL_2} = \frac{\CC[s_{ij}]}{(s_{ij}s_{k\ell}-s_{ik}s_{j\ell}+s_{i\ell}s_{jk})}, \quad \textrm{for } 1 \leq i < j <k < \ell \leq n.
$$
The inclusion morphism $\widetilde{A}_n^{\SL_2} \to \widetilde{A}_n$ is given by $s_{ij} \mapsto \langle a_i, a_j\rangle = a_i^1a_j^2-a_i^2a_j^1$.
\end{prop}

From this result, we directly get a description of the corresponding free character variety.

\begin{cor}\label{cor:iso-tildQ2}
    We have an isomorphism
    $$
        \X_{\widetilde{Q}_2}(\Free{n}) \cong \Spec \frac{\CC[s_{ij}]}{(s_{ij}s_{k\ell}-s_{ik}s_{j\ell}+s_{i\ell}s_{jk})}, \quad \textrm{for } 1 \leq i < j <k < \ell \leq n.
    $$
    The isomorphism is explicitly given by $[a_1, \ldots, a_n] \in \X_{\widetilde{Q}_2}(\Free{n}) \mapsto \left(s_{ij} = \langle a_i, a_j \rangle \right)_{i,j}$.
\end{cor}

In the following, let us denote by
\begin{equation}\label{eq:Bntilde}
    \widetilde{B}_n = \widetilde{A}_n^{\SL_2}= \frac{\CC[s_{ij}]}{(s_{ij}s_{k\ell}-s_{ik}s_{j\ell}+s_{i\ell}s_{jk})} \quad \textrm{for } 1 \leq i < j <k < \ell \leq n.
\end{equation}
the algebra of Proposition \ref{prop:tildeQ2-SL2}, so $\X_{\widetilde{Q}_2}(\Free{n}) = \Spec \widetilde{B}_n$.

\begin{rem}\label{rem:redundant-Bntilde}
    Notice that, in the description of $\widetilde{B}_n$ of Proposition \ref{prop:tildeQ2-SL2}, this is generated by variables $s_{ij}$ with $i < j$. However, often it is also useful to consider the auxiliary variables $s_{ij} = -s_{ji}$ for $i \geq j$ (in particular, $s_{ii} = 0$). In that case, we can rewrite (\ref{eq:Bntilde}) as
$$
    \widetilde{B}_n = \frac{\CC[s_{ij}]}{(s_{ij}s_{k\ell}-s_{ik}s_{j\ell}+s_{i\ell}s_{jk}, s_{ij} + s_{ji})} \quad \textrm{for } 1 \leq i, j, k, \ell \leq n.
$$
\end{rem}

\begin{ex}\label{ex:B2}
    In the case $n = 2$, we have $\widetilde{B}_2 = \CC[s]$, where $s = s_{12}$ and with no Pl\"ucker relations holding. In other words, $\X_{\widetilde{Q}_2}(\Free{2}) = \mathbb{A}^1$ is an affine line.
\end{ex}

\begin{rem}
    Corollary \ref{cor:iso-tildQ2} allows us to perform local calculations in $\X_{\widetilde{Q}_2}(\Free{n})$. Indeed, under the isomorphism of Corollary \ref{cor:iso-tildQ2}, we can consider the open sets $D(s_{pq}) \subseteq \X_{\widetilde{Q}_2}(\Free{n})$, corresponding to the points with $s_{pq} \neq 0$, for fixed $1 \leq p < q \leq n$. In that case, from the Pl\"ucker relations, we have for $(i,j) \neq (p,q)$,
    $$
        s_{ij}	= 	\frac{s_{pi}s_{qj}-s_{pj}s_{qi}}{s_{pq}}.
    $$
    This implies that the localization is $(\widetilde{B}_n)_{s_{pq}} = \CC\left[
	s_{pq}^{\pm1}, s_{pi},s_{qi}\right]$ for $i\neq p,q$. Consequently, $D(s_{pq}) = \left(\mathbb{A}^1-\{0\}\right) \times \mathbb{A}^{2n-4}$.
\end{rem}

Now, observe that $\widetilde{A}_n$ also admits an action of $\ZZ_2^n$ by change of sign, given by
$$
    (\epsilon_1, \ldots, \epsilon_n) \cdot a_k = \epsilon_k a_k,
$$
for $\epsilon_1, \ldots, \epsilon_n \in \{\pm1\} = \ZZ_2$ and $a_k \in \CC^2$ with $1 \leq k \leq n$. Under the inclusion $\widetilde{B}_n \to \widetilde{A}_n$, this action restricts to the action on $\widetilde{B}_n$ given, for $\epsilon_1, \ldots, \epsilon_n \in \{\pm1\} = \ZZ_2$ and $1 \leq i, j \leq n$, by
$$
    (\epsilon_1, \ldots, \epsilon_n) \cdot s_{ij} = \epsilon_i\epsilon_j s_{ij}.
$$
In this direction, we have the following result.

\begin{prop}[{\cite[Theorem 2.6]{drensky2003defining}}]\label{prop:Q2-SL2}
We have an identification
$$
\widetilde{B}_n^{\ZZ_2^n} =\CC[u_{ij}, t_{ijk}]/H_n, \quad \textrm{for } 1 \leq i < j < k \leq n,
    $$
for a certain ideal $H_n$ of $\CC[u_{ij}, t_{ijk}]$. The inclusion morphism $\widetilde{B}_n^{\ZZ_2^n} \to \widetilde{B}_n$ is given by $u_{ij} \mapsto s_{ij}^2$ and $t_{ijk} \mapsto s_{ij}s_{jk}s_{ik}$.
\end{prop}

\begin{cor}\label{cor:iso-Q2}
We have an isomorphism
    $$
        \X_{Q_2}(\Free{n}) \cong \Spec \CC[u_{ij}, t_{ijk}]/H_n, \quad \textrm{for } 1 \leq i < j < k \leq n,
    $$
for a certain ideal $H_n$ of $\CC[u_{ij}, t_{ijk}]$. Under this isomorphism, the map $\X_{\widetilde{Q}_2}(\Free{n}) \to \X_{Q_2}(\Free{n})$ becomes $(s_{ij}) \mapsto (u_{ij} = s_{ij}^2, t_{ijk} = s_{ij}s_{jk}s_{ik})$.
\end{cor}

\begin{proof}
    We have $\X_{Q_2}(\Free{n}) = (\CC^2/\ZZ_2)^n \sslash \SL_2 = \left((\CC^2)^{n} /\ZZ_2^n \right) \sslash \SL_2$, where each factor of $\ZZ_2^n$ acts on $(\CC^2)^{n}$ by changing the sign of the corresponding vector. However, both actions commute, so by Proposition \ref{prop:tildeQ2-SL2}, we have
    $$
        \X_{Q_2}(\Free{n}) = \left((\CC^2)^{n}\sslash \SL_2\right) /\ZZ_2^n = \left(\Spec \widetilde{B}_n \right)/\ZZ_2^n = \Spec \widetilde{B}_n^{\ZZ_2^n},
    $$
    and the result follows from Proposition \ref{prop:Q2-SL2}.
\end{proof}

\begin{rem}\label{rmk:ideal-Hn}
    The ideal $H_n$ of $\CC[u_{ij}, t_{ijk}]$ appearing in Proposition \ref{prop:Q2-SL2} is also explicitly described in \cite[Theorem 2.6]{drensky2003defining}. To be precise, it is the ideal generated by $u_{ij} - u_{ji} = 0$ and $u_{ii} = 0$ for $1 \leq i < j \leq n$, as well as
    \[
2t_{ijk}t_{i'j'k'}
-
\det\begin{pmatrix}
u_{ii'}&u_{ij'}&u_{ik'}\\
u_{ji'}&u_{jj'}&u_{jk'}\\
u_{ki'}&u_{kj'}&u_{kk'}
\end{pmatrix} \qquad \textrm{for } 1 \leq i < j < k \leq n \textrm{ and } 1 \leq i' < j' < k '\leq n,
\]
as well as
\[
u_{ai}t_{jk\ell}
-u_{aj}t_{ik\ell}
+u_{ak}t_{ij\ell}
-u_{a\ell}t_{ijk}, \qquad \textrm{for } 1 \leq i < j < k < \ell \leq n \textrm{ and } 1 \leq a \leq n.
\]
For instance, we have that $H_3 = (t_{123}^2-u_{12}u_{13}u_{23})$.
\end{rem}

As for $\X_{\widetilde{Q}_2}(\Free{n})$, in the following we shall denote by $B_n = \CC[u_{ij}, t_{ijk}]/H_n$ the algebra of regular functions on the quandle character variety $\X_{Q_2}(\Free{n})$.

Finally, we can also consider the action of $\ZZ_2^n$ directly on $\widetilde{A}_n = \CC[a_1^1, \ldots, a_n^2]$ by change of sign. In that case, we have the following result.

\begin{prop}\label{prop:Q2tilde-sign}
    We have an identification
$$
    \widetilde{A}_n^{\ZZ_2^n} = \frac{\CC[x_k, y_k, z_k]}{(x_kz_k-y_k^2)}, \qquad \textrm{for } 1 \leq k \leq n. 
$$
The inclusion morphism $\widetilde{A}_n^{\ZZ_2^n} \to \widetilde{A}_n$ is given by $x_k \mapsto (a_k^1)^2$, $y_k \mapsto a_k^1a_k^2$ and $z_k \mapsto (a_k^2)^2$.
\end{prop}

\begin{proof}
Notice that, we can decompose $\widetilde{A}_n$ as a tensor product of rings $\widetilde{A}_n = \CC[a^1, a^2]^{\otimes n}$, and each factor of $\ZZ_2^n$ acts on a factor of $\widetilde{A}_n$. In this manner, we have
$$
    \widetilde{A}_n^{\ZZ_2^n} = \left(\CC[a^1, a^2]^{\ZZ_2}\right)^{\otimes n}.
$$
Now, observe that the generators for the invariant monomials of $\CC[a^1, a^2]$ under the $\ZZ_2$-action are the products $(a^1)^2, (a^2)^2$ and $a^1a^2$. Therefore, $\CC[a^1, a^2]^{\ZZ_2} = \CC[(a^1)^2, (a^2)^2, a^1a^2]$ and the result follows by setting $x = (a^1)^2, z = (a^2)^2$, and $y =a^1a^2$, which are constrained to the relation $xz = y^2$.
\end{proof}

\begin{rem}\label{rmk:Q2tilde-sign-matrix}
There is an alternative way of describing the ring $\CC[a^1, a^2]^{\ZZ_2} = \CC[x, y, z]/(xz-y^2)$ of Proposition \ref{prop:Q2tilde-sign}. Indeed, if we consider the variety $M_0$ of complex symmetric $2 \times 2$ matrices of zero determinant, then $\CC[x, y, z]/(xz-y^2)$ is exactly its algebra of regular functions. In other words, the morphism $\widetilde{Q}_2 = \CC^2 = \Spec \CC[a^1, a^2] \to M_0$ given by
$$
    \widetilde{Q}_2 \ni a = \begin{pmatrix}
        a^1 \\
        a^2
    \end{pmatrix} \mapsto aa^t \in M_0
$$
is $\ZZ_2$-invariant and descends to an isomorphism $Q_2 = \widetilde{Q}_2 / \ZZ_2 = \Spec \CC[a^1, a^2]^{\ZZ_2} \stackrel{\cong}{\to} M_0$.
\end{rem}

As for the other rings considered in this section, we shall denote
$$
    A_n = \widetilde{A}_n^{\ZZ_2^n} = \frac{\CC[x_k, y_k, z_k]}{(x_kz_k-y_k^2)}, \qquad \textrm{for } 1 \leq k \leq n.
$$
Notice that this is the ring of regular functions on the representation variety $\R_{Q_2}(\Free{n})$. Therefore, we have a pair of dual commutative diagrams
\begin{equation}\label{eq:diagram-free-rings}
\begin{split}
\xymatrix{
\R_{\widetilde{Q}_2}(\Free{n}) \ar[r]\ar[d] & \R_{Q_2}(\Free{n}) \ar[d] \\
\X_{\widetilde{Q}_2}(\Free{n}) \ar[r] & \X_{Q_2}(\Free{n})
}\qquad\qquad \xymatrix{
\widetilde{A}_n  & A_n = \widetilde{A}_n^{\ZZ_2^n} \ar[l]  \\
 \widetilde{B}_n \ar[u] & B_n = \widetilde{B}_n^{\ZZ_2^n} = A_n^{\SL_2} \ar[l]  \ar[u]
}
\end{split}
\end{equation}

\subsection{The determinant map of a knot quandle}\label{sec:determinant-knot}

Now, let us come back to the case of a knot $K \subseteq S^3$, for which we have a regular morphism
\begin{equation}\label{eq:determinant-cover}
    \X_{\widetilde{Q}_2}(K) \longrightarrow \X_{Q_2}(K).
\end{equation}

Denote by $\widetilde{A}_K$ the ring of regular functions on the representation variety $\R_{\widetilde{Q}_2}(K)$ and, as in Section \ref{sec:determinant-free}, by $A_K$ the ring of regular functions on $\R_{Q_2}(K)$, which are naturally equipped with a ring homomorphism $A_K \to \widetilde{A}_K$ coming from the morphism $\R_{\widetilde{Q}_2}(K) \to \R_{Q_2}(K)$. The group $\SL_2$ naturally acts on $\widetilde{A}_K$ and $A_K$, so that $\widetilde{B}_K := \widetilde{A}_K^{\SL_2}$ and $B_K := A_K^{\SL_2}$ are the rings of regular functions on $\X_{\widetilde{Q}_2}(K)$ and $\X_{Q_2}(K)$, respectively. Notice that $B_K$ is exactly the Riley algebra introduced in Section \ref{sec:parabolic-character}.

Furthermore, if we consider a presentation of the knot group $\Gamma_K$ with $n$ generators, it provides embeddings $\R_{\widetilde{Q}_2}(K) \subseteq \R_{\widetilde{Q}_2}(\Free{n})$ and $\R_{Q_2}(K) \subseteq \R_{Q_2}(\Free{n})$ as closed subvarieties. Therefore, we get commutative diagrams
\begin{equation}\label{eq:diagram-free}
\begin{split}
\xymatrix{
\X_{\widetilde{Q}_2}(K) \ar[r]\ar@{^(-{>}}[d] & \X_{Q_2}(K) \ar@{^(-{>}}[d] \\
\X_{\widetilde{Q}_2}(\Free{n})  \ar[r] & \X_{Q_2}(\Free{n})
}\qquad\qquad \xymatrix{
\widetilde{B}_K = \widetilde{A}_K^{\SL_2}  & B_K = A_K^{\SL_2} \ar[l]  \\
 \widetilde{B}_n \ar@{-{>>}}[u] & B_n \ar[l]  \ar@{-{>>}}[u]
}
\end{split}
\end{equation}
where the horizontal arrows are the determinant maps and the right bottom one is given by $u_{ij} \mapsto s_{ij}^2$ and $t_{ijk} \mapsto s_{ij}s_{jk}s_{ik}$. 

This picture is particularly simple in the case $n = 2$, for instance, if $K$ is a rational knot and we are considering its Schubert presentation with two generators and one relation. In that case, by Example \ref{ex:B2} we have $\widetilde{B}_2 = \CC[s]$ and $\X_{\widetilde{Q}_2}(\Free{2}) = \mathbb{A}^1$. Moreover, by Corollary \ref{cor:iso-Q2}, $B_2 = \CC[u]$ and the morphism $B_2 \to \widetilde{B}_2$ is given by $u \mapsto s^2$. Equivalently, $\X_{{Q}_2}(\Free{2}) = \mathbb{A}^1$ and we have the double ramified cover
$$
    \X_{\widetilde{Q}_2}(\Free{2}) = \mathbb{A}^1 \to \X_{{Q}_2}(\Free{2}) = \mathbb{A}^1, \qquad s \mapsto s^2.
$$
Actually, the ramification point $s = 0$ corresponds precisely to the identity representation.

Furthermore, in this picture, $\X_{\widetilde{Q}_2}(K) \subseteq \X_{\widetilde{Q}_2}(\Free{2}) = \mathbb{A}^1$ and $\X_{{Q}_2}(K) \subseteq \X_{{Q}_2}(\Free{2}) = \mathbb{A}^1$ are zero-dimensional closed subschemes fitting in the diagram (\ref{eq:diagram-free}).This implies that
$$
    \X_{{Q}_2}(K) = \Spec \CC[u]/(P(u)), \qquad \X_{\widetilde{Q}_2}(K) = \Spec \CC[s]/(\widetilde{P}(s)),
$$
where $R(u) = P(u)/u \in \CC[u]$ is the Riley polynomial of the presentation and $\widetilde{P}(s) \in \CC[s]$.

\begin{rem}
    As a consequence of this discussion, if a knot $K$ has associated algebra $B_K = A_K^{\SL_2}$ with Krull dimension $\dim B_K > 0$, then $K$ is not a rational knot. For non-rational knots, the study of those knots whose parabolic character variety has components of positive dimension is a very active area of research, see for instance \cite{ChoiPortiYoon,ChenHigh}.
\end{rem}

\section{Parabolic representations of pretzel knots} \label{sec:parabolicpretzelknots}
In this section, we apply the functors $Z_{Q_2}$ and $Z_{\widetilde{Q}_2}$, developed in Sections \ref{sec:quandleTQFT} and \ref{sec:SL2repvarietyknots}, to a particular family of knots, the so-called pretzel knots. Recall that, by Proposition \ref{prop:iso-Q2}, the quandle $Q_2$ is isomorphic to the slice quandle $(\SL_2)_2 \subseteq \Conj \SL_2$, so we will compute (semi)parabolic representation varieties.

\begin{figure}[h]
  \centering
  \includegraphics[width=0.31\textwidth]{pretzel1.pdf}\hfill
  \includegraphics[width=0.31\textwidth]{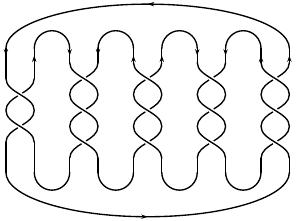}\hfill
  \includegraphics[width=0.31\textwidth]{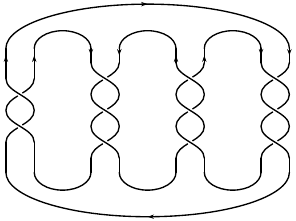}
  \caption{Oriented odd and even pretzel knots.}
  \label{fig:pretzelknot}
\end{figure}

Let us fix integers $q_1, \ldots, q_n \in \ZZ$. The pretzel knot or link $P(q_1,\ldots,q_n)$ is formed by cyclically connecting $n$ twisted bands, where each band is formed by $q_i$ half twists, as displayed in Figure \ref{fig:pretzelknot}. The sign of $q_i$ determines if the twist is positive or negative, and it is well-known that the construction produces a knot if and only if one of the following conditions holds:
\begin{itemize}
    \item Exactly one $q_i$ is even, being all other $q_j$ odd.
    \item The number of bands $n$ and all $q_1, \ldots, q_n$ are odd.
\end{itemize}
The isotopy class of any pretzel knot is invariant under cyclic permutation of $(q_1,\ldots,q_n)$, so we will assume that $q_1$ is even in the first case, and refer to this subfamily as \emph{even}, and consequently we will call it \emph{odd} in the second case. 

\begin{figure}[h]
\includegraphics[width=9cm]{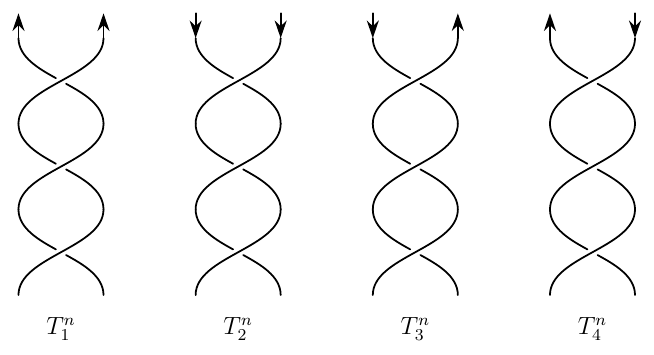}
\caption{Basic tangles appearing in a pretzel knot} \label{fig:twistedbraidings}
\end{figure}

After orienting the pretzel knot, each band in the construction is a $2$-braid that can have either both strands oriented in the same direction, or in opposite directions. In the former situation, we will say that the braid is \emph{codirected}, and \emph{alternating} in the latter. We will denote these tangles by $T_1^n$ and $T_2^n$ for the codirected ones, and $T_3^n$ and $T_4^n$ for the alternating tangles, as shown in Figure \ref{fig:twistedbraidings}. Here, $n$ denotes the number of twists in the tangle.

In this regard, in an oriented pretzel knot, all the twists of an even pretzel knot are codirected (except for $q_1$ when $n$ is odd) and every twist of an odd pretzel knot is alternating \cite{PretzelAlex2}. 

Recall from \eqref{eqn:symplectictangleC2} that $Q_2=\CC^2/\ZZ_2$ with the quandle operation $[x]\triangleright [y]=\left[y+\langle x,y\rangle x\right]$, where $\langle x,y\rangle = \det(x, y)$ is the standard symplectic form in $\CC^2$. However, for our purposes, it will be easier to consider the double cover $\widetilde{Q}_2 \to Q_2$, given by
$$
\widetilde{Q}_2=\CC^2, \quad x\triangleright y=y+\langle x,y\rangle x.
$$
As shown in (\ref{eqn:Xlinearmap}), the spans associated under $Z_{\widetilde{Q}_2}$ of the tangles $X_+, X_-: (+,+) \to (+,+)$ from Figure \ref{fig:tangles} are respectively given by the maps
\begin{align*}
f_+: \widetilde{Q}_2\times \widetilde{Q}_2 \rightarrow \widetilde{Q}_2\times \widetilde{Q}_2, & & f_+(a,b)&=\left(b+\langle a,b\rangle a,a\right),  \\
f_-: \widetilde{Q}_2\times \widetilde{Q}_2 \rightarrow \widetilde{Q}_2\times \widetilde{Q}_2, & & f_-(a,b)&=\left(b,a + \langle a,b\rangle b\right).
\end{align*}
A direct calculation shows that $f_+$ and $f_-$ are inverses of each other, due to the alternating bilinearity of the symplectic form. Notice that the spans for the original quandle $Q_2$ can be recovered from the spans for $\widetilde{Q}_2$ by taking equivalence classes.

Finally, notice that all the tangles $T_i^n$ for $1 \leq i \leq 4$ and $n \in \ZZ$ are braids from two points to two points. Therefore, the image of $T_i^n$ under the functor $Z_{\widetilde{Q}_2}$ will be a basic span $S_{f_{T_{i}^n}}$ of the form
$$
    \widetilde{Q}_2 \times \widetilde{Q}_2 \xleftarrow{\Id_{\widetilde{Q}_2 \times \widetilde{Q}_2}} \widetilde{Q}_2 \times \widetilde{Q}_2 \xrightarrow{f_{T_i^n}} \widetilde{Q}_2 \times \widetilde{Q}_2,
$$
where $f_{T_i^n}: \widetilde{Q}_2 \times \widetilde{Q}_2 \to \widetilde{Q}_2 \times \widetilde{Q}_2$ is a regular morphism. 

\begin{figure}[htbp]
    \centering
    \begin{minipage}[c]{0.47\linewidth}
        \centering
        \includegraphics[width=\linewidth]{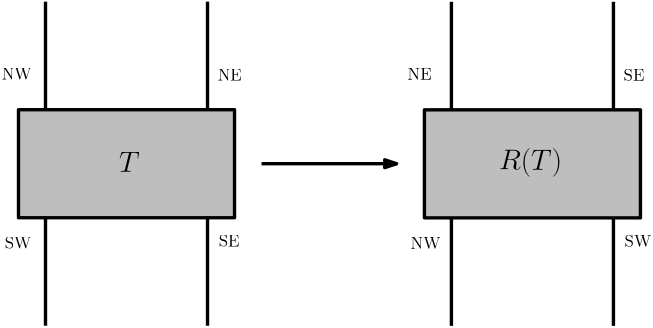}
    \end{minipage}%
    \hfill
    \begin{minipage}[c]{0.35\linewidth}
        \centering
        \includegraphics[width=\linewidth]{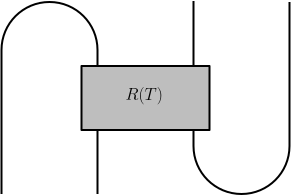}
    \end{minipage}

    \caption{A 2-tangle $T$ and its rotation by $\pi/2$, decomposed in terms of elementary tangles}
    \label{fig:tangle-rotation}
\end{figure}

In the upcoming calculations, it will be very useful to consider the \emph{rotation of a $2$-tangle}. This is an operation that takes a tangle $T$ between two points, and applies a counterclockwise rotation of $\pi/2$ radians  to it, as shown in Figure \ref{fig:tangle-rotation}.

\begin{lem}\label{lem:rotation}
    Let $T$ be a $2$-tangle, and let $Z_{\widetilde{Q}_2}(T)$ be the span
    $$
        \widetilde{Q}_2 \times \widetilde{Q}_2 \xleftarrow{(f_1, f_2)} \R_{\widetilde{Q}_2}(T) \xrightarrow{(g_1, g_2)}\widetilde{Q}_2 \times \widetilde{Q}_2,
    $$
    for some functions $f_1, f_2, g_1, g_2: \R_{\widetilde{Q}_2}(T) \to \widetilde{Q}_2$. Then, the rotated tangle $R(T)$ has associated span $Z_{\widetilde{Q}_2}(R(T))$ given by
    $$
        \widetilde{Q}_2 \times \widetilde{Q}_2 \xleftarrow{(g_1, f_1)} \R_{\widetilde{Q}_2}(T) \xrightarrow{(g_2, f_2)}\widetilde{Q}_2 \times \widetilde{Q}_2.
    $$
\end{lem}

\begin{proof}
The proof is just a calculation. As shown in Figure \ref{fig:tangle-rotation}, we can decompose $R(T)$ as
$$
R(T)= (\curvearrowleft\otimes \uparrow\otimes \uparrow)\circ (\downarrow\otimes T \otimes \uparrow)\circ (\downarrow\otimes \uparrow\otimes \upcurvearrowright)
$$
Now, observe that $Z_{\widetilde{Q}_2}(\downarrow\otimes \uparrow\otimes \upcurvearrowright)$ equals the span
$$
    \widetilde{Q}_2^2 \xleftarrow{\, (a,b,c)\mapsto(a,b) \,} \widetilde{Q}_2^3 \xrightarrow{\, (a,b,c)\mapsto(a,b,c,c) \,} \widetilde{Q}_2^4.
$$
whereas $Z_{\widetilde{Q}_2}(\curvearrowleft\otimes \uparrow\otimes \uparrow)$ is the span
$$
\widetilde{Q}_2^4 \xleftarrow{\, (a,b,c)\mapsto(a,a,b,c)  \,} \widetilde{Q}_2^3 \xrightarrow{\, (a,b,c)\mapsto(b,c)  \,} \widetilde{Q}_2^2.
$$

Therefore, applying $Z_{\widetilde{Q}_2}$ to the three factors and composing the resulting spans, we obtain the following diagram, all of whose squares are cartesian
\begingroup
\small
\[
\begin{tikzcd}[
  column sep=2.8em,
  row sep=5.0em,
  cells={nodes={inner sep=1.2pt}},
  every label/.append style={font=\scriptsize}
]
& & &
\R_{\widetilde{Q}_2}(T)
  \arrow[dl, "{\rho\mapsto(g_1(\rho),\rho)}"', pos=.52]
  \arrow[dr, "{\rho\mapsto(\rho,f_2(\rho))}", pos=.52]
& & & \\
& &
\widetilde{Q}_2\times\R_{\widetilde{Q}_2}(T)
  \arrow[dl,
    "{\substack{(a,\rho)\mapsto\\
    (a,f_1(\rho),f_2(\rho))}}"',
    pos=.50]
  \arrow[dr,
    "{\substack{(a,\rho)\mapsto\\
    (a,\rho,f_2(\rho))}}",
    pos=.28]
& &
\R_{\widetilde{Q}_2}(T)\times\widetilde{Q}_2
  \arrow[dl,
    "{\substack{(\rho,b)\mapsto\\
    (g_1(\rho),\rho,b)}}"',
    pos=.28]
  \arrow[dr,
    "{\substack{(\rho,b)\mapsto\\
    (g_1(\rho),g_2(\rho),b)}}",
    pos=.50]
& & \\
&
\widetilde{Q}_2^3
  \arrow[dl]
  \arrow[dr]
& &
\widetilde{Q}_2\times
\R_{\widetilde{Q}_2}(T)\times
\widetilde{Q}_2
  \arrow[dl,
    "{\substack{(x,\rho,y)\mapsto\\
    (x,f_1(\rho),f_2(\rho),y)}}"',
    pos=.49]
  \arrow[dr,
    "{\substack{(x,\rho,y)\mapsto\\
    (x,g_1(\rho),g_2(\rho),y)}}",
    pos=.49]
& &
\widetilde{Q}_2^3
  \arrow[dl]
  \arrow[dr]
& \\
|[xshift=-1.8em]|\widetilde{Q}_2^2
& &
\widetilde{Q}_2^4
& &
\widetilde{Q}_2^4
& &
|[xshift=1.8em]|\widetilde{Q}_2^2
\end{tikzcd}
\]
\endgroup
In particular, composing the morphisms from the apex of the diagram to the left-most and right-most legs, we obtain the desired span.
\end{proof}

\subsection{Codirected braidings}
\label{sec:codirected-braid}

Let us begin with the upwards codirected braiding. For $n \geq 0$, we have $T_1^n = X_+^n: (+,+) \to (+,+)$. By functoriality of the functor $Z_{\widetilde{Q}_2}: \Tang \to \Span \Varc$, we have that the span $Z_{\widetilde{Q}_2}(X_+^n): \widetilde{Q}_2 \times \widetilde{Q}_2 \to \widetilde{Q}_2 \times \widetilde{Q}_2$ is equivalent to the basic span $S_{f_+^n}$. In order to describe it explicitly, given $s\in \CC$, set
$$
F_+(s):=
\begin{pmatrix}
s&1\\
1&0
\end{pmatrix},
$$
and, more generally for $n\geq 1$, we define the product
$$
P^n_+(s):=F_+\bigl((-1)^{n-1}s\bigr) \cdots F_+(-s)F_+(s), \qquad P_+^0(s):=\id_{\CC^2}.
$$
\begin{prop}  \label{prop:closedformspowersf_-}
Let $a_0, b_0 \in \widetilde{Q}_2$, and denote $(a_n,b_n):=f_+^n(a_0,b_0)$. Then, we have
\begin{equation} \label{eqn:f-minus-n-matrixexpression}
\begin{pmatrix}
a_n\\ b_n
\end{pmatrix}
=
P^n_+(s_0)
\begin{pmatrix}
a_0\\ b_0
\end{pmatrix},
\end{equation}
where $s_0:= \langle a_0,b_0 \rangle$. In other words, for $n \geq 0$, $Z_{\widetilde{Q}_2}(T_1^n): \widetilde{Q}_2 \times \widetilde{Q}_2 \to \widetilde{Q}_2 \times \widetilde{Q}_2$ is the basic span associated to the function $f_+^n(a,b) = P^n_+(s_0)(a, b)$.
\end{prop}
\begin{proof}
From the very definition of $f_+$ we have $(a_{n+1},b_{n+1})=f_+(a_n,b_n)=(b_n+\langle a_n,b_n\rangle a_n, a_n)$. In particular, if we set $s_n:=\langle a_n,b_{n}\rangle$, the previous recursion implies that
$$
    s_{n+1}= \langle a_{n+1},b_{n+1}\rangle = \langle b_n+\langle a_n,b_n\rangle a_n, a_n \rangle = \langle b_n, a_n \rangle = -s_n.
$$
In this manner, $s_n = (-1)^ns_0$ and thus we have
$$
\begin{pmatrix} a_{n+1} \\ b_{n+1} \end{pmatrix} = F_+(s_n) \begin{pmatrix} a_n \\ b_n \end{pmatrix} = F_+\bigl((-1)^ns_0\bigr) \begin{pmatrix} a_n \\ b_n \end{pmatrix}.
$$
Iterating this recurrence we obtain \eqref{eqn:f-minus-n-matrixexpression}.
\end{proof}

A completely analogous argument gives the corresponding closed-form expressions for the braid with negative twisting $T_1^{-n} = X_-^n: (+,+) \to (+,+)$ for $n \geq 1$. Again, the span $Z_{\widetilde{Q}_2}(X_-^n)$ is equivalent to the basic span $S_{f_{-}^n}$. To identify this morphism, set
$$
F_-(s):=
\begin{pmatrix}
0&1\\
1&s
\end{pmatrix},
$$
and, for $n\geq 1$, define
\begin{equation}\label{eqn:def-P-plus}
P^n_-(s):=
F_-\bigl((-1)^{n-1}s\bigr)\cdots F_-(-s)F_-(s),
\qquad
P^0_-(s):=\id_{\CC^2}.
\end{equation}

\begin{cor}  \label{cor:closedformspowersf_+}
Let $a_0, b_0 \in \widetilde{Q}_2$ and denote $(a_n,b_n):=f^n_-(a_0,b_0)$. Then, we have
\begin{equation}\label{eq:f_-^nmatrixexpression}
\begin{pmatrix}
a_n\\ b_n
\end{pmatrix}
=
P^n_-(s_0)
\begin{pmatrix}
a_0\\ b_0
\end{pmatrix},
\end{equation}
where $s_0 = \langle a_0, b_0\rangle$.
\end{cor}
\begin{rem}
Although $f_+=f_-^{-1}$, the matrices $F_+(s)$ and $F_-(s)$ with the same parameter are not, in general, inverses of each other, as the maps are not linear. The correct identities are
\begin{equation}
F_+(s)^{-1}=F_-(-s),
\qquad
F_-(s)^{-1}=F_+(-s),
\end{equation}
which reflect that the scalar $\langle a,b\rangle$ changes sign after each application of either crossing map. 
\end{rem}

Now, in order to compute $T_2^n: (-,-) \to (-,-)$ for $n \geq 1$, observe that it can be decomposed as $\left(R^2(X_+)\right)^n$, where $R^2(X_+)$ is the result of applying twice the rotation of $X_+$, i.e.\ flipping the tangle. Since $Z_{\widetilde{Q}_2}(X_+)$ is the span
$$
    \widetilde{Q}_2 \times \widetilde{Q}_2 \xleftarrow{\id \times \id}  \widetilde{Q}_2 \times \widetilde{Q}_2 \xrightarrow{(f_+^1, f_+^2)} \widetilde{Q}_2 \times \widetilde{Q}_2,
$$
with $f_+^1(a,b) = b+\langle a,b\rangle a$ and $f_+^2(a,b)=a$ then, applying Lemma \ref{lem:rotation} twice, we obtain that $Z_{\widetilde{Q}_2}(R^2(X_+))$ is
$$
    \widetilde{Q}_2 \times \widetilde{Q}_2 \xleftarrow{(f_+^2, f_+^1)} \widetilde{Q}_2 \times \widetilde{Q}_2 \xrightarrow{(a,b) \mapsto (b,a)} \widetilde{Q}_2 \times \widetilde{Q}_2.
$$
In other words, denoting by $C: \widetilde{Q}_2 \times \widetilde{Q}_2 \to \widetilde{Q}_2 \times \widetilde{Q}_2$ the map $C(a,b) = (b,a)$, we have
$$
    \widetilde{Q}_2 \times \widetilde{Q}_2 \xleftarrow{C \circ f_+} \widetilde{Q}_2 \times \widetilde{Q}_2  \xrightarrow{C}\widetilde{Q}_2 \times \widetilde{Q}_2.
$$
This span is equivalent to the basic span
$$
    \widetilde{Q}_2 \times \widetilde{Q}_2 \xleftarrow{\id}\widetilde{Q}_2 \times \widetilde{Q}_2  \xrightarrow{C \circ f_+^{-1} \circ C^{-1}} \widetilde{Q}_2 \times \widetilde{Q}_2,
$$
which, using that $f_+^{-1} = f_-$ and $C^{-1}= C$, we have that $Z_{\widetilde{Q}_2}(R^2(X_+)) = S_{C \circ f_- \circ C}$ is the basic span associated to the morphism
$$
C \circ f_- \circ C: \widetilde{Q}_2 \times \widetilde{Q}_2 \to \widetilde{Q}_2 \times \widetilde{Q}_2, \qquad (a,b) \mapsto (b - \langle a, b \rangle a, a).
$$

In other words, we have that $Z_{\widetilde{Q}_2}(T_2^{n}) = S_{C \circ f_-^n \circ C}$ for $n \geq 1$. To identify this morphism, let $a_0, b_0 \in \widetilde{Q}_2$ and set $(a_n, b_n) = C \circ f_-^n \circ C(a_0, b_0)$ and $s_n = \langle a_n, b_n \rangle$. Arguing as in Proposition \ref{prop:closedformspowersf_-}, we have
$$
    s_{n+1} = \langle a_{n+1}, b_{n+1} \rangle = \langle b_n - \langle a_n, b_n \rangle a_n, a_n\rangle = \langle b_n, a_n \rangle = -s_n.
$$
Thus, again $s_n = (-1)^n s_0$ and
\begin{equation*}
\begin{pmatrix}
a_n\\ b_n
\end{pmatrix}
=\begin{pmatrix}
b_{n-1} - s_{n-1} a_{n-1} \\ a_{n-1}
\end{pmatrix} =
F_+(-s_{n-1})
\begin{pmatrix}
a_{n-1}\\ b_{n-1}
\end{pmatrix} = P_+^n(-s_0)
\begin{pmatrix}
a_{0}\\ b_{0}
\end{pmatrix}.
\end{equation*}

Analogously for $T_2^{-n}$ for $n \geq 1$, we have that $T_2^{-n} = (R^2(X_-))^n$. Repeating the calculation in this context, it is straightforward to check that $Z_{\widetilde{Q}_2}(R^2(X_-))$ is the basic span
$$
    \widetilde{Q}_2 \times \widetilde{Q}_2 \xleftarrow{\id}\widetilde{Q}_2 \times \widetilde{Q}_2  \xrightarrow{C \circ f_+ \circ C}\widetilde{Q}_2 \times \widetilde{Q}_2.
$$
Therefore, $Z_{\widetilde{Q}_2}(T_2^{-n}) = S_{C \circ f_+^n \circ C}$, the basic span associated to $C \circ f_+^n \circ C$. Setting again $s_0 = \langle a_0, b_0 \rangle $ and $(a_n, b_n) = C \circ f_+^n \circ C(a_0, b_0)$ for $a_0, b_0 \in \widetilde{Q}_2$, we have
\begin{equation*}
\begin{pmatrix}
a_n\\ b_n
\end{pmatrix}
=
P_-^n(-s_0)
\begin{pmatrix}
a_0\\ b_0
\end{pmatrix}.
\end{equation*}

\subsection{Alternating braidings}
\label{sec:alternating-braid}

In order to compute the alternating braidings, we will need to consider rotated tangles of the basic crossings $X_+$ and $X_-$. To do so, a direct application of Lemma \ref{lem:rotation} leads to the following result.

\begin{cor} \label{cor:rotatedbasicspans}
Let $C: \widetilde{Q}_2 \times \widetilde{Q}_2 \to \widetilde{Q}_2 \times \widetilde{Q}_2$ denote the transposition map $C(a,b)=(b,a)$ for $a, b \in \widetilde{Q}_2$. The spans of the elementary tangles $R(X_-)$ and $R^3(X_-)$ are, respectively, the basic spans associated to the maps
$$
r_1:=C\circ f_-\circ C
\qquad\text{and}\qquad
r_3:=f_-^{-1}=f_+,
$$
which equal
\begin{equation*}
 r_1(a,b)=\bigl(b -\langle a,b \rangle a, a\bigr),
 \qquad
 r_3(a,b)=\bigl(b +\langle a,b \rangle a, a\bigr).
\end{equation*}
\end{cor}
\begin{proof}
Writing $f_-(a,b) = (b, f_-^2(a,b))$, with $f_-^2(a,b) = a+\langle a,b\rangle b$, a direct calculation using Lemma \ref{lem:rotation} shows that $Z_{\widetilde{Q}_2}(R(X_-))$ is the span
$$
    \widetilde{Q}_2 \times \widetilde{Q}_2 \xleftarrow{(a,b) \mapsto (b,a)} \widetilde{Q}_2 \times \widetilde{Q}_2 \xrightarrow{(a,b) \mapsto (f_-^2(a,b),b)} \widetilde{Q}_2 \times \widetilde{Q}_2.
$$
The left-most map is exactly $C$, whereas the right most is $C \circ f_-$. Hence, using that $C^{-1} = C$, this span is equivalent to the basic span
$$
    \widetilde{Q}_2 \times \widetilde{Q}_2 \xleftarrow{\id}\widetilde{Q}_2 \times \widetilde{Q}_2 \xrightarrow{C \circ f_- \circ C}\widetilde{Q}_2 \times \widetilde{Q}_2.
$$

On the other hand, for $R^3(X_-)$, applying Lemma \ref{lem:rotation} three times we have that $Z_{\widetilde{Q}_2}(R^3(X_-))$ is the span
$$
    \widetilde{Q}_2 \times \widetilde{Q}_2 \xleftarrow{(a,b) \mapsto (b,f_-^2(a,b))} \widetilde{Q}_2 \times \widetilde{Q}_2 \xrightarrow{\id} \widetilde{Q}_2 \times \widetilde{Q}_2.
$$
The left-most map is exactly $f_-$, and thus the equivalent basic span is
$$
    \widetilde{Q}_2 \times \widetilde{Q}_2 \xleftarrow{\id} \widetilde{Q}_2 \times \widetilde{Q}_2 \xrightarrow{f_-^{-1} = f_+} \widetilde{Q}_2 \times \widetilde{Q}_2,
$$
as we wanted to prove.
\end{proof}

Now, for the alternating braidings, for $n\geq 0$, we define the alternating iterates recursively by
\begin{equation}\label{eq:def-hn-alternating}
h^0=\Id_{\widetilde{Q}_2\times \widetilde{Q}_2},
\qquad
h^{n+1}=
\begin{cases}
 r_3\circ h^n, & n \text{ even},\\
 r_1\circ h^n, & n \text{ odd}.
\end{cases}
\end{equation}

Now, observe that we can decompose the alternating braiding $T_3^n$ with $n \geq 1$ as $T_3^n = \ldots \circ R^3(X_-) \circ R(X_-) \circ R^3(X_-)$. This directly implies that the image of $T_3^n$ can be identified with the basic span
$$
    Z_{\widetilde{Q}_2}(T_3^n) = S_{h^n}.
$$

\begin{prop} \label{prop:alternatingtanglepowerformula}
Let $a_0, b_0 \in \widetilde{Q}_2$ and denote $s_0:= \langle a_0,b_0\rangle$. Then
\begin{equation} \label{eqn:f_+^nmatrixexpression}
h^n\begin{pmatrix} a_0 \\ b_0 \end{pmatrix}=F_+(s_0)^n \begin{pmatrix} a_0 \\ b_0\end{pmatrix}, \qquad \textrm{where } F_+(s):=
\begin{pmatrix}
 s & 1\\
 1 & 0
\end{pmatrix}.
\end{equation}
\end{prop}

\begin{proof}
Let $(a_n, b_n) := h^n(a_0, b_0)$ and set $s_n:= \langle a_n,b_n\rangle$. From Corollary \ref{cor:rotatedbasicspans}, since $h^{n+1} = r_3 \circ h^{n}$ or $h^{n+1} = r_1 \circ h^{n}$, depending on the parity of $n$, we get
$$
    \langle a_{n+1}, b_{n+1}\rangle = \left\langle b_n \pm \langle a_n, b_n\rangle a_n, a_n \right\rangle = -\langle a_{n},b_{n}\rangle.
$$
Therefore, $s_n=(-1)^ns_0$ for all $n \geq 0$. If $n$ is even, then $s_n=s_0$, and the recursive definition \eqref{eq:def-hn-alternating} gives
$$
(a_{n+1},b_{n+1})=r_3(a_{n},b_{n})=(b_n + s_0 a_n, a_n).
$$
But also, if $n$ is odd, $s_n=-s_0$, so
$$
(a_{n+1},b_{n+1})=r_1(a_{n},b_{n})=(b_n - s_n a_n, a_n)=(b_n + s_0 a_n, a_n).
$$
This implies that we get the same linear recurrence at every step, namely,
$$
\begin{pmatrix}
 a_{n+1}\\ b_{n+1}
\end{pmatrix}
=
F_+(s_0)
\begin{pmatrix}
 a_n\\ b_n
\end{pmatrix},
$$
as we wanted to prove.
\end{proof}

One may also consider the alternating composition that starts with \(r_1\)
rather than \(r_3\), which corresponds to the tangle $T_4^n$ with $n \geq 1$, which is the same alternating braiding with opposite orientation on each strand. In other words, $T_4^n = \ldots \circ R(X_-) \circ R^3(X_-) \circ R(X_-)$.

In this case, we define
\[
\widehat{h}^0=\Id_{\widetilde{Q}_2\times \widetilde{Q}_2},
\qquad
\widehat{h}^{n+1}=
\begin{cases}
r_1\circ\widehat{h}^n, & n \text{ even},\\[2pt]
r_3\circ\widehat{h}^n, & n \text{ odd}.
\end{cases}
\]
Then the same argument as in Proposition \ref{prop:alternatingtanglepowerformula} shows that, writing $s_0:=\langle a_0,b_0\rangle$,
$$
\widehat{h}^n\begin{pmatrix} a_0 \\ b_0 \end{pmatrix}
=
F_+(-s_0)^n
\begin{pmatrix}
a_0\\ b_0
\end{pmatrix}.
$$
In other words, reversing the alternating tangle amounts to replacing \(s_0\) by \(-s_0\) in the matrix
formula of Proposition~\ref{prop:alternatingtanglepowerformula}. In this manner, the image of $T_4^n$ can be identified with the basic span
$$
    Z_{\widetilde{Q}_2}(T_4^n) = S_{\widehat{h}^n}.
$$

\begin{rem}
Similarly, in the case of the other alternating braiding, the basic spans associated with the elementary tangles $R(X_+)$ and $R^3(X_+)$ are determined by maps
$$
r_1^+:=C \circ f_+ \circ C \qquad \text{and} \qquad r_3^+:= f_+^{-1} = f_-.
$$
If again $s_0= \langle a_0,b_0 \rangle$, then a similar calculation as above shows that
$$ h^n_+ \begin{pmatrix} a_0 \\ b_0 \end{pmatrix} = F_-(s_0)^n \begin{pmatrix} a_0\\ b_0 \end{pmatrix}, \qquad \widehat{h}^n_+ \begin{pmatrix} a_0 \\ b_0 \end{pmatrix} = F_-(-s_0)^n \begin{pmatrix} a_0\\ b_0 \end{pmatrix}
$$
where $h^n_+$ corresponds to the iterates where the iteration starts with $r_3^+$, and $\widehat{h}^n_+$ where it starts with $r^+_1$. Notice that they correspond respectively to $T_3^{-n}$ and $T_4^{-n}$ for $n \geq 1$.
\end{rem}

The summary of the results of this section is displayed in Table \ref{tab:summary-tangles}.

\begin{table}[H]
    \centering
    \begin{tabular}{|c|c|c|c|c|c|c|c|c|}
        \hline \textup{\textbf{Tangle}} & $T_1^n$ & $T_1^{-n}$ & $T_2^n$ & $T_2^{-n}$ & $T_3^n$ & $T_3^{-n}$ & $T_4^n$ & $T_4^{-n}$\\\hline
        \textup{\textbf{Morphism}} & $f_+^n$ & $f_-^n$ & $C \circ f_-^n \circ C$ & $C \circ f_+^n \circ C$  & $h^n$ & $h^n_+$ & $\widehat{h}^n$ & $\widehat{h}^n_+$ \\\hline
        \textup{\textbf{Matrix}} & $P^n_+(s_0)$ & $P^n_-(s_0)$ & $P_+^n(-s_0)$ & $P_-^n(-s_0)$ & $F_+(s_0)^n$ & $F_-(s_0)^n$ & $F_+(-s_0)^n$ & $F_-(-s_0)^n$ \\\hline
    \end{tabular}
    \caption{Morphisms and matrix representations of each of the tangles used to construct a pretzel knot. Here, we are considering $n \geq 1$.}
    \label{tab:summary-tangles}
\end{table}

\subsection{Closing the pretzel braiding}\label{sec:closing-pretzel}

With the descriptions obtained in the previous section, we are ready to explicitly compute the parabolic representation variety of a pretzel knot. 

\begin{figure}[h]
\includegraphics[width=12cm]{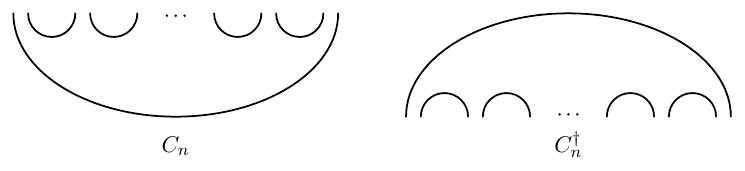}
\caption{The pretzel closure tangles $C_n$ and $C_n^{\dagger}$} \label{fig:pretzelclosuretangles}
\end{figure}

As previously explained, an oriented pretzel knot $P(q_1,\ldots,q_n)$, for $q_1, \ldots, q_n \in \ZZ$, is made of a juxtaposition of tangles $T_{i_1}^{q_1}, \ldots, T_{i_n}^{q_n}$, where $1 \leq i_j \leq 4$ are types of braidings uniquely determined by $q_1, \ldots, q_n$. More explicitly, $P(q_1,\ldots,q_n)$ can be decomposed as a tangle as
\begin{equation}\label{eq:decomposition-pretzel}
    P(q_1,\ldots,q_n) = C_n^\dagger \circ (T_{i_1}^{q_1} \otimes  \ldots \otimes  T_{i_n}^{q_n}) \circ C_n,
\end{equation}
where $C_n: \emptyset \to (\epsilon_1, \ldots, \epsilon_{2n})$ is the closure tangle and $C_n^\dagger: (\epsilon_1, \ldots, \epsilon_{2n}) \to \emptyset$ is the reversed closure, as in Figure \ref{fig:pretzelclosuretangles}.

Using the description of Section \ref{sec:TQFT-relative} and \ref{sec:explicit-functor}, it is straightforward to notice that $Z_{\widetilde{Q}_2}(C_n)$ and $Z_{\widetilde{Q}_2}(C_n^\dagger)$ are given, respectively, by the spans
\[
\begin{tikzcd}[
  column sep={5em,between origins},
  row sep={3.5em,between origins}
]
& \widetilde{Q}_2^n
  \arrow[dl]
  \arrow[dr,"\iota_n"] & \\
\star
& &
\left(\widetilde{Q}_2\times\widetilde{Q}_2\right)^n
\end{tikzcd}
\qquad\qquad
\begin{tikzcd}[
  column sep={5em,between origins},
  row sep={3.5em,between origins}
]
& \widetilde{Q}_2^n
  \arrow[dl,"\iota_n"']
  \arrow[dr] & \\
\left(\widetilde{Q}_2\times\widetilde{Q}_2\right)^n
& &
\star
\end{tikzcd}
\]
where $\iota_n$ is the morphism
$$
\iota_n\colon \widetilde{Q}_2^n\longrightarrow \left(\widetilde{Q}_2 \times \widetilde{Q}_2\right)^n,
\qquad
\iota_n(a_1,\ldots,a_{n})
=
\bigl((a_1,a_2),(a_2,a_3),\ldots,(a_{n},a_1)\bigr).
$$

\begin{thm}\label{thm:pullback-rep-var}
Let $P(q_1,\ldots,q_n) \subseteq S^3$ be an oriented pretzel knot made of braidings $T_{i_1}^{q_1}, \ldots, T_{i_n}^{q_n}$. Then, the $\widetilde{Q}_2$-representation variety of $P(q_1,\ldots,q_n)$ is given as the pullback
$$
\xymatrix{
\R_{\widetilde{Q}_2}(P(q_1,\ldots,q_n)) \ar[r]\ar[d] & \widetilde{Q}_2^n\ar[d]^{(f_{T_{i_1}^{q_1}} \times \ldots \times f_{T_{i_n}^{q_n}}) \circ \iota_n} \\
\widetilde{Q}_2^n \ar[r]_{\iota_n} & \left(\widetilde{Q}_2 \times \widetilde{Q}_2\right)^n
}
$$
where $f_{T_{i_1}^{q_1}}, \ldots, f_{T_{i_n}^{q_n}}: \widetilde{Q}_2 \times \widetilde{Q}_2 \to \widetilde{Q}_2 \times \widetilde{Q}_2$ are the morphisms corresponding to $T_{i_1}^{q_1}, \ldots, T_{i_n}^{q_n}$ respectively, as computed in Table \ref{tab:summary-tangles}.

\end{thm}
\begin{proof}
Notice that, seen as a tangle $P(q_1,\ldots,q_n): \emptyset \to \emptyset$, we have that $Z_{\widetilde{Q}_2}(P(q_1,\ldots,q_n))$ is the span
$$
    \star \longleftarrow \R_{\widetilde{Q}_2}(P(q_1,\ldots,q_n)) \longrightarrow \star
$$

Now, using the decomposition (\ref{eq:decomposition-pretzel}) and applying the functor $Z_{\widetilde{Q}_2}$, we get that $Z_{\widetilde{Q}_2}(P(q_1,\ldots,q_n))$ can also be computed as the composition of the spans
\begin{center}
\begin{tikzcd}[
  column sep={7.8em,between origins},
  row sep={4.4em,between origins},
  cells={nodes={inner sep=1.2pt}},
  every label/.append style={font=\scriptsize}
]
\widetilde{Q}_2^n
  \arrow[d]
  \arrow[dr,"\iota_n",pos=.36]
& &
\left(\widetilde{Q}_2\times\widetilde{Q}_2\right)^n
  \arrow[dl,"\Id"',pos=.36]
  \arrow[dr,
    "{f_{T_{i_1}^{q_1}}\times\cdots\times f_{T_{i_n}^{q_n}}}",
    pos=.36
  ]
& &
\widetilde{Q}_2^n
  \arrow[dl,"\iota_n"',pos=.36]
  \arrow[d]
\\
|[xshift=-3.8em]|\star
&
\left(\widetilde{Q}_2\times\widetilde{Q}_2\right)^n
& &
\left(\widetilde{Q}_2\times\widetilde{Q}_2\right)^n
&
|[xshift=3.8em]|\star
\end{tikzcd}
\end{center}

For the left-most pullback, since one of the morphisms is the identity, we get that the pullback is
\begin{center}
\begin{tikzcd}[
  column sep={7.0em,between origins},
  row sep={4.4em,between origins},
  cells={nodes={inner sep=1.2pt}},
  every label/.append style={font=\scriptsize}
]
& &
\widetilde{Q}_2^n
  \arrow[dl,"\Id"',pos=.45]
  \arrow[dr,"\iota_n",pos=.45]
& & \\
&
\widetilde{Q}_2^n
  \arrow[d]
  \arrow[dr,"\iota_n",pos=.36]
& &
\left(\widetilde{Q}_2\times\widetilde{Q}_2\right)^n
  \arrow[dl,"\Id"',pos=.36]
  \arrow[dr,
    "{f_{T_{i_1}^{q_1}}\times\cdots\times f_{T_{i_n}^{q_n}}}",
    pos=.36
  ]
& \\
&
|[xshift=-3.2em]|\star
&
\left(\widetilde{Q}_2\times\widetilde{Q}_2\right)^n
& &
\left(\widetilde{Q}_2\times\widetilde{Q}_2\right)^n
\end{tikzcd}
\end{center}
From this, we get that the final composition equals the desired pullback diagram.
\end{proof}

\subsubsection{Semiparabolic character varieties for pretzel knots}\label{sec:Q2tilde-pretzel}

Using the previous results, it is possible to explicitly describe the quandle representation variety and quandle character variety for $\widetilde{Q}_2$ and $Q_2$, for any pretzel knot. Throughout this section, we will use the notation and results of Sections \ref{sec:determinant-free} and \ref{sec:determinant-knot}.

Let us begin with the $\widetilde{Q}_2$-quandle representation variety. Let $K = P(q_1, \ldots, q_n) \subseteq S^3$ be a pretzel knot with associated morphism $f := f_{T_{i_1}^{q_1}} \times \ldots \times f_{T_{i_n}^{q_n}}: \widetilde{Q}_2^{2n} \to \widetilde{Q}_2^{2n}$. From the results of Sections \ref{sec:codirected-braid} and \ref{sec:alternating-braid}, the map $f$ can be represented by a $2n \times 2n$ matrix $M_{s_1, \ldots, s_n}$ with entries in $\ZZ[s_{1}, \ldots, s_n]$, which is made of $n$ matrices of order $2 \times 2$ as given by Table \ref{tab:summary-tangles}. In this manner, we have
$$
    f \circ \iota_n(a_1, \ldots, a_n) = M_{s_{1}, \ldots, s_{n}}\,\iota_n(a_1, \ldots, a_n) = M_{s_{1}, \ldots, s_{n}}(a_1, a_2, a_2, \ldots, a_n, a_n, a_1)^t,
$$
where $s_1 = \langle a_1, a_2\rangle, s_2 = \langle a_2, a_3\rangle, \ldots, s_n = \langle a_n, a_1\rangle$, and $x^t$ denotes the transpose of the vector $x$. In other words, considering the $2n \times n$ matrix
$$
    D_n = \left(\begin{array}{ccccc}
1 & 0 & 0 & \cdots & 0 \\
0 & 1 & 0 & \cdots  & 0 \\
0 & 1 & 0 & \cdots  & 0 \\
\vdots &  & \ddots & & \vdots  \\
\vdots &  &  & \ddots & \vdots \\
0 & 0 & 0 & \cdots & 1 \\
0 & 0 & 0 & \cdots & 1 \\
1 & 0 & 0 & \cdots & 0
\end{array}\right)
$$
we have
$$
 f \circ \iota_n(a_1, \ldots, a_n) = M_{s_{1}, \ldots, s_{n}}D_n(a_1, \ldots, a_n)^t.
$$

Notice that the vector $r = M_{s_{1}, \ldots, s_{n}}D_n(a_1, \ldots, a_n)^t$ has $2n$ components with polynomial entries in the variables $s_1, \ldots, s_n$ and linear in $a_1, \ldots, a_n$. Let us denote its components by $r = (r_1, \ldots, r_{2n})$, with $r_p \in \CC^2$ for $1 \leq p \leq 2n$. Setting $\Delta = \{(2,3), (4,5), \ldots, (2n-2, 2n-1), (1,2n)\}$, the pullback by $\iota_n: \widetilde{Q}_2^n \to \widetilde{Q}_2^{2n}$ is given by imposing $r_p = r_q$ for every $(p,q) \in \Delta$. Therefore, applying Theorem \ref{thm:pullback-rep-var}, we get the explicit description of the representation variety as follows.

\begin{prop}
For any pretzel knot $P(q_1, \ldots, q_n)$, we have
\begin{align*}
    \R_{\widetilde{Q}_2}(P(q_1, \ldots, q_n)) &= \left\{(a_1, \ldots, a_n) \in (\CC^{2})^n \,\mid\, r_p = r_q \textrm{ for }(p,q) \in \Delta\right\} \\
    & = \Spec \frac{\CC[a_{1}^1, a_{1}^2, a_{2}^1, \ldots a_{n}^1, a_{n}^2]}{(r_p - r_q)_{(p,q) \in \Delta}}.
\end{align*}
\end{prop}
Here, we are writing each $a_k \in \CC^2$ as $a_k = (a_{k}^1, a_{k}^2)$ for $1 \leq k \leq n$, and the relations $r_p - r_q$ must be written as polynomials in $a_{1}^1, a_1^2, a_2^1, \ldots, a_n^2$ by rewriting the variables $s_k$ appearing in $r$ in terms of the corresponding determinants, and considering the elements of $\CC[a_{1}^1, a_{1}^2, a_{2}^1, \ldots a_{n}^1, a_{n}^2]$ given by each of their two components. In other words, in the notation of Sections \ref{sec:determinant-free} and \ref{sec:determinant-knot}, we have
$$
    \widetilde{A}_K = \frac{\CC[a_{1}^1, a_{1}^2, a_{2}^1, \ldots a_{n}^1, a_{n}^2]}{\widetilde{I}_K}, \qquad \widetilde{I}_K = (r_p - r_q)_{(p,q) \in \Delta}.
$$
Notice that $\widetilde{A}_n = \CC[a_{1}^1, a_{1}^2, a_{2}^1, \ldots a_{n}^1, a_{n}^2]$ and the associated map $\widetilde{A}_n \to \widetilde{A}_K = \widetilde{A}_n/\widetilde{I}_K$ is the quotient morphism.

From this expression, we can also compute the quandle character variety $\mathcal{X}_{\widetilde{Q}_2}(P(q_1, \ldots, q_n)) = \R_{\widetilde{Q}_2}(P(q_1, \ldots, q_n)) \sslash \SL_2$, where the action of $\SL_2$ is $P \cdot a_k = Pa_k$ for $P \in \SL_2$ and $1 \leq k \leq n$. In order to do so, first observe that $\mathcal{X}_{\widetilde{Q}_2}(P(q_1, \ldots, q_n)) \subseteq \mathcal{X}_{\widetilde{Q}_2}(\Free{n})$. The description of the latter variety is given by Proposition \ref{prop:tildeQ2-SL2} to be $\X_{\widetilde{Q}_2}(\Free{n}) = \Spec \widetilde{B}_n$, where
$$
        \widetilde{B}_n = \frac{\CC[s_{ij}]}{(s_{ij}s_{k\ell}-s_{ik}s_{j\ell}+s_{i\ell}s_{jk}, s_{ij} + s_{ji})}, \quad \textrm{for } 1 \leq i, j \leq n \textrm{ and } k,\ell \neq i,j.
$$
Here, $s_{ij} = \langle a_i, a_j\rangle$ for $1 \leq i, j \leq n$. Recall that, as in Remark \ref{rem:redundant-Bntilde}, for simplicity we have introduced new variables $s_{ij}$ with $i \geq j$ satisfying $s_{ij} = - s_{ji}$. In particular, with respect to the original variables, we have $s_1 = s_{12}, s_2 = s_{23}, \ldots, s_{n-1} = s_{(n-1)n}$ and $s_n = s_{n1} = -s_{1n}$.

To give equations for $\mathcal{X}_{\widetilde{Q}_2}(P(q_1, \ldots, q_n))$ inside $\mathcal{X}_{\widetilde{Q}_2}(\Free{n})$ let us consider the regular functions on $\R_{\widetilde{Q}_2}(P(q_1, \ldots, q_n))$ given by
$$
    f_{pk} := \langle r_p, a_k\rangle, \qquad \textrm{for } 1 \leq p \leq 2n \textrm{ and }1 \leq k \leq n.
$$
These functions are invariant with respect to the action of $\SL_2$. Indeed, recall that each component $r_p \in \ZZ[s_{12}, \ldots, s_{1n}, a_{1}, \ldots, a_n]$ is polynomial in $s_{12}, \ldots, s_{1n}$, which are invariant under the $\SL_2$-action, and has degree $1$ in $a_{1}, \ldots, a_n$. This implies that the functions $r_p$ are $\SL_2$-equivariant since, for $P \in \SL_2$, we have that $(P \cdot r_p)(s_{12}, \ldots, s_{1n}, a_{1}, \ldots, a_n) = r_p(s_{12}, \ldots, s_{1n}, P^{-1}a_{1}, \ldots, P^{-1}a_n) = P^{-1}r_p(s_{12}, \ldots, s_{1n}, a_{1}, \ldots, a_n)$. Thus, for the functions $f_{pk}$ we have
$$
    P\cdot f_{pk} = \langle P \cdot r_p, P^{-1}a_k\rangle = \langle P^{-1}r_p, P^{-1}a_k\rangle = \langle r_p, a_k\rangle = f_{pk}.
$$

In this manner, since we have $r_p - r_q = 0$ for all $(p, q) \in \Delta$, this implies that
$$
    f_{pk} - f_{qk}= \langle r_p - r_q, a_k\rangle = 0, \textrm{for } 1 \leq k \leq n \textrm{ and } (p,q) \in \Delta.
$$
It turns out that these these functions generate the ideal of invariant functions $\R_{\widetilde{Q}_2}(\Free{n})$ that vanish on $\R_{\widetilde{Q}_2}(P(q_1, \ldots, q_n))$, as the following result shows.

\begin{prop}\label{prop:invariant-I_K}
    Set $\widetilde{A}_n = \CC[a_1^1, \ldots, a_n^2]$ with its natural action by $\SL_2$. Given a pretzel knot $K = P(q_1, \ldots, q_n) \subseteq S^3$, let $\widetilde{A}_K = \widetilde{A}_n/\widetilde{I}_K$  be the ring of regular functions on $\R_{\widetilde{Q}_2}(K)$, with $\widetilde{I}_K  = (r_p-r_q)_{(p,q) \in \Delta} \subseteq \widetilde{A}_n$. Then, $\widetilde{I}_K$ is closed under the action of $\SL_2$ and the subideal of invariant elements of $\widetilde{I}_K$ is generated over $\widetilde{B}_n = \widetilde{A}_n^{\SL_2}$ by the functions $f_{pk}-f_{qk}$ for $(p,q) \in \Delta$ and $1 \leq k \leq n$.
\end{prop}

\begin{proof}
    To prove that $\widetilde{I}_K$ is closed for the $\SL_2$-action, notice that, for each generator $r_p - r_q$, by the calculation above we have
    $$
        P \cdot (r_p - r_q) = P^{-1}(r_p - r_q) = \begin{pmatrix}
            p_{22}(r_p^1 - r_q^1) - p_{12}(r_p^2 - r_q^2) \\
            -p_{21}(r_p^1 - r_q^1) + p_{11}(r_p^2 - r_q^2)
        \end{pmatrix}.
    $$
    Here, we write $r_p = (r_{p}^1, r_{p}^2)$ for each component of $r_p$ and $P =(p_{ij})$ for $i,j = 1,2$. 
    Hence, each component of $P \cdot (r_p - r_q)$ is a linear combination of the components of $r_p - r_q$, and thus also belongs to $\widetilde{I}_K$.

    For the invariant elements of $\widetilde{I}_K$, suppose that $f \in \widetilde{I}_K$ is an invariant element, and write it as
    $$
        f = \sum_{(p,q) \in \Delta} c_{pq}^1(r_{p}^1 - r_{q}^1) + \sum_{(p,q) \in \Delta} c_{pq}^2(r_{p}^2 - r_{q}^2),
    $$
    for some polynomial coefficients $c_{pq}^1, c_{pq}^2 \in \CC[a_{1}^1, \ldots, a_{n}^2]$. For simplicity, let us denote $\bm{a} = (a_{1}^1, \ldots, a_{n}^2)$ and $\bm{x} = (x_{p,q}^1, x_{p,q}^2)_{(p,q) \in \Delta}$, and consider the auxiliary ring $\CC[\bm{a}, \bm{x}]$ and the auxiliary element
    $$
        f' = \sum_{\alpha = 1,2}\sum_{(p,q) \in \Delta} c_{pq}^\alpha\,x_{pq}^\alpha \in \CC[\bm{a}, \bm{x}].
    $$
    This element has degree $1$ in the variables $x_{pq}^\alpha$. Furthermore, for the ring homomorphism $\psi: \CC[\bm{a}, \bm{x}] \to \CC[\bm{a}] = \widetilde{A}_n$ sending $a_{k}^\alpha \mapsto a_k^\alpha$ and $x_{pq}^\alpha \mapsto r_{p}^\alpha - r_{q}^\alpha$, we have $\psi(f') =f$.

    Let $d$ be greater than the total degree of $f'$ in the variables $\bm{a}$ and denote by $\CC^{\leq d}[\bm{a}]\langle\bm{x}\rangle \subseteq \CC[\bm{a}, \bm{x}]$ the complex vector space of polynomials of degree equal to or less than $d$ in the variables $\bm{a}$ and linear in $\bm{x}$. Restricting to this subspace, we get a linear map
    $$
        \psi_{d}: \CC^{\leq d}[\bm{a}]\langle\bm{x}\rangle \to \CC[\bm{a}] \quad \textrm{satisfying } \psi_d(f') = f.
    $$
    Furthermore, $\CC[\bm{a}, \bm{x}]$, and also $\CC^{\leq d}[\bm{a}]\langle\bm{x}\rangle$ by restriction, can be endowed with a natural action of $\SL_2$ on both the variables $a_{k}^\alpha$ and $x_{pq}^\alpha$. In this manner, since $r_p - r_q$ are $\SL_2$-equivariant, the morphism $\psi_d$ also becomes equivariant and thus $\ker(\psi_d) \subseteq \CC^{\leq d}[\bm{a}]\langle\bm{x}\rangle$ is a subrepresentation.
    
    Now, since $\SL_2$ is semisimple and we are working in characteristic zero, every finite-dimensional representation is totally reducible and, thus, there exists a complement subrepresentation $V \subseteq \CC^{\leq d}[\bm{a}]\langle\bm{x}\rangle$ such that $\CC^{\leq d}[\bm{a}]\langle\bm{x}\rangle = \ker(\psi_d) \oplus V$, leading to an isomorphism
    $$
        \psi_d|_{V}: V \stackrel{\cong}{\longrightarrow} \Im(\psi_d).
    $$
    In particular, since $f = \psi_d(f') \in \Im(\psi_d)$ this implies that there exists a unique $f_0 \in V \subseteq \CC^{\leq d}[\bm{a}]\langle\bm{x}\rangle$ such that $\psi(f_0) = f$. The key difference is that now $f_0$ is invariant under the $\SL_2$-action since, for any $P \in \SL_2$, we have
    $$
        \psi_d(P \cdot f_0) = P \cdot \psi_d(f_0) = P \cdot f = f.
    $$
    This implies that $P \cdot f_0 \in V$ is another preimage of $f$ but, by uniqueness, it must be $P \cdot f_0 = f_0$, as claimed.

    Now, by Proposition \ref{prop:tildeQ2-SL2}, we have that $\CC[\bm{a}, \bm{x}]^{\SL_2}$
    is generated by the invariant elements
    $$
        \langle a_{k}, a_{k'} \rangle, \qquad \langle x_{pq}, a_{k} \rangle, \qquad \langle x_{pq}, x_{p'q'} \rangle, \qquad \textrm{for } 1 \leq k,k' \leq n \textrm{ and } (p,q),(p',q') \in \Delta. 
    $$
    This means that $f_0 \in \CC^{\leq d}[\bm{a}]\langle\bm{x}\rangle^{\SL_2}$ must be written in terms of these elements. However, as $f_0$ has degree $1$ in the variables $x_{ij}^\alpha$, and the only generators with degree $1$ in these variables are $\langle x_{pq}, a_k \rangle$, we have
    $$
        f_0 = \sum_{(p,q) \in \Delta}\sum_{k = 1}^n d_{pqk}\,\langle x_{pq}, a_{k} \rangle,
    $$
    for some $d_{pqk} \in \CC[\bm{a}]^{\SL_2} = \widetilde{B}_n$. Therefore, we have
    $$
        f = \psi(f_0) = \sum_{(p,q) \in \Delta} \sum_{k=1}^n d_{pqk}\,\langle r_p-r_q, a_{k} \rangle = \sum_{(p,q) \in \Delta} \sum_{k=1}^n d_{pqk}(f_{pk}-f_{qk}),
    $$
    which clearly belongs to the ideal generated over $\widetilde{B}_n = \CC[\bm{a}]^{\SL_2}$ by the functions $f_{pk}-f_{qk}$, as we wanted to prove.
\end{proof}

\begin{rem}\label{rmk:f_pk-rewrite}
The functions $f_{pk} = \langle r_p, a_k\rangle$ can be easily written in terms of the variables $s_{ij}$ to represent an element of $\widetilde{B}_n$. For this purpose, write $r_p = \sum_\ell c_{p\ell} a_\ell$ for $c_{p\ell}\in \CC[s_{ij}]$. Therefore, we have
$$
    f_{pk} = \langle r_p, a_k\rangle = \sum_\ell\langle c_{p\ell} a_\ell, a_k\rangle =  \sum_\ell c_{p\ell} s_{\ell k},
$$
where all the terms are functions in the variables $s_{ij}$.
\end{rem}

In this manner, using Propositions \ref{prop:tildeQ2-SL2} and \ref{prop:invariant-I_K}, we have that
$$
    \widetilde{B}_K =\widetilde{A}_{K}^{\SL_2} =\left(\frac{\widetilde{A}_n}{\widetilde{I}_K}\right)^{\SL_2} = \frac{\widetilde{A}_n^{\SL_2}}{\widetilde{A}_n^{\SL_2} \cap \widetilde{I}_K} = \frac{\CC[s_{ij}]}{\left(s_{ij}s_{k\ell}-s_{ik}s_{j\ell}+s_{i\ell}s_{jk}, s_{ij} + s_{ji}\right) +\widetilde{J}_K},
$$
where $\widetilde{J}_K = (f_{pk}-f_{qk})$, for $1 \leq i,j,k, \ell \leq n$ and $(p,q) \in \Delta$. Therefore, we have obtained the desired description of the character variety.

\begin{cor}\label{cor:char-var-Q2tilde-pretzel}
    For any pretzel knot $P(q_1, \ldots, q_n) \subseteq S^3$, we have
    $$
        \mathcal{X}_{\widetilde{Q}_2}(P(q_1, \ldots, q_n)) = \Spec \frac{\CC[s_{ij}]}{\left(s_{ij}s_{k\ell}-s_{ik}s_{j\ell}+s_{i\ell}s_{jk}, s_{ij} + s_{ji}, f_{pk}-f_{qk}\right)},
    $$
    for $1 \leq i,j,k, \ell \leq n$ and $(p,q) \in \Delta$.
\end{cor}

Notice in particular that the map $\widetilde{B}_n \to \widetilde{B}_K$ induced by the inclusion $\mathcal{X}_{\widetilde{Q}_2}(P(q_1, \ldots, q_n)) \hookrightarrow\mathcal{X}_{\widetilde{Q}_2}(\Free{n})$ is the quotient map by $\widetilde{J}_K$.

\subsubsection{Parabolic character varieties for pretzel knots}\label{sec:Q2-pretzel}

Now, let us consider the case of $Q_2 = \widetilde{Q}_2 / \ZZ_2 = \CC^2/\ZZ_2$. In this setting, the calculation works verbatim as in the case of $\widetilde{Q}_2$ described in Section \ref{sec:Q2tilde-pretzel}, except that the equalities $r_i = r_j$ are no longer on the nose, but only up to sign. However, there is a simple way of equating two vectors up to sign, namely, given $x = (x^1, x^2) , y = (y^1, y^2) \in \CC^2$, we have $[x] = [y]$ if and only if the following equality of symmetric $2 \times 2$ matrices holds
$$
    xx^t = yy^t \Longleftrightarrow \begin{pmatrix}
        (x^1)^2 & x^1x^2 \\
        x^1x^2 & (x^2)^2
    \end{pmatrix} = \begin{pmatrix}
        (y^1)^2 & y^1y^2 \\
        y^1y^2 & (y^2)^2
    \end{pmatrix}.
$$
Here, $x$ and $y$ are seen as column vectors and $x^t$ and $y^t$ denote their transposed row vectors. In this way, for any $(p,q) \in \Delta$, we have $[r_p] = [r_q]$ if and only if
\begin{equation}\label{eq:equality-sign}
    r_pr_p^t = r_qr_q^t.
\end{equation}
Using the fact that each vector $r_p$ comes from the matrices of Table \ref{tab:summary-tangles}, it is straightforward to check that each term of $r_pr_p^t$ has one of the two following forms
$$
    c\,a_ka_k^t \qquad \textrm{ or } \qquad c\,s_{k\ell}a_ka_{\ell}^t ,
$$
with $c \in \CC[s_{ij}^2]$, for some $1 \leq k, \ell \leq n$ with $k \neq \ell$. In the former case, it is immediate that the term is invariant under the $\ZZ_2$-action on $a_k$, whereas in the latter, it is also invariant since sending $a_k \mapsto -a_k$ produces
$$
    c\,s_{k\ell}a_ka_{\ell}^t \mapsto c\,(-s_{k\ell})(-a_k)a_{\ell}^t  = c\,s_{k\ell}a_ka_{\ell}^t,
$$
and similarly for $a_\ell \mapsto -a_\ell$.

\begin{rem}\label{rem:a_k-a_kt-rewrite}
We can write these terms in terms of the basic invariant functions $x_k = (a_k^1)^2, z_k= (a_k^2)^2$ and $y_k=a_k^1a_k^2$ for $a_k = (a_k^1, a_k^2)$. Indeed, denoting as in Proposition \ref{prop:Q2tilde-sign} and Remark \ref{rmk:Q2tilde-sign-matrix} by
$$
    a_ka_k^t = \begin{pmatrix}
        (a_k^1)^2 & a_k^1a_k^2 \\
        a_k^1a_k^2 & (a_k^2)^2 
    \end{pmatrix} = \begin{pmatrix}
        x_k & y_k \\
        y_k & z_k 
    \end{pmatrix}
$$
then we have
$$
s_{k\ell}a_ka_\ell^t = \begin{pmatrix}
        (a_k^1)^2a_\ell^1a_\ell^2 - (a_\ell^1)^2a_k^1a_k^2 & (a_k^1)^2(a_\ell^2)^2-a_k^1a_k^2a_\ell^1a_\ell^2 \\
        a_k^1a_k^2a_\ell^1a_\ell^2-(a_k^2)^2(a_\ell^1)^2 & (a_\ell^2)^2a_k^1a_k^2-(a_k^2)^2a_\ell^1a_\ell^2
    \end{pmatrix} = \begin{pmatrix}
        x_ky_\ell - x_\ell y_k & x_kz_\ell-y_ky_\ell \\
        y_ky_\ell -z_k x_\ell & z_\ell y_k-z_ky_\ell
    \end{pmatrix}.
$$
Additionally, notice that $s_{ij}^2 = (a_i^1a_j^2-a_i^2a_j^1)^2 = x_iz_j+z_ix_j - 2y_iy_j$.
\end{rem}

In this manner, we get the following result.

\begin{prop} \label{prop:Q2pretzelrepvariety}
For any pretzel knot $P(q_1, \ldots, q_n) \subseteq S^3$, the $Q_2$-quandle representation variety is given by
\begin{align*}
    \R_{Q_2}(P(q_1, \ldots, q_n)) & = \left\{([a_1], \ldots, [a_n]) \in (\CC^2/\ZZ_2)^n\,\mid\, r_pr_p^t = r_qr_q^t \textrm{ for } (p,q) \in \Delta\right\} \\
    & = \Spec \frac{\CC[x_k, y_k, z_k]}{(r_pr_p^t-r_qr_q^t, x_kz_k - y_k^2)},
\end{align*}
for $1 \leq k \leq n$ and $(p,q) \in \Delta$.
\end{prop}

In the previous result, $r_pr_p^t-r_qr_q^t$ must be rewritten in terms of the variables $x_k = (a_k^1)^2, z_k = (a_k^2)^2$ and $y_k =a_k^1a_k^2$ as explained above. In other words, 
$$
    A_K = \frac{\CC[x_k, y_k, z_k]}{(x_kz_k - y_k^2)+I_K}
$$
for the ideal $I_K = (r_pr_p^t-r_qr_q^t)_{(p,q) \in \Delta}$, and the morphism $A_n = {\CC[x_k, y_k, z_k]}/(x_kz_k-y_k^2) \to A_K$ is the quotient map by $I_K$.

Finally, for the $Q_2$-quandle character variety $\mathcal{X}_{Q_2}(P(q_1, \ldots, q_n)) = \cR_{Q_2}(P(q_1, \ldots, q_n)) \sslash \SL_2$, we can consider the same invariant functions as in the case of $\widetilde{Q}_2$, given by
$$
    f_{pk} = \langle r_p, a_k\rangle, \qquad \textrm{for } 1 \leq p \leq 2n \textrm{ and }1 \leq k \leq n.
$$
In this case, observe that it does not directly hold $f_{pk} = f_{qk}$ on $\R_{Q_2}(P(q_1, \ldots, q_n))$ since $r_p$ may be different from $r_q$. However, since (\ref{eq:equality-sign}) holds, we have instead
$$
    f_{pk}f_{p\ell} - f_{qk}f_{q\ell} = \langle r_p, a_k\rangle\langle r_p, a_\ell\rangle - \langle r_q, a_k\rangle\langle r_q, a_\ell\rangle = 0, \quad \textrm{for } (p,q) \in \Delta \textrm{ and } 1 \leq k, \ell \leq n.
$$

Recall from Proposition \ref{prop:Q2-SL2} that, if we denote by $\widetilde{B}_n = \CC[s_{ij}]/(s_{ij}s_{k\ell}-s_{ik}s_{j\ell}+s_{i\ell}s_{jk})$ the ring of regular functions on $\X_{\widetilde{Q}_2}(\Free{n})$, then
$$
    B_n =\widetilde{B}_n^{\ZZ_2^n} = \CC[u_{ij}, t_{ijk}]/H_n
$$
is the ring of regular functions on $\mathcal{X}_{Q_2}(\Free{n}) = \mathcal{X}_{\widetilde{Q}_2}(\Free{n}) \sslash \ZZ_2^n$, where $H_n$ is the ideal of Remark \ref{rmk:ideal-Hn}. However, notice that the functions $f_{pk}f_{p\ell} - f_{qk}f_{q\ell}$ with $k \neq \ell$ are not $\ZZ_2^n$-invariant, since $a_k \mapsto -a_k$ sends 
$$
f_{pk}f_{p\ell} - f_{qk}f_{q\ell} \mapsto -f_{pk}f_{p\ell} + f_{qk}f_{q\ell},
$$
and similarly for $a_\ell \mapsto -a_\ell$, but remain invariant for $a_{m} \mapsto -a_m$ for $m \neq k, \ell$. For this reason, in order to get an invariant function, we need to multiply it by a function $h_{k\ell} \in \widetilde{B}_n$ with the symmetry that $h_{k\ell} \mapsto -h_{k\ell}$ for $a_k \mapsto -a_k$ and $a_\ell \mapsto -a_\ell$ but remains invariant for $a_{m} \mapsto -a_m$ for $m \neq k, \ell$. It is straightforward to check that the only monomials with this property have the form
\begin{equation}\label{eq:paths}
    h_{k\ell} = c\,s_{kv_1}s_{v_1v_2} \cdots s_{v_s\ell}
\end{equation}
for $c \in \CC[s_{ij}^2]$ and $v_1, \ldots, v_s \in \{1, \ldots ,n\} - \{k,\ell\}$ different. 

\begin{rem}
There exists a combinatorial way of understanding the products of (\ref{eq:paths}). Consider the complete graph $K_n$ of $n$-vertices, i.e.\ the graph with $n$ vertices, denoted $\{1, \ldots, n\}$, and all the possible simple edges between them. A simple path $\gamma$ between $k$ and $\ell$ is thus an ordered sequence of vertices of the form $\gamma = (k, v_1, v_2, \ldots, v_s, \ell)$, with all the vertices different. Then, to the simple path $\gamma$, we associate the element of $\widetilde{B}_n$ given by $s_{\gamma} := s_{kv_1}s_{v_1v_2} \cdots s_{v_s\ell}$.
\end{rem}

In this manner, the functions to consider have the form
$$
    s_{\gamma}(f_{pk}f_{p\ell} - f_{qk}f_{q\ell}), \qquad f_{pk}^2 - f_{qk}^2
$$
for $(p,q) \in \Delta$, $1 \leq k < \ell \leq n$ and $\gamma$ a simple path between $k$ and $\ell$. By construction, these functions are invariant for the $\ZZ_2^n$-action. Now, recall from Proposition \ref{prop:Q2-SL2} that the ring $\widetilde{B}_n^{\ZZ_2^n}$ is generated by $u_{ij} = s_{ij}^2$ and $t_{ijk} = s_{ij}s_{jk}s_{ik}$. Since each function $f_{pk}$ can be easily written as a polynomial in the variables $s_{ij}$ as in Section \ref{sec:Q2tilde-pretzel}, we have that the only combinations of variables appearing in each monomial of $s_{\gamma}(f_{pk}f_{p\ell} - f_{qk}f_{q\ell})$ and $f_{pk}^2 - f_{qk}^2$ are those appearing in the definition of $u_{ij}$ and $t_{ijk}$. This implies that these invariant functions can be directly written in terms of the variables $u_{ij}$ and $t_{ijk}$ to define elements of $\widetilde{B}_n^{\ZZ_2^n}$.

From this observation, working as in Proposition \ref{prop:invariant-I_K}, it is straightforward to prove that
$$
    (f_{pk}f_{p\ell} - f_{qk}f_{q\ell}) \cap \widetilde{B}_n^{\ZZ_2^n} = \left(s_{\gamma}(f_{pk}f_{p\ell} - f_{qk}f_{q\ell}), f_{pk}^2 - f_{qk}^2\right).
$$
We therefore get the following result analogous to Corollary \ref{cor:char-var-Q2tilde-pretzel}.

\begin{cor}\label{cor:char-var-Q2-pretzel}
    For any pretzel knot $K=P(q_1, \ldots, q_n) \subseteq S^3$, we have
    $$
        \mathcal{X}_{{Q}_2}(P(q_1, \ldots, q_n)) = \Spec \frac{\CC[u_{ij}, t_{ijk}]}{H_n + \left((s_{\gamma}(f_{pk}f_{p\ell} - f_{qk}f_{q\ell}), f_{pk}^2 - f_{qk}^2\right)},
    $$
    for $1 \leq i,j,k, \ell \leq n, (p,q) \in \Delta$ and $\gamma$ a simple path between $k$ and $\ell$.
\end{cor}

In other words, Corollary \ref{cor:char-var-Q2-pretzel} tells us that $B_K = A_K^{\SL_2}$ is given by
$$
B_K = \frac{\CC[u_{ij}, t_{ijk}]}{H_n + J_K},
$$
with $J_K = ((s_{\gamma}(f_{pk}f_{p\ell} - f_{qk}f_{q\ell}), f_{pk}^2 - f_{qk}^2)$ for $1 \leq k, \ell \leq n$, $(p,q) \in \Delta$ and $\gamma$ a simple path between $k$ and $\ell$. 
In particular, the morphism $B_n = \CC[u_{ij}, t_{ijk}]/H_n \to B_K$ coming from the inclusion of the character variety $ \mathcal{X}_{{Q}_2}(P(q_1, \ldots, q_n)) \hookrightarrow  \mathcal{X}_{{Q}_2}(\Free{n})$ is given by taking the quotient by $J_K$.

\subsection{An explicit example}

In this section, we shall explicitly compute all the representation and character varieties for the quandles $\widetilde{Q}_2$ and $Q_2$ associated to a pretzel knot, as described in Sections \ref{sec:Q2tilde-pretzel} and \ref{sec:Q2-pretzel}. For this purpose, the algorithms described in these Sections have been implemented in SageMath \cite{sagemath}. A publicly available repository with this code can be found in \cite{code_repository}.

\begin{figure}[htbp]
  \centering
  \includegraphics[width=0.35\textwidth]
    {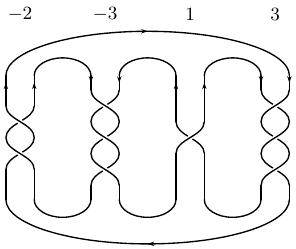}
  \caption{Pretzel knot $P(-2,-3,1,3)$.}
  \label{fig:pretzel-example}
\end{figure}

Concretely, consider the non-rational pretzel knot $K=P(-2,-3,1,3)$, as depicted in Figure \ref{fig:pretzel-example} for which all four twist regions are codirected. We orient $K$ so that the first twist region is oriented upwards and thus, with the conventions introduced above, we have
$$
    P(-2,-3,1,3) = C_4^\dagger \circ (T_{1}^{-2} \otimes T_2^{-3} \otimes T_1 \otimes   T_{2}^{3}) \circ C_4.
$$

Using Table \ref{tab:summary-tangles}, the corresponding maps for $T_{1}^{-2}, T_2^{-3}, T_1$ and $T_{2}^{3}$ are respectively
$$
\begin{array}{ll}
    P_-^2(s_1) = F_-(-s_1)F_-(s_1), & P_-^3(-s_2)=F_-(-s_2)F_-(s_2)F_-(-s_2),\\
    P^1_+(s_3) = F_+(s_3),& P_+^3(-s_4)=F_+(-s_4)F_+(s_4)F_+(-s_4),
\end{array}
$$
where $s_1 = \langle a_1, b_1\rangle$, $s_2 = \langle a_2, b_2\rangle$, $s_3 = \langle a_3, b_3\rangle$ and $s_4 = \langle a_4, b_4\rangle$ for elements $a_1, \ldots, a_4, b_1, \ldots, b_4 \in \widetilde{Q}_2$. Performing the corresponding calculations, these maps correspond respectively to the matrices
$$
\begin{pmatrix}1&s_1\\-s_1&1-s_1^2\end{pmatrix},
\quad
\begin{pmatrix}s_2&1-s_2^2\\1-s_2^2&s_2(s_2^2-2)\end{pmatrix},
\quad
\begin{pmatrix}s_3&1\\1&0\end{pmatrix},
\quad
\begin{pmatrix}
s_4(s_4^2-2)&1-s_4^2\\
1-s_4^2&s_4
\end{pmatrix}.
$$
In other words, $Z_{\widetilde{Q}_2}(T_{1}^{-2} \otimes T_2^{-3} \otimes T_1 \otimes   T_{2}^{3})$ is a basic span with associated morphism $f_{T_{1}^{-2}} \times f_{T_2^{-3}} \times f_{T_1} \times f_{T_{2}^{3}}: \widetilde{Q}_2^8 \to \widetilde{Q}_2^8$, which is given by
$$
\small
\left(\begin{array}{c}
a_1 \\ b_1 \\ a_2 \\ b_2 \\a_3 \\ b_3 \\a_4 \\ b_4
\end{array}\right) \mapsto 
\left(\begin{array}{rrrrrrrr}
1 & s_{1} & 0 & 0 & 0 & 0 & 0 & 0 \\
-s_{1} & -s_{1}^{2} + 1 & 0 & 0 & 0 & 0 & 0 & 0 \\
 0 & 0 & s_{2} & -s_{2}^{2} + 1 & 0 & 0 & 0 & 0 \\
0 & 0 & -s_{2}^{2} + 1 & {\left(s_{2}^{2} - 1\right)} s_{2} - s_{2} & 0 & 0 & 0 & 0 \\
 0 & 0 & 0 & 0 & s_{3} & 1 & 0 & 0 \\
0 & 0 & 0 & 0 & 1 & 0 & 0 & 0 \\
 0 & 0 & 0 & 0 & 0 & 0 & {\left(s_{4}^{2} - 1\right)} s_{4} - s_{4} & -s_{4}^{2} + 1 \\
0 & 0 & 0 & 0 & 0 & 0 & -s_{4}^{2} + 1 & s_{4}
\end{array}\right) \left(\begin{array}{c}
a_1 \\ b_1 \\ a_2 \\ b_2 \\a_3 \\ b_3 \\a_4 \\ b_4
\end{array}\right)
$$
Here, the previous matrix corresponds to the $8 \times 8$ matrix $M_{s_1s_2s_3s_4}$ described in Section \ref{sec:Q2tilde-pretzel}.

Now, observe that $\iota_4: \widetilde{Q}_2^4 \to \widetilde{Q}_2^8$ is the morphism
$$
    \iota_4(a_1, a_2, a_3, a_4) = (a_1, a_2, a_2, a_3, a_3, a_4, a_4, a_1).
$$
Hence, the morphism $(f_{T_{1}^{-2}} \times f_{T_2^{-3}} \times f_{T_1} \times f_{T_{2}^{3}}) \circ \iota_4: \widetilde{Q}_2^4 \to \widetilde{Q}_2^8$ is given by
\begin{align*}
\small
\left(\begin{array}{c}
a_1 \\ a_2 \\ a_3 \\ a_4
\end{array}\right) \mapsto & 
\left(\begin{array}{rrrrrrrr}
1 & s_{1} & 0 & 0 & 0 & 0 & 0 & 0 \\
-s_{1} & -s_{1}^{2} + 1 & 0 & 0 & 0 & 0 & 0 & 0 \\
 0 & 0 & s_{2} & -s_{2}^{2} + 1 & 0 & 0 & 0 & 0 \\
0 & 0 & -s_{2}^{2} + 1 & {\left(s_{2}^{2} - 1\right)} s_{2} - s_{2} & 0 & 0 & 0 & 0 \\
 0 & 0 & 0 & 0 & s_{3} & 1 & 0 & 0 \\
0 & 0 & 0 & 0 & 1 & 0 & 0 & 0 \\
 0 & 0 & 0 & 0 & 0 & 0 & {\left(s_{4}^{2} - 1\right)} s_{4} - s_{4} & -s_{4}^{2} + 1 \\
0 & 0 & 0 & 0 & 0 & 0 & -s_{4}^{2} + 1 & s_{4}
\end{array}\right) \left(\begin{array}{c}
a_1 \\ a_2 \\ a_2 \\ a_3 \\ a_3 \\ a_4 \\a_4 \\ a_1
\end{array}\right) \\
& = \left(\begin{array}{c}
a_{2} s_{1} + a_{1} \\
-{\left(s_{1}^{2} - 1\right)} a_{2} - a_{1} s_{1} \\
-{\left(s_{2}^{2} - 1\right)} a_{3} + a_{2} s_{2} \\
-{\left(s_{2}^{2} - 1\right)} a_{2} + {\left({\left(s_{2}^{2} - 1\right)} s_{2} - s_{2}\right)} a_{3} \\
a_{3} s_{3} + a_{4} \\
a_{3} \\
-{\left(s_{4}^{2} - 1\right)} a_{1} + {\left({\left(s_{4}^{2} - 1\right)} s_{4} - s_{4}\right)} a_{4} \\
-{\left(s_{4}^{2} - 1\right)} a_{4} + a_{1} s_{4}
\end{array}\right)
\end{align*}
with $s_1 = \langle a_1, a_2\rangle$, $s_2 = \langle a_2, a_3\rangle$, $s_3 = \langle a_3, a_4\rangle$ and $s_4 = \langle a_4, a_1\rangle$. In other words, as in Section \ref{sec:Q2tilde-pretzel}, considering the closing matrix
$$
    D_4 = \left(\begin{array}{ccccc}
1 & 0 & 0 & 0 \\
0 & 1 & 0 & 0 \\
0 & 1 & 0 & 0 \\
0 & 0 & 1 & 0 \\
0 & 0 & 1 & 0 \\
0 & 0 & 0 & 1 \\
0 & 0 & 0 & 1 \\
1 & 0 & 0 & 0
\end{array}\right)
$$
we have that
$$
    \left(\begin{array}{c}
r_1 \\ r_2 \\ r_3 \\ r_4 \\ r_5 \\ r_6 \\ r_7 \\ r_8
\end{array}\right) = r = M_{s_1s_2s_3s_4}D_4\left(\begin{array}{c}
a_1 \\ a_2 \\ a_3 \\ a_4
\end{array}\right) = \left(\begin{array}{c}
a_{2} s_{1} + a_{1} \\
-{\left(s_{1}^{2} - 1\right)} a_{2} - a_{1} s_{1} \\
-{\left(s_{2}^{2} - 1\right)} a_{3} + a_{2} s_{2} \\
-{\left(s_{2}^{2} - 1\right)} a_{2} + {\left({\left(s_{2}^{2} - 1\right)} s_{2} - s_{2}\right)} a_{3} \\
a_{3} s_{3} + a_{4} \\
a_{3} \\
-{\left(s_{4}^{2} - 1\right)} a_{1} + {\left({\left(s_{4}^{2} - 1\right)} s_{4} - s_{4}\right)} a_{4} \\
-{\left(s_{4}^{2} - 1\right)} a_{4} + a_{1} s_{4}
\end{array}\right).
$$

Now, taking the pullback with respect to $\iota_4$ amounts to equating the rows corresponding to the set $\Delta = \{(2,3),(4,5),(6,7),(1,8)\}$. Therefore, using Theorem \ref{thm:pullback-rep-var}, we finally obtain that
$$
    \R_{\widetilde{Q}_2}(P(-2,-3,1,3)) = \left\{(a_1, a_2, a_3, a_4) \in (\CC^2)^4 \,\left|\, \begin{array}{rcl}
a_{4} s_{4}^{2} + a_{2} s_{1} - a_{1} s_{4} + a_{1} - a_{4} & = & 0 \\
-a_{2} s_{1}^{2} + a_{3} s_{2}^{2} - a_{1} s_{1} - a_{2} s_{2} + a_{2} - a_{3}& = & 0 \\
a_{3} s_{2}^{3} - a_{2} s_{2}^{2} - 2 \, a_{3} s_{2} - a_{3} s_{3} + a_{2} - a_{4}& = & 0 \\
-a_{4} s_{4}^{3} + a_{1} s_{4}^{2} + 2 \, a_{4} s_{4} - a_{1} + a_{3}& = & 0
\end{array}\right.\right\}
$$

In other words, denoting by $a_k = (a_k^1, a_k^2)$ and $r_i = (r_i^1, r_i^2)$ the components of the vectors $a_k, r_i \in \CC^2$, we have
$$
    \widetilde{A}_4 = \CC[a_1^1, a_1^2, a_2^1, a_2^2, a_3^1, a_3^2, a_4^1, a_4^2], \qquad \widetilde{A}_K = \frac{\widetilde{A}_4}{\widetilde{I}_K},
$$
where $\widetilde{I}_K = \left(r_2^1-r_3^1, r_2^2-r_3^2, r_4^1-r_5^1, r_4^2-r_5^2, r_6^1-r_7^1, r_6^2-r_7^2, r_1^1-r_8^1, r_1^2-r_8^2\right)$, and thus the representation variety is $\R_{\widetilde{Q}_2}(P(-2,-3,1,3)) = \Spec \widetilde{A}_K$.

Furthermore, using Corollary \ref{cor:char-var-Q2tilde-pretzel}, we can also give an explicit description of the $\widetilde{Q}_2$-character variety. For this reason, we consider the functions
$$
    f_{11} = \langle r_1, a_1\rangle, \quad f_{12} = \langle r_1, a_2\rangle, \quad  \ldots\quad f_{21} = \langle r_2, a_1\rangle, \quad \ldots \quad f_{84} = \langle r_8, a_4\rangle.
$$
As described in Remark \ref{rmk:f_pk-rewrite}, these functions can be fully written in terms of the functions $s_{ij} = \langle a_i, a_j\rangle$ for $1 \leq i <j \leq 4$. For instance, we have
$$
    f_{11} = \langle r_1, a_1\rangle = \langle a_{2} s_{1} + a_{1}, a_1\rangle = s_{1}\langle a_{2}, a_1\rangle + \langle a_{1}, a_1\rangle = -s_{12}^2.
$$
Here, we are using that $s_1 = s_{12}$, $s_2 = s_{23}$, $s_3 = s_{34}$ and $s_4 = -s_{14}$.

In this way, we can compute the ideal $\widetilde{J}_K = \left(f_{pk}-f_{qk}\right)$ for $(p,q) \in \Delta = \{(2,3),(4,5),(6,7),(1,8)\}$ and $1 \leq k \leq 4$, obtaining the ideal
\begin{align*}
\widetilde{J}_K = & (-s_{14}^{3} - s_{12}^{2} + s_{14}, -s_{14}^{2} s_{24} + s_{12} s_{14} + s_{12} + s_{24}, -s_{23}^{3} - s_{12}^{2} + s_{23}, \\
&-s_{14}^{2} s_{34} + s_{13} s_{14} + s_{12} s_{23} + s_{13} + s_{34}, s_{14}^{2} + s_{12} s_{24} + s_{14}, s_{12}^{3} - s_{13} s_{23}^{2} + s_{12} s_{23} - s_{12} + s_{13}, \\
& -s_{12}^{2} s_{23} - s_{12} s_{13} - s_{23}^{2} + s_{23}, -s_{12}^{2} s_{24} + s_{23}^{2} s_{34} - s_{12} s_{14} - s_{23} s_{24} + s_{24} - s_{34}, \\
& -s_{13} s_{23}^{3} + s_{12} s_{23}^{2} + 2 s_{13} s_{23} + s_{13} s_{34} - s_{12} + s_{14}, -s_{23}^{4} + 2 s_{23}^{2} + s_{23} s_{34} + s_{24}, -s_{23}^{3} + s_{23} + s_{34},\\
&\; s_{23}^{3} s_{34} - s_{23}^{2} s_{24} - 2 s_{23} s_{34} - s_{34}^{2} + s_{24}, -s_{14}^{4} + 2 s_{14}^{2} - s_{13}, -s_{14}^{3} s_{24} + s_{12} s_{14}^{2} + 2 s_{14} s_{24} - s_{12} - s_{23},\\
&-s_{14}^{3} s_{34} + s_{13} s_{14}^{2} + 2 s_{14} s_{34} - s_{13}, s_{14}^{3} - s_{14} + s_{34}).
\end{align*}
Therefore, since the only Pl\"ucker relation for $n=4$ is $s_{14} s_{23} - s_{13} s_{24} + s_{12} s_{34} = 0$, by Corollary \ref{cor:char-var-Q2tilde-pretzel}, we have
$$
    \widetilde{B}_K = \widetilde{A}_K^{\SL_2} = \frac{\CC[s_{12},s_{23},s_{34},s_{14},s_{13},s_{24}]}{(s_{14} s_{23} - s_{13} s_{24} + s_{12} s_{34})+\widetilde{J}_K}, \qquad \X_{\widetilde{Q}_2}(K) = \Spec \widetilde{B}_K.
$$

Using a computer algebra system such as SageMath \cite{sagemath}, it is possible to compute a Gr\"obner basis for the ideal $(s_{14} s_{23} - s_{13} s_{24} + s_{12} s_{34})+\widetilde{J}_K$, leading to the alternative presentation
$$
    (s_{14} s_{23} - s_{13} s_{24} + s_{12} s_{34})+\widetilde{J}_K = \left(s_{24}^{2} + s_{24}, s_{12}, s_{13} + s_{24}, s_{14} - s_{24}, s_{23} + s_{24}, s_{34}\right).
$$
Thus, we get
$$
\widetilde{B}_K = \widetilde{A}_K^{\SL_2} = \frac{\CC[s_{24}]}{(s_{24}^{2} + s_{24})}, \qquad \X_{\widetilde{Q}_2}(K) = \{0, -1\}.
$$
In particular, $\X_{\widetilde{Q}_2}(K)$ is zero dimensional.

Regarding the $Q_2$-quandle representation variety, from Proposition \ref{prop:Q2tilde-sign} we know that
$$
    A_4 = \widetilde{A}_4^{\ZZ_2^4} = \frac{\CC[x_1, y_1, z_1, \ldots, x_4, y_4, z_4]}{(x_1z_1-y_1^2, \ldots, x_4z_4-y_4^2)}.
$$
For the ideal $I_K = (r_pr_p^t - r_qr_q^t)$ for $(p,q) \in \Delta = \{(2,3),(4,5),(6,7),(1,8)\}$, performing the corresponding products we get
\begin{align*}
I_K = & \Big({a_2} {a^t_2} s_{12}^{4} - {a_3} {a^t_3} s_{23}^{4} + {\left({a_2} {a^t_1} + {a_1} {a^t_2}\right)} s_{12}^{3} + {\left({a_3} {a^t_2} + {a_2} {a^t_3}\right)} s_{23}^{3} + {\left({a_1} {a^t_1} - 2  {a_2} {a^t_2}\right)} s_{12}^{2} \\
& \qquad - {\left({a_2} {a^t_2} - 2  {a_3} {a^t_3}\right)} s_{23}^{2} + {a_2} {a^t_2} - {a_3} {a^t_3} - {\left({a_2} {a^t_1} + {a_1} {a^t_2}\right)} s_{12} - {\left({a_3} {a^t_2} + {a_2} {a^t_3}\right)} s_{23},\\
& \quad {a_3} {a^t_3} s_{23}^{6} - {\left({a_3} {a^t_2} + {a_2} {a^t_3}\right)} s_{23}^{5} + {\left({a_2} {a^t_2} - 4  {a_3} {a^t_3}\right)} s_{23}^{4} + 3  {\left({a_3} {a^t_2} + {a_2} {a^t_3}\right)} s_{23}^{3} - {a_3} {a^t_3} s_{34}^{2} \\
& \qquad - 2  {\left({a_2} {a^t_2} - 2  {a_3} {a^t_3}\right)} s_{23}^{2} + {a_2} {a^t_2} - {a_4} {a^t_4} - 2  {\left({a_3} {a^t_2} + {a_2} {a^t_3}\right)} s_{23} - {\left({a_4} {a^t_3} + {a_3} {a^t_4}\right)} s_{34},\\
& \quad -{a_4} {a^t_4} s_{14}^{6} - {\left({a_4} {a^t_1} + {a_1} {a^t_4}\right)} s_{14}^{5} - {\left({a_1} {a^t_1} - 4  {a_4} {a^t_4}\right)} s_{14}^{4} + 3  {\left({a_4} {a^t_1} + {a_1} {a^t_4}\right)} s_{14}^{3} \\
& \qquad + 2  {\left({a_1} {a^t_1} - 2  {a_4} {a^t_4}\right)} s_{14}^{2} - {a_1} {a^t_1} + {a_3} {a^t_3} - 2 {\left({a_4} {a^t_1} + {a_1} {a^t_4}\right)} s_{14},\\
& \quad -{{\left(s_{14}^{2} - 1\right)}s_{14} {a_1a^t_4}  - {\left(s_{14}^{2} - 1\right)} {a_1a^t_1} + s_{12}{a_1a^t_2} }  + s_{12}^{2}{a_2a^t_2} + s_{12}{a_2a^t_1} \\
& \qquad - {{\left(s_{14}^{2} - 1\right)}^{2} {a_4a^t_4} - s_{14}{\left(s_{14}^{2} - 1\right)} {a_4a^t_1}} \Big).
\end{align*}

It is worth mentioning that, as explained in Remark \ref{rem:a_k-a_kt-rewrite}, every entry of $I_K$ can be written in terms of the variables $x_k = (a_k^1)^2$, $y_k = a_k^1a_k^2$ and $z_k = (a_k^2)^2$. For instance, using that $s_1 = \langle a_1, a_2\rangle = a_1^1a_2^2-a_2^1a_1^2$ and $s_1^2 = (a_1^1)^2(a_2^2)^2+(a_2^1)^2(a_1^2)^2 - 2 a_1^1a_1^2a_2^2a_2^1 = x_1z_2 + x_2z_1 - 2 y_1y_2$, we have
\begin{align*}
    r_1r_1^t & = (a_{2} s_{1} + a_{1})(a_{2} s_{1} + a_{1})^t = s_1^2a_2a_2^t + s_1a_2a_1^t + s_1 a_1a_2^t + a_1a_1^t \\
    & = s_1^2 \begin{pmatrix}
        (a_2^1)^2 & a_2^1a_2^2 \\
        a_2^1a_2^2 & (a_2^2)^2
    \end{pmatrix} + s_1 \begin{pmatrix}
        a_1^1a_2^1 & a_1^2a_2^1 \\
        a_1^1a_2^2 & a_1^2a_2^2
    \end{pmatrix} + s_1 \begin{pmatrix}
        a_1^1a_2^1 & a_1^1a_2^2 \\
        a_1^2a_2^1 & a_1^2a_2^2
    \end{pmatrix}  + \begin{pmatrix}
        (a_1^1)^2 & a_1^1a_1^2 \\
        a_1^1a_1^2 & (a_1^2)^2
    \end{pmatrix} \\
    & = s_1^2 \begin{pmatrix}
        (a_2^1)^2 & a_2^1a_2^2 \\
        a_2^1a_2^2 & (a_2^2)^2
    \end{pmatrix} +  \begin{pmatrix}
        (a_1^1)^2a_2^1a_2^2 - a_1^1a_1^2(a_2^1)^2 & a_1^1a_1^2a_2^1a_2^2-(a_1^2)^2(a_2^1)^2 \\
        (a_1^1)^2(a_2^2)^2-a_1^1a_1^2a_2^1a_2^2 & a_1^1a_1^2(a_2^2)^2-(a_1^2)^2a_2^1a_2^2
    \end{pmatrix} \\
    & \qquad + \begin{pmatrix}
        (a_1^1)^2a_2^1a_2^2 - a_1^1a_1^2(a_2^1)^2 & (a_1^1)^2(a_2^2)^2-a_1^1a_1^2a_2^1a_2^2\\
        a_1^1a_1^2a_2^1a_2^2-(a_1^2)^2(a_2^1)^2  & a_1^1a_1^2(a_2^2)^2-(a_1^2)^2a_2^1a_2^2
    \end{pmatrix} + \begin{pmatrix}
        (a_1^1)^2 & a_1^1a_1^2 \\
        a_1^1a_1^2 & (a_1^2)^2
    \end{pmatrix} \\
    & = (x_1z_2 + x_2z_1 - 2 y_1y_2)\begin{pmatrix}
        x_2 & y_2 \\
        y_2 & z_2
    \end{pmatrix} +  \begin{pmatrix}
        x_1y_2 - y_1x_2 & y_1y_2-z_1x_2 \\
        x_1z_2-y_1y_2 & y_1z_2-z_1y_2
    \end{pmatrix} \\
    & \qquad + \begin{pmatrix}
        x_1y_2 - y_1x_2 & x_1z_2-y_1y_2 \\
        y_1y_2-z_1x_2 & y_1z_2-z_1y_2
    \end{pmatrix} + \begin{pmatrix}
        x_1 & y_1 \\
        y_1 & z_1
    \end{pmatrix}.
\end{align*}
Performing this calculation for every product, we can write $I_K \subseteq A_4$ and therefore
$$
    A_K = \frac{A_4}{I_K} = \frac{\CC[x_1, y_1, z_1, \ldots, x_4, y_4, z_4]}{(x_1z_1-y_1^2, \ldots, x_4z_4-y_4^2) + I_K}, \qquad \R_{Q_2}(K) = \Spec A_K.
$$
The form of $I_K$ in these variables is very lengthy and we do not reproduce it here. However, it can be checked in the accompanying script \cite{code_repository}.

Finally, we can also obtain an explicit description for the associated parabolic character variety. For this purpose, recall that
$$
    B_4 = A_4^{\SL_2} = \widetilde{B}_4^{\ZZ_2^4} = \frac{\CC[u_{12},u_{23}, u_{34}, u_{14}, u_{13}, u_{24}, t_{123}, t_{124}, t_{134}, t_{234}]}{H_4}.
$$
The ideal $H_4$ can also be computed explicitly using the techniques from \cite{drensky2003defining} and we get
\begin{align*}
    H_4  = \bigl(
&u_{12}t_{134}-u_{13}t_{124}+u_{14}t_{123},-u_{12}t_{234}-u_{23}t_{124}+u_{24}t_{123},\\
&u_{23}t_{134}-u_{13}t_{234}+u_{34}t_{123}, -u_{14}t_{234}+u_{24}t_{134}-u_{34}t_{124},\\
&t_{123}^{2}-u_{12}u_{13}u_{23}, t_{124}^{2}-u_{12}u_{14}u_{24}, t_{134}^{2}-u_{13}u_{14}u_{34}, t_{234}^{2}-u_{23}u_{24}u_{34},\\
&2t_{123}t_{124}
 +u_{12}^{2}u_{34}
 -u_{12}u_{13}u_{24}
 -u_{12}u_{14}u_{23}, 2t_{123}t_{134}
 +u_{12}u_{13}u_{34}
 -u_{13}^{2}u_{24}
 +u_{13}u_{14}u_{23},\\
&2t_{123}t_{234}
 -u_{12}u_{23}u_{34}
 -u_{13}u_{23}u_{24}
 +u_{14}u_{23}^{2}, 2t_{124}t_{134}
 -u_{12}u_{14}u_{34}
 -u_{13}u_{14}u_{24}
 +u_{14}^{2}u_{23},\\
&2t_{124}t_{234}
 +u_{12}u_{24}u_{34}
 -u_{13}u_{24}^{2}
 +u_{14}u_{23}u_{24}, 2t_{134}t_{234}
 +u_{12}u_{34}^{2}
 -u_{13}u_{24}u_{34}
 -u_{14}u_{23}u_{34}
\bigr).
\end{align*}

In this manner, we have that the parabolic character variety is given by
$$
    B_K = \frac{B_4}{J_K} = \frac{\CC[u_{12},u_{23}, u_{34}, u_{14}, u_{13}, u_{24}, t_{123}, t_{124}, t_{134}, t_{234}]}{H_4 + J_K}, \qquad \X_{Q_2}(K) = \Spec B_K.
$$
The ideal $J_K$ is given as in Corollary \ref{cor:char-var-Q2-pretzel} by
$$
    J_K = (s_{\gamma}(f_{pk}f_{p\ell}-f_{qk}f_{q\ell}), f_{pk}^2-f_{qk}^2),
$$
for $1 \leq k < \ell \leq 4$, $(p,q) \in \Delta$ and $\gamma$ a simple path between $k$ and $\ell$. All these terms can be explicitly written in terms of the variables $u_{ij}$ and $t_{ijk}$. For instance, for $(p,q) = (1,8)$, $k = 1$ and $\ell = 4$, performing the calculation as above we straightforwardly find that $f_{14} = s_{12}s_{24} + s_{14}$, $f_{81} = (s_{14}^2-1)s_{14}$ and $f_{84} = -s_{14}^2$ and thus
$$
    f_{11}f_{14} - f_{81}f_{84} = (-s_{12}^2)(s_{12}s_{24} + s_{14}) - ((s_{14}^2-1)s_{14})(-s_{14}^2) = -s_{12}^3s_{24} -s_{12}^2s_{14} + (s_{14}^2-1)s_{14}^3.
$$

From this expression, we observe that this function satisfies $f_{11}f_{14} - f_{81}f_{84} \mapsto -(f_{11}f_{14} - f_{81}f_{84})$ for $a_1 \mapsto -a_1$ and $a_4 \mapsto -a_4$, but remains invariant for $a_2 \mapsto -a_2$ and $a_3 \mapsto -a_3$. Therefore, in order to make it invariant, we must multiply it by $s_\gamma$ for a simple path $\gamma$ joining $1$ and $4$. There are five possibilities, leading to the following five coefficients
$$
    \begin{array}{lclclcl}
        \gamma=(1\,4) & \leadsto & s_{\gamma} = s_{14}, & \qquad & \gamma=(1\,2\,4) & \leadsto & s_{\gamma} = s_{12}s_{24}, \\
        \gamma=(1\,3\,4) & \leadsto & s_{\gamma} = s_{13}s_{34}, & \qquad & \gamma=(1\,2\,3\,4) & \leadsto & s_{\gamma} = s_{12}s_{23}s_{34}, \\
        \gamma=(1\,3\,2\,4) & \leadsto & s_{\gamma} = s_{13}s_{32}s_{24} = -s_{13}s_{23}s_{24}. & & & &  \\
    \end{array}
$$
Multiplying each of these coefficients by $f_{11}f_{14} - f_{81}f_{84}$ we get a $\ZZ_2^4$-invariant function and, thus, it can be written in terms of the variables $u_{ij} = s_{ij}^2$ and $t_{ijk} = s_{ij}s_{jk}s_{ik}$. For instance, for $\gamma=(1\,3\,4)$ we have $s_{\gamma} = s_{13}s_{34}$ and therefore
\begin{align*}
    s_\gamma(f_{11}f_{14} - f_{81}f_{84}) & = s_{13}s_{34}(-s_{12}^3s_{24} -s_{12}^2s_{14} + (s_{14}^2-1)s_{14}^3) \\
    & = -s_{12}^3s_{24}s_{13}s_{34} -s_{12}^2s_{14}s_{13}s_{34} + (s_{14}^2-1)s_{14}^3s_{13}s_{34} \\
    & = -u_{12}s_{12}s_{24}s_{13}s_{34} - u_{12}t_{134} + (u_{14}-1)u_{14}t_{134}.
\end{align*}
Now, observe that $s_{12}s_{34}-s_{13}s_{24}+s_{14}s_{23}=0$ by the Pl\"ucker relation, so $s_{14}s_{23}=s_{13}s_{24}-s_{12}s_{34}$ and thus taking the square
$$
    s_{13}s_{24}s_{12}s_{34} = \frac{1}{2}(s_{13}^2s_{24}^2 + s_{12}^2s_{34}^2 - s_{14}^2s_{23}^2) = \frac{1}{2}(u_{13}u_{24} + u_{12}u_{34} - u_{14}u_{23}).
$$
Hence, finally we get 
\begin{align*}
    s_\gamma(f_{11}f_{14} - f_{81}f_{84}) = -\frac{1}{2}(u_{12}u_{13}u_{24} + u_{12}^2u_{34} - u_{12}u_{14}u_{23}) - u_{12}t_{134} + (u_{14}-1)u_{14}t_{134},
\end{align*}
as desired.

In this manner, proceeding similarly with all the terms and possible paths, we get an explicit expression of the ideal $J_K$ with $136$ generators. We do not replicate the complete description for simplicity, but the expression can be checked here. Furthermore, computing a Gr\"obner basis for $H_4 + J_K$, we get that
\begin{align*}
    H_4 + J_K =& \Big(u_{12}^{2} - u_{34}, u_{12} u_{24} - u_{12} + u_{34} + 4 t_{124} + t_{134} + 3 t_{234}, \\
    & \quad u_{24}^{2} - u_{12} - u_{24} + 4 u_{34} + 10 t_{124} + 3 t_{134} + 8 t_{234}, u_{12} u_{34} + 2 u_{12} + u_{34} - 4 t_{234},\\
    & \quad  u_{24} u_{34} + 4 u_{12} - 4 u_{34} - 12 t_{124} - 3 t_{134} - 13 t_{234}, u_{34}^{2} - 6 u_{12} - 11 u_{34} - 24 t_{124} - 8 t_{134},\\
    & \quad u_{12} t_{124} - u_{12} + u_{34} + 4 t_{124} + t_{134} + 4 t_{234}, u_{24} t_{124} - 2 u_{12} + 4 u_{34} + 11 t_{124} + 3 t_{134} + 9 t_{234}, \\
    & \quad u_{34} t_{124} + 5 u_{12} - u_{34} - 6 t_{124} - t_{134} - 12 t_{234}, t_{124}^{2} - 3 u_{12} + 4 u_{34} + 12 t_{124} + 3 t_{134} + 11 t_{234},\\
    & \quad u_{12} t_{134} + u_{12} + 7 u_{34} + 14 t_{124} + 5 t_{134} + 4 t_{234}, u_{24} t_{134} + 6 u_{12} - 4 t_{124} - 13 t_{234},\\
    & \quad u_{34} t_{134} - u_{12} - 43 u_{34} - 102 t_{124} - 31 t_{134} - 44 t_{234},\\
    & \quad t_{124} t_{134} + 6 u_{12} + 9 u_{34} + 18 t_{124} + 7 t_{134} - 3 t_{234},\\
    & \quad t_{134}^{2} + 16 u_{12} - 107 u_{34} - 272 t_{124} - 81 t_{134} - 153 t_{234},\\
    & \quad u_{12} t_{234} - u_{12} - 3 u_{34} - 6 t_{124} - 2 t_{134} - t_{234}, u_{24} t_{234} - 3 u_{12} + 2 t_{124} + 6 t_{234},\\
    & \quad u_{34} t_{234} + u_{12} + 19 u_{34} + 46 t_{124} + 14 t_{134} + 19 t_{234},\\
    & \quad t_{124} t_{234} - 3 u_{12} - 4 u_{34} - 8 t_{124} - 3 t_{134} + 2 t_{234},\\
    & \quad t_{134} t_{234} - 7 u_{12} + 50 u_{34} + 126 t_{124} + 37 t_{134} + 70 t_{234},\\
    & \quad t_{234}^{2} + 3 u_{12} - 23 u_{34} - 58 t_{124} - 17 t_{134} - 32 t_{234}, u_{13} - u_{24} + 4 u_{34} + 4 t_{124} + 4 t_{234},\\
    & \quad u_{14} - u_{24} + 3 u_{34} + 4 t_{124} + t_{134} + 4 t_{234}, u_{23} - u_{24} - u_{34} - t_{234}, t_{123} + t_{124} + t_{134} + t_{234}\Big)
\end{align*}
In particular, from this description, we can easily obtain that $\X_{Q_2}(K) = \Spec \CC[u_{ij}, t_{ijk}]/(H_4 + J_K)$ is a variety of dimension $0$ comprised of $7$ different points.

\subsection{Pretzel knots with parabolic character varieties of large dimension}

We finish this work by pointing out that the parabolic character variety of a pretzel knot need not be $0$-dimensional, as it was shown in general for Montesinos knots of Kinoshita-Terasaka type in \cite{PortiPaoluzzi2013}. However, using the techniques developed in this paper, it is very easy to show that the parabolic character variety can have arbitrarily large dimension by finding an explicit parametrization of a locally closed subvariety of arbitrary dimension on a family of pretzel knots, matching \cite{PortiPaoluzzi2013}.

\begin{prop}\label{prop:families-dimensions}
Let $d\geq 2$, $m\geq 1$ and consider the $2d$-strand pretzel knot
$$
K_{d,m}=P(\underbrace{3,\ldots,3}_{2d-1},2m).
$$
Then the parabolic character variety $\X_{Q_2}(K_{d,m})$ has an irreducible component of dimension at least $2d-3$.
\end{prop}

\begin{proof}
    It is easy to check that $K_{d,m}$ is made of $2d$ codirected braids with alternating orientations. Therefore, orienting the first braid upwards, we have
    $$
        K_{d,m} = C_{2d}^\dagger \circ (T_1^{3} \otimes T_2^{3} \otimes T_1^3 \otimes \ldots \otimes T_1^3 \otimes T_2^{2m}) \circ C_{2d}.
    $$
    Hence, $Z_{\widetilde{Q}_2}(T_1^{3} \otimes T_2^{3} \otimes T_1^3 \otimes \ldots \otimes T_1^3 \otimes T_2^{2m})$ is a basic span with associated morphism $f: \widetilde{Q}_2^{4d} \to \widetilde{Q}_2^{4d}$ given by
    $$
        f = f_{T_1^{3}} \times f_{T_2^{3}} \times \ldots \times f_{T_1^3} \times f_{T_2^{2m}}.
    $$

    Now, given $\lambda = (\lambda_2, \ldots, \lambda_{2d-2}) \in \CC^{2d-3}$, we construct the tuple $(a_1^\lambda, a_2^\lambda, \ldots, a_{2d}^\lambda) \in \widetilde{Q}_2^{2d}$ as follows. First, set $a_1^\lambda = (1,0)$ and $a_2^\lambda = (0,1)$, and now we define inductively
    $$
        a_{k+1}^\lambda = \lambda_ka_k^\lambda + a_{k-1}^\lambda, \qquad \textrm{for }2 \leq k \leq 2d-2,
    $$
    and, if $a_{2d-1}^\lambda = (a_{2d-1}^{\lambda, 1}, a_{2d-1}^{\lambda, 2})$ with $a_{2d-1}^{\lambda, 2} \neq 0$, then set $a_{2d}^\lambda = \left(-1/a_{2d-1}^{\lambda, 2}, 0\right)$.

    For this choice, if we consider the determinants $s_{i,j}^\lambda = \langle a_i^\lambda, a_j^\lambda\rangle$ for $1 \leq i < j \leq 2d$, a direct check from the definition shows that, for all $\lambda$,
    $$
        s_{k,k+1}^\lambda = (-1)^{k+1} \textrm{ for }1 \leq k \leq 2d-1, \qquad s_{1,2d}^\lambda = 0. 
    $$
    In the same vein, we have $s_{k, k+2} = (-1)^{k+1}\lambda_{k+1}$ for $1 \leq k \leq 2d-3$.
    
    In particular, by Table \ref{tab:summary-tangles}, the matrices corresponding to the functions $f_{T_1^{3}}$, $f_{T_2^{3}}$ and $f_{T_2^{2m}}$ appearing in $f$, evaluated at the chosen family, are respectively
    $$
        f_{T_1^{3}} = P_+^3(1) = \begin{pmatrix}
            1 & 0 \\
            0 & -1
        \end{pmatrix}, \qquad f_{T_2^{3}} = P_+^3(-(-1)) = \begin{pmatrix}
            1 & 0 \\
            0 & -1
        \end{pmatrix}, \qquad f_{T_2^{2m}} = P_+^{2m}(0) = \begin{pmatrix}
            1 & 0 \\
            0 & 1
        \end{pmatrix}.
    $$
    
    In other words, the matrix $M_{s_{1,2}^\lambda, \ldots, s_{1,2d}^\lambda}$ associated to $f = f_{T_1^{3}} \times f_{T_2^{3}} \times \ldots \times f_{T_1^3} \times f_{T_2^{2m}}$ is
    $$
        M_{s_{1,2}^\lambda, \ldots, s_{1,2d}^\lambda} = \textup{diag}(\underbrace{1,-1,1,-1, \ldots, 1, -1}_{4d-2}, 1, 1).
    $$
    and thus the vector $r^\lambda = M_{s_{1,2}^\lambda, \ldots, s_{1,2d}^\lambda}D_{2d}(a_1^\lambda, \ldots, a_{2d}^\lambda)^t$ is
    $$
        r^\lambda = (a_1^\lambda, -a_2^\lambda, a_2^\lambda, \ldots, -a_{2d}^\lambda, a_{2d}^\lambda, a_1^\lambda).
    $$
    In particular, the functions $f_{p,k}^\lambda = \langle r_{p}^\lambda, a_{k}^\lambda \rangle$ for $1 \leq p \leq 4d$ and $1 \leq k \leq 2d$ are
    $$
        f_{2s,k}^\lambda = -s_{s+1,k}^\lambda, \qquad f_{2s+1,k}^\lambda = s_{s+1,k}^\lambda, \qquad \textrm{for } 1 \leq s \leq 2d-1,
    $$
    and $f_{1,k}^\lambda = f_{4d,k}^\lambda = s_{1,k}^\lambda$. 
    
    Now, observe that the set of pairs to be identified is $\Delta = \{(2,3), (4,5), \ldots, (1,4d)\}$, which directly implies that
    $$
        f_{2s,k}^\lambda f_{2s,\ell}^\lambda - f_{2s+1,k}^\lambda f_{2s+1,\ell}^\lambda = (-s_{s+1,k}^\lambda)(-s_{s+1,\ell}^\lambda) - s_{s+1,k}^\lambda s_{s+1,\ell}^\lambda = 0,
    $$
    and $f_{1,k}^\lambda f_{1,\ell}^\lambda - f_{4d,k}^\lambda f_{4d,\ell}^\lambda = s_{1,k}^\lambda s_{1,\ell}^\lambda-s_{1,k}^\lambda s_{1,\ell}^\lambda = 0$. In other words, $(a_1^\lambda, a_2^\lambda, \ldots, a_{2d}^\lambda)$ satisfies the equations of $\X_{Q_2}(K_{d,m})$ for all $\lambda$. 

    Consider $B_{2d} = \CC[u_{ij}, t_{ijk}]/H_{2d}$ and consider the open set of $\CC^{2d-3}$ given by
    $$
        \Lambda = \{\lambda =(\lambda_2, \ldots, \lambda_{2d-2}) \in \CC^{2d-3}\,\mid\, a_{2d-1}^{\lambda, 2} \neq 0\}.
    $$
    This set is non-empty since, for instance, $\lambda = (0,0,\ldots,0, 1) \in \Lambda$. The argument above shows that the regular morphism
    $$
        \Lambda \stackrel{\varphi}{\longrightarrow} \X_{Q_2}(\Free{2d}) = \Spec B_{2d} \quad \lambda =(\lambda_2, \ldots, \lambda_{2d-2}) \mapsto (u_{ij}^\lambda = (s_{ij}^\lambda)^2, t_{ijk} = s_{ij}^\lambda s_{jk}^\lambda s_{ik}^\lambda)
    $$
    factorizes through the inclusion $\X_{Q_2}(K_{d,m}) = \Spec B_{2d}/J_{K_{d,k}} \hookrightarrow \X_{Q_2}(\Free{2d})$. Furthermore, the quotient map $\CC[u_{1,3}, u_{2,4}, \ldots, u_{2d-3,2d-1}] \to B_{2d}$ provides a regular morphism given by projection onto the variables $u_{1,3}, u_{2,4}, \ldots, u_{2d-3,2d-1}$,
    $$
        X_{Q_2}(\Free{2d}) \stackrel{\pi}{\longrightarrow} \mathbb{A}^{2d-3}.
    $$
    
    In this way, we have that the composition $\pi \circ \varphi: \Lambda \to \mathbb{A}^{2d-3}$ is given by $\pi \circ \varphi(\lambda_2, \ldots, \lambda_{2d-2}) = (u_{1,3}^{\lambda}, u_{2,4}^{\lambda}, \ldots, u_{2d-3,2d-1}^\lambda) = ((s_{1,3}^{\lambda})^2, (s_{2,4}^{\lambda})^2, \ldots, (s_{2d-3,2d-1}^\lambda)^2) = (\lambda_2^2, \lambda_3^2, \ldots, \lambda_{2d-2}^2)$. Therefore, $\pi \circ \varphi$ is a dominant map, and so is $\pi$. Thus, $\dim \X_{Q_2}(K_{d,m}) \geq \dim \pi(\X_{Q_2}(K_{d,m})) = 2d-3$, as we wanted to prove.
\end{proof}

\begin{rem}
    In the proof of Proposition \ref{prop:families-dimensions}, it is possible to give a closed formula for all the variables $s_{ij}^\lambda = \langle a_i^\lambda, a_j^\lambda\rangle$ for $1 \leq i < j < 2d$ as
    $$
        s_{ij}^\lambda = (-1)^{i+1}K_{j-i-1}(\lambda_{i+1}, \ldots, \lambda_{j-1}).
    $$
    Here $K_{m}(x_1, \ldots, x_m)$ is the so-called continuant polynomial, defined recursively as $K_0 = 1$, $K_1(x_1) = x_1$ and, for $m \geq 2$,
    $$
        K_m(x_1,\ldots,x_m) =x_mK_{m-1}(x_1,\ldots,x_{m-1}) + K_{m-2}(x_1,\ldots,x_{m-2}).
    $$
\end{rem}

\bibliographystyle{abbrv}
\bibliography{TQFTSL2}

\end{document}